\documentclass{article}

\input{tcilatex}
\begin{document}

\[
\text{\textbf{Some conjectured formulas for }}\frac{\mathbf{1}}{\mathbf{\pi }%
\text{ }}\text{ \textbf{coming from polytopes, K3-surfaces and Moonshine.}} 
\]%
\[
\text{\textsl{Zur Erinnerung an Max Kreuzer aus Gmunden}} 
\]%
\[
\text{Gert Almkvist} 
\]

\textbf{Calabi-Yau, Ramanujan and Guillera.}

In the papers [5], [6] Jesus Guillera and the author have established the
close connection between Calabi-Yau differential equations and formulas for $%
\dfrac{1}{\pi ^{2}}.$ We will see that the same is true for Ramanujan-like
formulas for $\dfrac{1}{\pi }.$In [10],[17],[22] the following definition is
given:

\textbf{Definition. }Let the differential operator%
\[
L=\sum_{k=0}^{n}a_{n-k}(x)D^{n-k} 
\]%
with $a_{n}=1$ and $D=\dfrac{d}{dx}$ , be given. Then%
\[
L^{\ast }=D^{n}+\sum_{k=1}^{\infty }(-1)^{k}D^{n-k}a_{n-k} 
\]%
is the \textit{adjoint }of $L.$Let \ $\alpha $ be defined by 
\[
\frac{\alpha ^{\prime }}{\alpha }=-\frac{2}{n}a_{n-1}\text{ \ \ and \ }%
\alpha (0)=1. 
\]%
Then $L$ is \textit{self-adjoint }if%
\[
L^{\ast }=\alpha ^{-1}L\alpha 
\]

\textbf{Theorem. }$L$ is self-adjoint if and only if \ $W(n,k)=0$ \ for $%
k=0,1,2,...,n-1$ , where%
\[
W(n,k)=\sum_{j=k}^{n}\dbinom{j}{k}a_{j}\frac{\alpha ^{(j-k)}}{\alpha }%
-\sum_{j=k}^{n-1}(-1)^{n-j}a_{j}^{(j-k)} 
\]%
Here we define $\dfrac{\alpha ^{(0)}}{\alpha }=1.$

\textbf{Proof. }We have%
\[
D^{k}a=\sum_{j=0}^{k}\dbinom{k}{j}a^{(j)}D^{k-j} 
\]%
which leads to%
\[
L^{\ast }=D^{n}+\sum_{k=1}^{n}(-1)^{k}\dbinom{n-k}{j}a_{n-k}^{(j)}D^{n-k-j} 
\]%
For $\dfrac{\alpha ^{(k)}}{\alpha }$ we get the recursion formula%
\[
\dfrac{\alpha ^{(k)}}{\alpha }=\frac{d}{dx}(\dfrac{\alpha ^{(k-1)}}{\alpha }%
)-\frac{2}{n}a_{n-1}\dfrac{\alpha ^{(k-1)}}{\alpha } 
\]%
for $k\geq 1$ and $\dfrac{\alpha ^{(0)}}{\alpha }=1$ by definition. Now%
\[
L\alpha =\sum_{k=0}^{n}a_{n-k}D^{n-k}\alpha
=\sum_{k=0}^{n}a_{n-k}\sum_{j=0}^{n-k}\dbinom{n-k}{j}\alpha ^{(j)}D^{n-k-j} 
\]%
and%
\[
\alpha ^{-1}L\alpha =\sum_{k=0}^{n}a_{n-k}\sum_{j=0}^{n-k}\dbinom{n-k}{j}%
\frac{\alpha ^{(j)}}{\alpha }D^{n-k-j} 
\]%
Comparing the two sides of \ $L^{\ast }=\alpha ^{-1}L\alpha $ we are done.

In general we find \ $W(n,n-1)=W(n,n-2)\equiv 0$ \ and%
\[
W(n,n-3)=2a_{n-3}-(2-\frac{4}{n})a_{n-1}a_{n-2}+\frac{4}{n^{3}}\dbinom{n}{3}%
a_{n-1}^{3}-(n-2)a_{n-2}^{\prime }+\frac{6}{n^{2}}\dbinom{n}{3}%
a_{n-1}a_{n-1}^{\prime }+\frac{1}{n}\dbinom{n}{3}a_{n-1}^{\prime \prime } 
\]

\textbf{n=3.}%
\[
W(3,0)=2a_{1}-\frac{2}{3}a_{2}a_{3}+\frac{4}{27}a_{2}^{3}-a_{1}^{\prime }+%
\frac{2}{3}a_{2}a_{2}^{\prime }+\frac{1}{3}a_{2}^{\prime \prime }=0 
\]%
This condition is exactly the condition for that the analytic solution $%
y_{0} $ should be a square, $y_{0}=u^{2}$ , where $u$ is the solution of a
second order differential equation (see [1],[7],[8]). This condition is
crucial for existence of formulas for $\dfrac{1}{\pi }.$

\textbf{n=4.}%
\[
W(4,1)=2a_{1}-a_{2}a_{3}+\frac{1}{4}a_{2}^{3}-2a_{2}^{\prime }+\frac{3}{2}%
a_{2}a_{2}^{\prime }+\frac{1}{3}a_{3}^{\prime \prime }=0 
\]%
This is the condition, first found in [7], for that a 4-th order
differential equation should be Calabi-Yau. $W(4,0)$ contains $a_{1}$ and $%
a_{1}^{\prime }$ and is identically zero when $W(4,1)=0$ is substituted.

\textbf{n=5.}

Here we have two conditions%
\[
W(5,2)=2a_{2}-\frac{6}{5}a_{3}a_{4}+\frac{8}{25}a_{4}^{3}-3a_{4}^{\prime }+%
\frac{12}{5}a_{4}a_{4}^{\prime }+2a_{4}^{\prime \prime }=0 
\]%
\[
W(5,0)=2a_{0}-\frac{2}{5}a_{1}a_{4}+\frac{4}{25}a_{2}a_{4}^{2}-\frac{8}{125}%
a_{3}a_{4}^{3}+\frac{48}{3125}a_{4}^{5}-a_{1}^{\prime }-\frac{2}{5}%
a_{2}a_{4}^{\prime }-\frac{16}{125}a_{4}^{3}a_{4}^{\prime }+\frac{12}{25}%
a_{3}a_{4}a_{4}^{\prime } 
\]%
\[
-\frac{12}{25}a_{4}a_{4}^{\prime 2}+a_{2}^{\prime \prime }-\frac{2}{5}%
a_{3}a_{4}^{\prime \prime }+\frac{8}{5}a_{4}^{\prime }a_{4}^{\prime \prime
}-a_{3}^{\prime \prime \prime }+\frac{2}{5}a_{4}a_{4}^{\prime \prime \prime
}+\frac{3}{5}a_{4}^{\prime \prime \prime \prime }=0 
\]%
Again we have $W(5,1)=0$ follows from $W(5,2)=0.$

\textbf{n=7.}%
\[
W(7,4)=2a_{4}-\frac{10}{7}a_{5}a_{6}+\frac{20}{49}a_{6}^{3}-5a_{5}^{\prime }+%
\frac{30}{7}a_{6}a_{6}^{\prime }+5a_{6}^{\prime \prime }=0 
\]%
We do not write out $W(7,2)$ and $W(7,0),$since they are rather long. Only
one formula for $\dfrac{1}{\pi ^{3}}$ is known (found by B.Gourevitch)%
\[
\sum_{n=0}^{\infty }\dbinom{2n}{n}^{7}(1+14n+76n^{2}+168n^{3})\frac{1}{%
2^{20n}}=\frac{32}{\pi ^{3}} 
\]%
The 7-th order differential equation satisfied by%
\[
y_{0}=\sum_{n=0}^{\infty }\dbinom{2n}{n}^{7}x^{n} 
\]%
satisfies the three conditions above.

\textbf{Symmetric and exterior powers of spaces of solutions to differential
equations.}

Start with a (MUM) 4-th order differential equation with solutions $%
y_{0},y_{1},y_{2},y_{3}$. Then the six wronskians $w_{0}=x(y_{0}y_{1}^{%
\prime }-y_{0}^{\prime }y_{1}),$ etc satisfy a 5-th order differential
equation if and only if \ $W(4,1)=0$ (and then the 5-th order equation
satisfies $W(5,2)=W(5,0)=0$ ). This was the start of my collaboration with
Zudilin and lead to the Calabi-Yau industry.

Another way to get a 5-th order equation is to take symmetric squares of the
solutions of a 3-rd order equation, i.e. $%
w_{0}=y_{0}^{2},w_{1}=y_{0}y_{1},...$Then you get a 5-th order equation if
and if $W(3,0)=0$ (see [1],[7],[8]). This was exploited by Zudilin in [23],
who found formulas for $\dfrac{1}{\pi ^{2}}$ by "squaring" formulas for $%
\dfrac{1}{\pi }$ (also all formulas in Baruah-Berndt [9] can be obtained by
"squaring").

One can also take the wronskians of the solutions of a 3-rd order equation.
Since $\dbinom{3}{2}=3$ you get back a 3-rd order equation. Some examples of
this will be given in Examples 33-37.

\textbf{Introduction.}

Almost a hundred years ago, Ramanujan found 17 series for $\dfrac{1}{\pi }$
e.g.%
\[
\sum_{n=0}^{\infty }\binom{2n}{n}\binom{3n}{n}\binom{6n}{3n}(10177+261702n)%
\frac{1}{(-5280^{3})^{n}}=\frac{880^{2}\sqrt{330}}{\pi } 
\]%
If one knows the argument \ $-\dfrac{1}{5280^{3}}$ \ then the rest of the
formula can be found by using PSLQ in Maple. But also \ $-\dfrac{1}{5280^{3}}
$ \ can be found as \ $-5280^{3}=J(-\exp (-\pi \sqrt{67}))$ \ where%
\[
J(q)=\frac{1}{q}+744+196884q+21493760q^{2}+... 
\]%
is the modular invariant. There are many definitions of \ $J(q)$ \ but the
most elementary is the following. The power series%
\[
y_{0}=\sum_{n=0}^{\infty }\binom{2n}{n}\binom{3n}{n}\binom{6n}{3n}x^{n} 
\]%
satisfies the hypergeometric differential equation%
\[
(\theta ^{3}-1728x(\theta +\frac{1}{2})(\theta +\frac{1}{6})(\theta +\frac{5%
}{6}))y_{0}=0 
\]%
where \ $\theta =x\dfrac{d}{dx}.$ A second solution is given by%
\[
y_{1}=y_{0}\log (x)+744x+562932x^{2}+570443360x^{3}+... 
\]%
Define%
\[
q=\exp (\frac{y_{1}}{y_{0}})=x+744x^{2}+750420x^{3}+872769632x^{4}+... 
\]%
and solving for \ $x$%
\[
x=x(q)=q-744q^{2}+356652q^{3}-140361152q^{4}+... 
\]%
Then%
\[
J(q)=\frac{1}{x(q)}=\frac{1}{q}+744+196884q+21493760q^{2}+... 
\]%
Similarly can all of Ramanujan's formulas for \ $\dfrac{1}{\pi }$ \ be found
using the corresponding \ $J-$functions belonging to the third order
differential equations of the hypergeometric functions \ $_{2}F_{1}(\dfrac{1%
}{2},s,1-s;1,1;x)$ \ where \ $s=\dfrac{1}{2},\dfrac{1}{3},\dfrac{1}{4}.$

The simplest non-hypergeometric differential equations of order three are of
degree two and listed as \ $\alpha ,\beta ,\gamma ,...,\kappa $ \ in [1],[2]
,where they are used to construct Calabi-Yau differential equations of order
four. We computed the corresponding J-functions and listed \ $I=$ the number
of integer values, $S=$ the number of values containing square roots, and $%
C= $ the number of cube roots of the form \ $c_{0}+c_{1}m^{1/3}+c_{2}m^{2/3}$
with $c_{0,}c_{1},c_{2},m$ integers. We found only five values for $\
m=2,3,10,28,98.$ In the last column we give the number of the $\ J-$function
in Conway-Norton's list [13].%
\[
\text{\textbf{Table 1.}} 
\]%
\[
\begin{tabular}{|l|l|l|l|l|l|}
\hline
Name & $A_{n}$ & $I$ & $S$ & $C$ & \# \\ \hline
$\alpha $ & $\sum_{k}\binom{n}{k}^{2}\binom{2k}{k}\binom{2n-2k}{n-k}$ & $6$
& $17$ & $3$ & 6C \\ \hline
$\beta $ & $\sum_{k}\binom{2k}{k}^{2}\binom{2n-2k}{n-k}^{2}$ & $2$ & $7$ & $%
- $ & 4C \\ \hline
$\gamma $ & $\sum_{k}\binom{n}{k}^{2}\binom{n+k}{n}^{2}$ & $1$ & $12$ & $1$
& 6B \\ \hline
$\delta $ & $\sum_{k}(-1)^{k}3^{n-3k}\binom{n}{3k}\binom{n+k}{n}\frac{(3k)!}{%
k!^{3}}$ & $4$ & $12$ & $1$ & 6D \\ \hline
$\epsilon $ & $\sum_{k}\binom{n}{k}^{2}\binom{2k}{n}^{2}$ & $-$ & $13$ & $-$
& 8A \\ \hline
$\zeta $ & $\sum_{j,k}\binom{n}{j}^{2}\binom{n}{k}\binom{j}{k}\binom{j+k}{n}$
& $1$ & $1$ & $1$ & 9A \\ \hline
$\eta $ & $\sum_{k}(-1)^{k}\binom{n}{k}^{3}\left\{ \binom{4n-5k-1}{3n}+%
\binom{4n-5k}{3n}\right\} $ & $-$ & $11$ & $-$ & 5B \\ \hline
$\vartheta $ & $64^{n}\sum_{k}\binom{-1/4}{k}^{2}\binom{-3/4}{n-k}^{2}$ & $6$
& $21$ & $1$ & 2B \\ \hline
$\iota $ & $27^{n}\sum_{k}\binom{-1/3}{k}^{2}\binom{-2/3}{n-k}^{2}$ & $3$ & $%
14$ & $2$ & 2B \\ \hline
$\kappa $ & $432^{n}\sum_{k}\binom{-1/6}{k}^{2}\binom{-5/6}{n-k}$ & $-$ & $%
11 $ & $-$ & $-$ \\ \hline
\end{tabular}%
\]

Some of the smallest $J-$values might lead to divergent series.

Another building block for Calabi-Yau equations is of order two (denoted by
(a), (b),(c),...,(j) in [1]). If 
\[
u_{0}=\sum_{n=0}^{\infty }A_{n}x^{n} 
\]%
is the solution then 
\[
y_{0}=u_{0}^{2}=\sum_{n=0}^{\infty }c_{n}x^{n} 
\]%
satisfies a third order differential equation with the same J-function. We
get%
\[
\text{\textbf{Table 2.}} 
\]%
\[
\begin{tabular}{|l|l|l|l|l|l|}
\hline
Name & $A_{n}$ & $I$ & $S$ & $C$ & \# \\ \hline
(a) & $\sum_{k}\binom{n}{k}^{3}$ & $1$ & $5$ & $1$ & 12B \\ \hline
(b) & $\sum_{k}\binom{n}{k}^{2}\binom{n+k}{n}$ & $-$ & $-$ & $-$ & $-$ \\ 
\hline
(c) & $\sum_{k}\binom{n}{k}^{2}\binom{2k}{k}$ & $1$ & $5$ & $1$ & 12B \\ 
\hline
(d) & $\sum_{k}\binom{n}{k}\binom{2k}{k}\binom{2n-k}{n-k}$ & $-$ & $2$ & $-$
& 8D \\ \hline
(e) & $\sum_{k}4^{n-k}\binom{2k}{k}^{2}\binom{2n-k}{n-k}$ & $2$ & $8$ & $-$
& 4C \\ \hline
(f) & $\sum_{k}(-1)^{k}3^{n-3k}\binom{n}{3k}\frac{(3k)!}{k!^{3}}$ & $1$ & $0$
& $2$ & $-$ \\ \hline
(g) & $\sum_{j,k}(-1)^{j}8^{n-j}\binom{n}{j}\binom{j}{k}^{3}$ & $1$ & $5$ & $%
1$ & 12B \\ \hline
(h) & $27^{n}\sum_{k}\binom{-2/3}{k}\binom{-1/3}{n-k}^{2}$ & $3$ & $14$ & $2$
& 3B \\ \hline
(i) & $64^{n}\sum_{k}\binom{-3/4}{k}\binom{-1/4}{n-k}^{2}$ & $6$ & $21$ & $1$
& 2B \\ \hline
(j) & $432^{n}\sum_{k}\binom{-5/6}{k}\binom{-1/6}{n-k}^{2}$ & $1$ & $11$ & $%
- $ & $-$ \\ \hline
\end{tabular}%
\]

Here we notice that the J-functions of (a),(c) and (g) are translations of
each other:%
\[
J(\text{a,c,g})=\frac{1}{q}%
+(3,4,5)+6q+4q^{2}-3q^{3}-12q^{4}-8q^{5}+12q^{6}-... 
\]

Another way to get a third order differential equation is to form the
Hadamard product by multiplying$~$the coefficients $A_{n}$ by \ $\dbinom{2n}{%
n}$. We get%
\[
\text{\textbf{Table 3.}} 
\]%
\[
\begin{tabular}{|l|l|l|l|l|}
\hline
$A_{n}$ & $I$ & $S$ & $C$ & \# \\ \hline
$\dbinom{2n}{n}\ast $(a) & $14$ & $22$ & $3$ & 6A \\ \hline
$\dbinom{2n}{n}\ast $(b) & $6$ & $34$ & $-$ & 5A \\ \hline
$\dbinom{2n}{n}\ast $(c) & $14$ & $22$ & $3$ & 6A \\ \hline
$\dbinom{2n}{n}\ast $(d) & $11$ & $28$ & $-$ & 4B \\ \hline
$\dbinom{2n}{n}\ast $(e) & $5$ & $20$ & $1$ & 2B \\ \hline
$\dbinom{2n}{n}\ast $(f) & $1$ & $-$ & $2$ & 3C \\ \hline
$\dbinom{2n}{n}\ast $(g) & $14$ & $22$ & $3$ & 6A \\ \hline
$\dbinom{2n}{n}\ast $(h) & $10$ & $24$ & $3$ & 3A \\ \hline
$\dbinom{2n}{n}\ast $(i) & $14$ & $38$ & $1$ & 2A \\ \hline
$\dbinom{2n}{n}\ast $(j) & $12$ & $44$ & $-$ & 1A \\ \hline
\end{tabular}%
\]

The formulas for \ $\dfrac{1}{\pi }$ for the integer values are easiast
found by using \ PSLQ([a0\symbol{94}2,a0*a1,a1\symbol{94}2,1/Pi\symbol{94}%
2]); where

ak=sum(A(n)*n\symbol{94}k*x0\symbol{94}n,n=0..100) for k=0,1. The formulas
for the square and cube roots are trickier, sometimes it is easier to find
the "squared" formula for $\dfrac{1}{\pi ^{2}}$

\textbf{Reflexive Polytopes.}

Reflexive polytopes with only one inner point are important in string theory
(see [18] for a definition). In dimension 2 there are only 16 of them (for a
picture see [18]). We will only consider the following case:

\textbf{Example 0.}

Let the vertices be given by (case B$_{15}$ in [18])%
\[
\left( 
\begin{tabular}{llll}
$0$ & $2$ & $-1$ & $-1$ \\ 
$1$ & $-1$ & $1$ & $-1$%
\end{tabular}%
\right) 
\]%
The corresponding Laurent polynomial is%
\[
S=y+\frac{x^{2}}{y}+\frac{y}{x}+\frac{1}{xy} 
\]%
We form the power series 
\[
u=\sum_{n=0}^{\infty }a_{n}x^{n} 
\]%
where%
\[
a_{n}=\text{Constant Term}(S^{2n}) 
\]%
We have%
\[
S^{m}=\sum_{i+j+k+l=m}\frac{m!}{i!j!k!l!}y^{i}(\frac{x^{2}}{y})^{j}(\frac{y}{%
x})^{k}(\frac{1}{xy})^{l} 
\]%
\[
=\sum_{i+j+k+l=m}\frac{m!}{i!j!k!l!}x^{2j-k-l}y^{i-j+k-l} 
\]%
So we get the system of equations

$\ \ \ \ \ \ \ \ \ \ \ \ \ \ \ \ \ \ \ \ \ \ \ \ \ \ \ \ \ \ \ \ \ \ \ \ \ \
\ \ \ \ \ \ \ \ \ \ \ \ \ \ \ \ \ \ \ \ \ \ \ \ \ \ \ \ \ \ \ \ \ \ \ \ \ \
\ \ \ \ \ \ \ 2j-k-l=0$

$\ \ \ \ \ \ \ \ \ \ \ \ \ \ \ \ \ \ \ \ \ \ \ \ \ \ \ \ \ \ \ \ \ \ \ \ \ \
\ \ \ \ \ \ \ \ \ \ \ \ \ \ \ \ \ \ \ \ \ \ \ \ \ \ \ \ \ \ \ \ \ \ \ \ \ \
\ \ \ \ \ \ \ i-j+k-l=0$

$\ \ \ \ \ \ \ \ \ \ \ \ \ \ \ \ \ \ \ \ \ \ \ \ \ \ \ \ \ \ \ \ \ \ \ \ \ \
\ \ \ \ \ \ \ \ \ \ \ \ \ \ \ \ \ \ \ \ \ \ \ \ \ \ \ \ \ \ \ \ \ \ \ \ \ \
\ \ \ \ \ \ \ i+j+k+l=m$

with solution \ ( $m=2n$ )

$\ \ \ \ \ \ \ \ \ \ \ \ \ \ \ \ \ \ \ \ \ \ \ \ \ \ \ \ \ \ \ \ \ \ \ \ \ \
\ \ \ \ \ \ \ \ \ \ \ \ \ \ \ \ \ \ \ \ \ \ \ \ \ \ \ \ \ \ \ \ \ \ \ \ \ \
\ \ \ \ \ \ \ i=2n-j$

$\ \ \ \ \ \ \ \ \ \ \ \ \ \ \ \ \ \ \ \ \ \ \ \ \ \ \ \ \ \ \ \ \ \ \ \ \ \
\ \ \ \ \ \ \ \ \ \ \ \ \ \ \ \ \ \ \ \ \ \ \ \ \ \ \ \ \ \ \ \ \ \ \ \ \ \
\ \ \ \ \ \ \ k=3j-n$

$\ \ \ \ \ \ \ \ \ \ \ \ \ \ \ \ \ \ \ \ \ \ \ \ \ \ \ \ \ \ \ \ \ \ \ \ \ \
\ \ \ \ \ \ \ \ \ \ \ \ \ \ \ \ \ \ \ \ \ \ \ \ \ \ \ \ \ \ \ \ \ \ \ \ \ \
\ \ \ \ \ \ \ l=n-j$

We get%
\[
a_{n}=\sum_{j=0}^{n}\frac{(2n)!}{(2n-3j)!j!(3j-n)!(n-j)!}=\sum_{j=0}^{n}%
\binom{2n}{3j}\binom{3j}{n}\binom{n}{j} 
\]%
so%
\[
u_{0}=1+24x^{2}+120x^{3}+2520x^{4}+... 
\]%
with differential operator%
\[
\theta ^{2}-x(13\theta +1)-x^{2}(102\theta ^{2}+288\theta
+96)+216x^{3}(2\theta +1)^{2} 
\]%
The second solution is%
\[
u_{1}=u_{0}\log (x)+x+\frac{103}{2}x^{2}+\frac{985}{3}x^{3}+... 
\]%
and we get the $J-$function%
\[
J_{0}=\frac{1}{q}%
+1+51q+202q^{2}+177q^{3}-8202q^{4}-77876q^{5}-100500q^{6}+... 
\]%
Now we form%
\[
y_{0}=u_{0}^{2}=1+48x^{2}+240x^{3}+5616x^{4}+56160x^{5}+... 
\]%
which satisfies a third order differential equation whose $J-$fuction is the
same since \ $y_{1}=u_{0}u_{1}$ \ and \ $\dfrac{y_{1}}{y_{0}}=\dfrac{u_{1}}{%
u_{0}}$. We have \ $y_{0}=\sum_{0}^{\infty }c_{n}x^{n}$ \ where%
\[
c_{n}=\sum_{j,k,l=0}^{n}\binom{2j}{3k}\binom{3k}{j}\binom{j}{k}\binom{2n-2j}{%
3l}\binom{3l}{n-j}\binom{n-j}{l} 
\]%
Now we are ready to find the formulas for \ $\dfrac{1}{\pi }$

\textbf{0.1. }%
\[
J_{0}(\exp (-\pi \sqrt{2})=8(5+\sqrt{34}) 
\]%
\[
\sum_{n=0}^{\infty }c_{n}\left\{ -68+19\sqrt{34}+(34+28\sqrt{34})n\right\} 
\frac{1}{(8(5+\sqrt{34}))^{n}}=\frac{102\sqrt{2}}{\pi } 
\]

\textbf{0.2.}%
\[
J_{0}(\exp (-2\pi )=24(11+\sqrt{129}) 
\]%
\[
\sum_{n=0}^{\infty }c_{n}\left\{ -1+\sqrt{129}+(31+3\sqrt{129})n\right\} 
\frac{1}{(24(11+\sqrt{129}))^{n}}=\frac{129+3\sqrt{129}}{5\pi } 
\]

\textbf{0.3.}%
\[
J_{0}(-\exp (-\pi )=-24 
\]%
\[
\sum_{n=0}^{\infty }c_{n}\left\{ 1+5n\right\} \frac{1}{(-24)^{n}}=\frac{6}{%
\pi } 
\]

\textbf{0.4.}%
\[
J_{0}(-\exp (-\pi \sqrt{3})=-24(5+\sqrt{21}) 
\]%
\[
\sum_{n=0}^{\infty }c_{n}\left\{ 10+(23+7\sqrt{21})n\right\} \frac{1}{(-24(5+%
\sqrt{21}))^{n}}=\frac{12\sqrt{7}}{5\pi } 
\]

\textbf{0.5.}%
\[
J_{0}(-\exp (-\pi \sqrt{7})=-24(255+\sqrt{7161}) 
\]

\[
\sum_{n=0}^{\infty }c_{n}\left\{ -16709+849\sqrt{7161}+(231539+2681\sqrt{7161%
})n\right\} \frac{1}{(-24(255+\sqrt{7161}))^{n}}=\frac{65472\sqrt{7}}{\pi } 
\]

\[
\]

\textbf{Dimension 3.}

\textbf{Example 1.}

We will consider the following polytope, an octagon,%
\[
\left( 
\begin{tabular}{llllll}
$1$ & $0$ & $0$ & $-1$ & $0$ & $0$ \\ 
$0$ & $1$ & $0$ & $0$ & $-1$ & $0$ \\ 
$0$ & $0$ & $1$ & $0$ & $0$ & $-1$%
\end{tabular}%
\right) 
\]%
with Laurent polynomial%
\[
S=x+\frac{1}{x}+y+\frac{1}{y}+z+\frac{1}{z} 
\]%
which gives%
\[
A_{n}=C.T.(S^{2n})=\sum_{i+j+k=n}\frac{(2n)!}{i!^{2}j!^{2}k!^{2}} 
\]%
Of some strange reason, Maple takes forever already to compute,say \ $A_{50}$
\ ,so we found a formula with binomial coefficients instead%
\[
A_{n}=\sum_{i,j}\binom{2n}{i}\binom{2n-i}{i}\binom{2n-2i}{j}\binom{2n-2i-j}{j%
}\binom{2n-2i-2j}{n-i-j} 
\]%
which instantly computes \ $A_{50}.$

Then%
\[
y_{0}=\sum_{n=0}^{\infty }A_{n}x^{n}=1+6x+90x^{2}+1860x^{3}+44730x^{4}+... 
\]%
with differential operator%
\[
\theta ^{3}-2x(2\theta +1)(10\theta ^{2}+10\theta +3)+36x^{2}(\theta
+1)(2\theta +1)(2\theta +3) 
\]%
which gives%
\[
J_{1}=\frac{1}{q}+14+79q+352q^{2}+1431q^{3}+4160q^{4}+...\text{ \ 6A} 
\]%
Here we note that the coefficients of%
\[
J_{1}=\frac{1}{q}+\sum_{n=0}^{\infty }c(n)q^{n} 
\]%
satisfy the nonlinear identities of replicable functions, e.g. \ $%
c(10)=c(6)+c(2)c(3)+c(1)c(4)$ . This is not the case for the two-dimensional
case.

We find the following formulas for \ $\dfrac{1}{\pi }$

\textbf{1.1.}%
\[
J_{1}(\exp (-\pi \sqrt{4/3}))=54 
\]%
\[
\sum_{n=0}^{\infty }A_{n}(2+10n)\frac{1}{54^{n}}=\frac{9\sqrt{3}}{\pi } 
\]

\textbf{1.2. }%
\[
J_{1}(\exp (-\pi \sqrt{2}))=100 
\]%
\[
\sum_{n=0}^{\infty }A_{n}(3+16n)\frac{1}{100^{n}}=\frac{25}{\pi \sqrt{3}} 
\]

\textbf{1.3.}%
\[
J_{1}(\exp (-\pi \sqrt{10/3}))=324 
\]%
\[
\sum_{n=0}^{\infty }A_{n}(13+80n)\frac{1}{324^{n}}=\frac{27\sqrt{3}}{\pi } 
\]

\textbf{1.4.}%
\[
J_{1}(\exp (-\pi \sqrt{14/3}))=900 
\]%
\[
\sum_{n=0}^{\infty }A_{n}(16+112n)\frac{1}{900^{n}}=\frac{75}{\pi \sqrt{2}} 
\]

\textbf{1.5.}%
\[
J_{1}(\exp (-\pi \sqrt{26/3}))=10404 
\]%
\[
\sum_{n=0}^{\infty }A_{n}(112+1040n)\frac{1}{10404^{n}}=\frac{289}{\pi }%
\sqrt{3/2} 
\]

\textbf{1.6.}%
\[
J_{1}(\exp (-\pi \sqrt{34/3}))=39204 
\]%
\[
\sum_{n=0}^{\infty }A_{n}(1798+19040n)\frac{1}{39204^{n}}=\frac{3267}{\pi }%
\sqrt{3} 
\]

\textbf{1.7.}%
\[
J_{1}(-\exp (-\pi \sqrt{5/3}))=-45 
\]%
\[
\sum_{n=0}^{\infty }A_{n}(2+7n)\frac{1}{(-45)^{n}}=\frac{\sqrt{15}}{\pi } 
\]

\textbf{1.8.}%
\[
J_{1}(-\exp (-\pi \sqrt{7/3}))=-108 
\]%
\[
\sum_{n=0}^{\infty }A_{n}(13+56n)\frac{1}{(-108)^{n}}=\frac{18\sqrt{3}}{\pi }
\]

\textbf{1.9.}%
\[
J_{1}(-\exp (-\pi \sqrt{11/3}))=-396 
\]%
\[
\sum_{n=0}^{\infty }A_{n}(7+40n)\frac{1}{(-396)^{n}}=\frac{6\sqrt{11}}{\pi } 
\]

\textbf{1.10.}%
\[
J_{1}(-\exp (-\pi \sqrt{19/3}))=-2700 
\]%
\[
\sum_{n=0}^{\infty }A_{n}(253+1976n)\frac{1}{(-2700)^{n}}=\frac{450\sqrt{3}}{%
\pi } 
\]

\textbf{1.11.}%
\[
J_{1}(-\exp (-\pi \sqrt{31/3}))=-24300 
\]%
\[
\sum_{n=0}^{\infty }A_{n}(2239+22568n)\frac{1}{(-24300)^{n}}=\frac{4050\sqrt{%
3}}{\pi } 
\]

\textbf{1.12.}%
\[
J_{1}(-\exp (-\pi \sqrt{59/3}))=-1123596 
\]%
\[
\sum_{n=0}^{\infty }A_{n}(2587+36040n)\frac{1}{(-1123596)^{n}}=\frac{1058%
\sqrt{59}}{\pi } 
\]

\textbf{1.13.}%
\[
J_{1}(-\exp (-\pi \sqrt{5}))=-4(139+80\sqrt{3}) 
\]%
\[
\sum_{n=0}^{\infty }A_{n}\{-1773+1104\sqrt{3}+(-3480+2560\sqrt{3})n\}\frac{1%
}{(-4(139+80\sqrt{3}))^{n}}=\frac{242\sqrt{3}}{\pi } 
\]

\textbf{1.14.}%
\[
J_{1}(-\exp (-\pi \sqrt{\frac{119}{3}}))=-195890796-113097600\sqrt{3} 
\]%
\[
\sum_{n=0}^{\infty }A_{n}\frac{-7511+5125\sqrt{3}+(-6440+19320\sqrt{3})n}{%
(-195890796-113097600\sqrt{3})^{n}}=\frac{282\sqrt{17}+438\sqrt{51}}{\pi } 
\]

\textbf{1.15.}%
\[
J_{1}(-\exp (-\pi \sqrt{3}))=-4(2+3\cdot 2^{1/3}+2^{2/3})^{2} 
\]%
\[
\sum_{n=0}^{\infty }A_{n}\frac{\{3-9\cdot 2^{1/3}+3\cdot 2^{2/3}+(2-26\cdot
2^{1/3}+8\cdot 2^{2/3})n\}}{(-4(2+3\cdot 2^{1/3}+2^{2/3})^{2})^{n}}=\frac{%
\sqrt{3}}{\pi }(4-6\cdot 2^{2/3}) 
\]

\textbf{1.16.}%
\[
J_{1}(-\exp (-\pi \sqrt{\frac{25}{3}}))=-2880-1368\cdot 10^{1/3}-612\cdot
10^{2/3} 
\]%
\[
\sum_{n=0}^{\infty }A_{n}\frac{-170+115\cdot 10^{1/3}+22\cdot
10^{2/3}+(-520+590\cdot 10^{1/3}+188\cdot 10^{2/3})n}{(-2880-1368\cdot
10^{1/3}-612\cdot 10^{2/3})^{n}}=\frac{18^{2}\sqrt{3}}{\pi } 
\]

\textbf{1.17}%
\[
J_{1}(-\exp (-\pi \sqrt{\frac{49}{3}}))=-108864-35964\cdot
28^{1/3}-11772\cdot 28^{2/3} 
\]%
\[
\sum_{n=0}^{\infty }A_{n}\frac{-119868+43036\cdot 28^{1/3}+11853\cdot
28^{2/3}+(-507696+344792\cdot 28^{1/3}+106866\cdot 28^{2/3})n}{%
(-108864-35964\cdot 28^{1/3}-11772\cdot 28^{2/3})^{n}} 
\]%
\[
=\frac{2\cdot 330^{2}\sqrt{3}}{\pi } 
\]

\textbf{1.18.}%
\[
J_{1}(-\exp (-\pi \sqrt{\frac{119}{3}}))=-195890796-113097600\sqrt{3} 
\]%
\[
\sum_{n=0}^{\infty }A_{n}\frac{-7511+5125\sqrt{3}+(-6440+19320\sqrt{3})n}{%
(-195890796-113097600\sqrt{3})^{n}}=\frac{282\sqrt{17}+438\sqrt{51}}{\pi } 
\]

\textbf{1.19.}%
\[
J_{1}(\exp (\frac{\pi i}{3}-\pi \sqrt{\frac{17}{9}}))=52+64i 
\]%
\[
\sum_{n=0}^{\infty }A_{n}\frac{15+6i+(72+16i)n}{(52+64i)^{n}}=\frac{%
(1350-300i)\sqrt{17}}{\pi \sqrt{7787+4016i}} 
\]

\textbf{1.20.}%
\[
J_{1}(\exp (\frac{2\pi i}{3}-\pi \sqrt{\frac{14}{9}}))=-12+16i\sqrt{7} 
\]%
\[
\sum_{n=0}^{\infty }A_{n}\frac{3i+9\sqrt{7}+32n\sqrt{7}}{(-12+16i\sqrt{7}%
)^{n}}=\frac{1452\sqrt{7}}{\pi \sqrt{3794+106i\sqrt{7}}} 
\]

There are at least 20 more values of \ $J_{1}$ leading to quadratic equations%
$.$Also \ $J_{1}(-\exp (-\pi ))=-12$ \ leads to a divergent series.

After writing this we found out that%
\[
A_{n}=\binom{2n}{n}\sum_{n=0}^{n}\binom{n}{k}^{2}\binom{2k}{k} 
\]%
and this case had been treated in [11] where the \ \textbf{1.1}$-$\textbf{1.6%
} \ are proved.

\bigskip

\textbf{Example 2,}

Consider the Laurent polynomial%
\[
S=x+\frac{1}{x}+y+\frac{1}{y}+z+\frac{xy}{z} 
\]%
We find 
\[
A_{n}=\binom{2n}{n}\sum_{i,j}\binom{n}{i}\binom{n}{j}\binom{i}{j}\binom{n-j}{%
i} 
\]%
and 
\[
y_{0}=\sum_{n=0}^{\infty }A_{n}x^{n} 
\]%
with the differential operator%
\[
\theta ^{3}-2x(2\theta +1)(7\theta ^{2}+7\theta +2)-32x^{2}(\theta
+1)(2\theta +1)(2\theta +3) 
\]%
This operator is wellknown and we recognize%
\[
A_{n}=\binom{2n}{n}\sum_{n=0}^{n}\binom{n}{k}^{3} 
\]%
which is also in \ [11]. Then%
\[
y_{0}=1+4x+60x^{2}+1120x^{3}+24220x^{4}+... 
\]%
and%
\[
J_{2}=\frac{1}{q}+10+79q+352q^{2}+1431q^{3}+4160q^{4}+...=J_{1}-4\text{ \ \
\ \ \ 6A} 
\]

Hence the denominator in \ $x_{0}$ \ will be \ $4$ \ less than in example 1.

\textbf{2.1.}%
\[
J_{2}(\exp (-\pi \sqrt{4/3}))=50 
\]%
\[
\sum_{n=0}^{\infty }A_{n}(4+18n)\frac{1}{50^{n}}=\frac{25}{\pi } 
\]

\textbf{2.2. }%
\[
J_{2}(\exp (-\pi \sqrt{2}))=96 
\]%
\[
\sum_{n=0}^{\infty }A_{n}(1+5n)\frac{1}{96^{n}}=\frac{3\sqrt{2}}{\pi } 
\]

\textbf{2.3}%
\[
J_{2}(\exp (-\pi \sqrt{10/3}))=320 
\]%
\[
\sum_{n=0}^{\infty }A_{n}(1+6n)\frac{1}{320^{n}}=\frac{8}{3\pi }\sqrt{\frac{5%
}{3}} 
\]

\textbf{2.4.}%
\[
J_{2}(\exp (-\pi \sqrt{14/3}))=896 
\]%
\[
\sum_{n=0}^{\infty }A_{n}(13+90n)\frac{1}{896^{n}}=\frac{16\sqrt{7}}{\pi } 
\]

\textbf{2.5.}%
\[
J_{2}(\exp (-\pi \sqrt{26/3}))=10400 
\]%
\[
\sum_{n=0}^{\infty }A_{n}(11+102n)\frac{1}{10400^{n}}=\frac{50}{3\pi }\sqrt{%
\frac{13}{3}} 
\]

\textbf{2.6.}%
\[
J_{2}(\exp (-\pi \sqrt{34/3}))=39200 
\]%
\[
\sum_{n=0}^{\infty }A_{n}(159+1683n)\frac{1}{39200^{n}}=\frac{1225}{\pi 
\sqrt{6}} 
\]

\textbf{2.7.}%
\[
J_{1}(-\exp (-\pi \sqrt{5/3}))=-49 
\]%
\[
\sum_{n=0}^{\infty }A_{n}(4+15n)\frac{1}{(-49)^{n}}=\frac{49}{3\pi \sqrt{3}} 
\]

\textbf{2.8.}%
\[
J_{2}(-\exp (-\pi \sqrt{7/3}))=-112 
\]%
\[
\sum_{n=0}^{\infty }A_{n}(2+9n)\frac{1}{(-112)^{n}}=\frac{2\sqrt{7}}{\pi } 
\]

\textbf{2.9.}%
\[
J_{2}(-\exp (-\pi \sqrt{11/3}))=-400 
\]%
\[
\sum_{n=0}^{\infty }A_{n}(17+99n)\frac{1}{(-400)^{n}}=\frac{50}{\pi } 
\]

\textbf{2.10.}%
\[
J_{2}(-\exp (-\pi \sqrt{19/3}))=-2704 
\]%
\[
\sum_{n=0}^{\infty }A_{n}(109+855n)\frac{1}{(-2704)^{n}}=\frac{338}{\pi } 
\]

\textbf{2.11.}%
\[
J_{2}(-\exp (-\pi \sqrt{31/3}))=-24304 
\]%
\[
\sum_{n=0}^{\infty }A_{n}(174+1755n)\frac{1}{(-24304)^{n}}=\frac{98\sqrt{31}%
}{\pi } 
\]

\textbf{2.12.}%
\[
J_{2}(-\exp (-\pi \sqrt{59/3}))=-1123600 
\]%
\[
\sum_{n=0}^{\infty }A_{n}(44709+622863n)\frac{1}{(-1123600)^{n}}=\frac{140450%
}{\pi } 
\]

\textbf{2.13.}%
\[
J_{2}(-\exp (-\pi \sqrt{5}))=-80(7+4\sqrt{3}) 
\]%
\[
\sum_{n=0}^{\infty }A_{n}\{2-2\sqrt{3}+(2-7\sqrt{3})n\}\frac{1}{(-80(7+4%
\sqrt{3}))^{n}}=\frac{2\sqrt{5}}{\pi } 
\]

\textbf{2.14.}%
\[
J_{2}(-\exp (-\pi \sqrt{3}))=-4(3+2^{1/3}+2\cdot 2^{2/3})^{2} 
\]%
\[
\sum_{n=0}^{\infty }A_{n}\frac{\{24+12\cdot 2^{1/3}-15\cdot
2^{2/3}+(78+39\cdot 2^{1/3}-30\cdot 2^{2/3})n\}}{(-4(3+2^{1/3}+2\cdot
2^{2/3})^{2})^{n}}=\frac{25\sqrt{3}}{\pi } 
\]

\textbf{2.15}%
\[
J_{2}(-\exp (-\pi \sqrt{\frac{25}{3}}))=-\frac{1}{3}(106+46\cdot
10^{1/3}+16\cdot 10^{2/3})^{2} 
\]%
\[
\sum_{n=0}^{\infty }A_{n}\frac{\{241200-408620\cdot 10^{1/3}-18175\cdot
10^{2/3}+(13093650-530775\cdot 10^{1/3}+244710\cdot 10^{2/3})n\}}{(-\dfrac{1%
}{3}(106+46\cdot 10^{1/3}+16\cdot 10^{2/3})^{2})^{n}} 
\]%
\[
=\frac{23^{2}47^{2}\sqrt{15}}{\pi } 
\]

\textbf{2.16.}%
\[
J_{2}(-\exp (-\pi \sqrt{\frac{49}{3}}))=-\frac{1}{3}(538+184\cdot
28^{1/3}+67\cdot 28^{2/3})^{2} 
\]%
\[
\sum_{n=0}^{\infty }A_{n}\frac{a+bn}{(-\dfrac{1}{3}(538+184\cdot
28^{1/3}+67\cdot 28^{2/3})^{2})^{n}}=\frac{133334450\sqrt{7}}{\pi } 
\]%
where%
\[
a=187181568-14524874\cdot 28^{1/3}-3336620\cdot 28^{2/3} 
\]%
\[
b=1425786576-6644988\cdot 28^{1/3}+2185785\cdot 28^{2/3} 
\]

\textbf{2.17}%
\[
J_{2}(\exp (\frac{\pi i}{3}-\pi \sqrt{\frac{17}{9}}))=48+64i 
\]%
\[
\sum_{n=0}^{\infty }A_{n}\frac{8-3i+17(2-i)n}{(48+64i)^{n}}=\frac{24-18i}{%
\pi } 
\]

\textbf{2,18.}%
\[
J_{2}(\exp (\frac{2\pi i}{3}-\pi \sqrt{\frac{14}{9}}))=16+16i\sqrt{7} 
\]%
\[
\sum_{n=0}^{\infty }A_{n}\frac{5i+7\sqrt{7}+(14i+26\sqrt{7})n}{(16+16i\sqrt{7%
})^{n}}=\frac{48}{\pi } 
\]%
\[
\]%
\textbf{\ \ }

There are at least 20 more values of \ $J_{2}$ leading to quadratic equations%
$.$Also \ $J_{1}(-\exp (-\pi ))=-16$ \ leads to a divergent series.

\bigskip

\textbf{Example 3.}

Consider%
\[
S=x+\frac{1}{x}+y+\frac{1}{y}+z+\frac{1}{z}+\frac{xy}{z} 
\]%
Then%
\[
A_{n}=\binom{2n}{n}\sum_{i,j,k}\binom{n}{i}\binom{n}{j}\binom{n-j}{i}\binom{i%
}{k}\binom{n-i}{j+k}=\binom{2n}{n}\sum_{k=0}^{n}\binom{n}{k}^{2}\binom{n+k}{n%
} 
\]%
with%
\[
y_{0}=\sum_{n=0}^{\infty }A_{n}x^{n}=1+6x+114x^{2}+2940x^{3}+87570x^{4}+... 
\]%
and%
\[
J_{3}=\frac{1}{q}+16+134q+760q^{2}+3345q^{3}+12256q^{4}+...\text{ \ \ \ \ 5A}
\]%
\[
\theta ^{3}-2x(2\theta +1)(11\theta ^{2}+11\theta +3)-4x^{2}(\theta
+1)(2\theta +1)(2\theta +3) 
\]

\textbf{3.1.}%
\[
J_{3}(\exp (-\sqrt{\frac{8}{5}}))=72 
\]%
\[
\sum_{n=0}^{\infty }A_{n}(1+5n)\frac{1}{72^{n}}=\frac{9}{\pi \sqrt{2}} 
\]

\textbf{3.2.}%
\[
J_{3}(\exp (-\sqrt{\frac{12}{5}}))=147 
\]%
\[
\sum_{n=0}^{\infty }A_{n}(2+11n)\frac{1}{147^{n}}=\frac{49\sqrt{3}}{10\pi } 
\]

\textbf{3.3.}%
\[
J_{3}(-\exp (-\sqrt{\frac{23}{5}}))=-828 
\]%
\[
\sum_{n=0}^{\infty }A_{n}(29+180n)\frac{1}{(-828)^{n}}=\frac{18\sqrt{23}}{%
\pi } 
\]

\textbf{3.4.}%
\[
J_{3}(-\exp (-\sqrt{\frac{47}{5}}))=-15228 
\]%
\[
\sum_{n=0}^{\infty }A_{n}(355+3410n)\frac{1}{(-15228)^{n}}=\frac{162\sqrt{47}%
}{\pi } 
\]

\textbf{3.5.}%
\[
J_{3}(-\exp (-\sqrt{\frac{311}{5}}))=-108(800002423+357771960\sqrt{5}) 
\]%
\[
\sum_{n=0}^{\infty }A_{n}\frac{\left\{ 128732-48919\sqrt{5}+(484000-2090%
\sqrt{5})n\right\} }{(-108(800002423+357771960\sqrt{5}))^{n}}=\frac{18}{\pi }%
\sqrt{5696283+2551048\sqrt{5}} 
\]

\bigskip There are at least 30 more formulas containing square roots.

\textbf{Example 4.}

Here we have 
\[
S=x+\frac{1}{x}+yz+\frac{1}{yz}+xyz+\frac{1}{xyz} 
\]%
with%
\[
A_{n}= 
\]%
\[
\sum_{i,j,k}\binom{n}{i+2j+2k}\binom{n}{2i+j+2k}\binom{2i+2j+3k}{n}\binom{%
2i+j+2k}{i+k}\binom{i+j+k}{i}\binom{j+k}{k}\binom{2i+2j+3k}{i+k}^{-1} 
\]%
and%
\[
u_{0}=\sum_{n=0}^{\infty }A_{n}x^{n}=1+6x^{2}+12x^{3}+90x^{4}+360x^{5}+... 
\]%
satisfying%
\[
\theta ^{2}-x\theta (\theta +1)-24x^{2}(\theta +1)^{2}-36x^{3}(\theta
+1)(\theta +2) 
\]%
\[
J_{4}=\frac{1}{q}+1+6q+4q^{2}-3q^{3}-12q^{4}-8q^{5}+12q^{6}+...\text{ \ \ \
\ 12B} 
\]%
We form%
\[
y_{0}=u_{0}^{2}=\sum_{n=0}^{\infty }c_{n}x^{n} 
\]%
Then we have

\textbf{4.1.}%
\[
J_{4}(-\exp (-\pi \sqrt{\frac{5}{3}}))=-30-12\sqrt{5} 
\]%
\[
\sum_{n=0}^{\infty }c_{n}\frac{\left\{ 1+(2+\sqrt{5})n\right\} }{(-6(5+2%
\sqrt{5}))^{n}}=\frac{15\sqrt{3}}{8\pi } 
\]

\textbf{4.2.}%
\[
J_{4}(-\exp (-\pi \sqrt{3}))=-78-60\cdot 2^{1/3}-48\cdot 2^{2/3} 
\]%
\[
\sum_{n=0}^{\infty }c_{n}\frac{\{-144+252\cdot 2^{1/3}-66\cdot
2^{2/3}+(18+156\cdot 2^{1/3}+102\cdot 2^{2/3})n\}}{(-2(39+30\cdot
2^{1/3}+24\cdot 2^{2/3}))^{n}}=\frac{125\sqrt{3}}{\pi } 
\]%
There are four more cases containing square roots. This example coincides
with the hexagon in dimension two.

\bigskip

\textbf{Example 5.}

Let%
\[
S=x+y+z+\frac{yz}{x}+\frac{x}{yz}+\frac{xz}{y}+\frac{y}{xz} 
\]%
with%
\[
A_{n}=\text{C.T.}(S^{2n})=\binom{2n}{n}\sum_{j,k}\binom{n}{k}\binom{n}{2k-2j}%
\binom{n+2j-2k}{j} 
\]%
and%
\[
u_{0}=\sum_{n=0}^{\infty }A_{n}x^{n}=1+4x+48x^{2}+760x^{3}+13720x^{4}+... 
\]%
The third order differential equation factors into%
\[
\theta ^{2}-x(45\theta ^{2}+17\theta +4)+x^{2}(594\theta ^{2}+320\theta
+72)-64x^{3}(2\theta +1)(16\theta +5)-672x^{4}(2\theta +1)(2\theta +3) 
\]%
and we get%
\[
J_{5}=\frac{1}{q}%
+9+51q+202q^{2}+177q^{3}-8202q^{4}-77876q^{5}-100500q^{6}+...=J_{0}+8 
\]%
Let%
\[
y_{0}=u_{0}^{2}=\sum_{n=0}^{\infty }c_{n}x^{n} 
\]

\textbf{5.1.}%
\[
J_{5}(-\exp (-\pi \sqrt{3}))=-112-24\sqrt{21} 
\]%
\[
\sum_{n=0}^{\infty }c_{n}\{325-70\sqrt{21}+(570-120\sqrt{21})n\}\frac{1}{%
(-8(14+3\sqrt{21}))^{n}}=\frac{49\sqrt{3}-28\sqrt{7}}{\pi } 
\]

\textbf{5.2.}%
\[
J_{5}(-\exp (-\pi \sqrt{7}))=-2032-24\sqrt{7161} 
\]%
\[
\sum_{n=0}^{\infty }c_{n}\{-50680+720\sqrt{7161}+84000n\}\frac{1}{(-8(254+3%
\sqrt{7161}))^{n}}=\frac{13981\sqrt{7}-24\sqrt{7\cdot 7161}}{\pi } 
\]%
There are of course three more formulas just like in the first
two-dimensional example. Examples 4 and 5 are degenerated three-dimensional
polytopes which is indicated by that the differential equations are
factorable.

\bigskip

\textbf{Example 6.}

Let%
\[
S=x+\frac{1}{x}+y+\frac{1}{y}+z+\frac{1}{z}+xyx+\frac{1}{xyz} 
\]%
Then%
\[
A_{n}=\binom{2n}{n}^{3} 
\]%
so we have the case \ $s=\dfrac{1}{2}$ \ of Ramanujan. We get 
\[
J_{6}=\frac{1}{q}+24+276q+2048q^{2}+11202q^{3}+49152q^{4}+...\text{ \ \ \ 2B}
\]%
Computing \ $J(\pm \exp (-\pi \sqrt{n}))$ \ we find \ $5$ \ integer values, $%
19$ \ with square roots and $6$ satisfying third degree equations. One of
the last is%
\[
J_{6}(\exp (-\pi \sqrt{27}))=4^{4}(1+2^{1/3}+2^{2/3})^{8} 
\]%
and we find the formula%
\[
\sum_{n=0}^{\infty }\binom{2n}{n}^{3}\frac{\{27-30\cdot 2^{1/3}+6\cdot
2^{2/3}+(66-84\cdot 2^{1/3}+12\cdot 2^{2/3})n\}}{(2(1+2^{1/3}+2^{2/3}))^{8n}}%
=\frac{4}{\pi } 
\]

\textbf{Example 7.}

Let 
\[
S=x+y+z+\frac{1}{xz}+\frac{1}{yz} 
\]%
Then%
\[
A_{n}=\text{C.T.}(S^{3n})=\binom{2n}{n}^{2}\binom{3n}{n} 
\]%
so we are in Ramanujan's case \ $s=\dfrac{1}{3}$. We have 
\[
J_{7}=\frac{1}{q}+42+783q+8672q^{2}+65367q^{3}+371520q^{4}+...\text{ \ \ \ 3A%
} 
\]%
with $11$ integer values, $37$ square roots and $6$ cube roots of which six
are simple

\textbf{7.1}%
\[
J_{7}(-\exp (-3\pi ))=-64(64002+44377\cdot 3^{1/3}+30771\cdot 3^{2/3}) 
\]%
\[
\sum_{n=0}^{\infty }\binom{2n}{n}^{2}\binom{3n}{n}\frac{\{-2889+1053\cdot
3^{1/3}+1719\cdot 3^{2/3}+(-9027+12879\cdot 3^{1/3}+12717\cdot 3^{2/3})n\}}{%
(-64(64002+44377\cdot 3^{1/3}+30771\cdot 3^{2/3}))^{n}}=\frac{4000\sqrt{3}}{%
\pi } 
\]

\textbf{7.2.}%
\[
J_{7}(\exp (-\pi \sqrt{\frac{36}{3}}))=12(23+17\cdot 2^{1/3}+14\cdot
2^{2/3})^{2} 
\]%
\[
\sum_{n=0}^{\infty }\binom{2n}{n}^{2}\binom{3n}{n}\frac{12\cdot
2^{1/3}-6\cdot 2^{2/3}+(12+54\cdot 2^{1/3}-12\cdot 2^{2/3})n}{(12(23+17\cdot
2^{1/3}+14\cdot 2^{2/3})^{2})^{n}}=\frac{(7+2\cdot 2^{2/3})\sqrt{3}}{\pi } 
\]

\textbf{7.3.}%
\[
J_{7}(\exp (-\pi \sqrt{\frac{100}{3}}))=972(93+43\cdot 10^{1/3}+20\cdot
10^{2/3})^{2} 
\]%
\[
\sum_{n=0}^{\infty }\binom{2n}{n}^{2}\binom{3n}{n}\frac{\{160-68\cdot
10^{1/3}+10\cdot 10^{2/3}+(1020-150\cdot 10^{1/3}+84\cdot 10^{2/3})n\}}{%
(972(93+43\cdot 10^{1/3}+20\cdot 10^{2/3})^{2})^{n}}=\frac{243}{\pi }\sqrt{%
\frac{3}{5}} 
\]

\textbf{7.4.}%
\[
J_{7}(\exp (-\pi \sqrt{\frac{196}{3}}))=\frac{108}{7}(2305025863+499958298%
\cdot 98^{1/3}+108440559\cdot 98^{2/3}) 
\]%
\[
\sum_{n=0}^{\infty }\binom{2n}{n}^{2}\binom{3n}{n}\frac{-28812-2156\cdot
28^{1/3}+6272\cdot 28^{2/3}+(-24696+27342\cdot 28^{1/3}+55566\cdot 28^{2/3})n%
}{(\dfrac{108}{7}(2305025863+499958298\cdot 98^{1/3}+108440559\cdot
98^{2/3}))^{n}} 
\]%
\[
=\frac{(5235\cdot 98^{1/3}+120\cdot 98^{2/3})\sqrt{7}}{\pi } 
\]%
\textbf{.}

\textbf{7.5.}%
\[
J_{7}(-\exp (-\pi \sqrt{\frac{289}{3}}))=1728(4749270926+1847044269\cdot
17^{1/3}+718336053\cdot 17^{2/3}) 
\]%
\[
\sum_{n=0}^{\infty }\binom{2n}{n}^{2}\binom{3n}{n}\frac{\{-173689+72607\cdot
17^{1/3}+2159\cdot 17^{2/3}+(-486387+378981\cdot 17^{1/3}+53397\cdot
17^{2/3})n\}}{(1728(4749270926+1847044269\cdot 17^{1/3}+718336053\cdot
17^{2/3}))^{n}} 
\]%
\[
=\frac{12000\sqrt{51}}{\pi } 
\]

\textbf{7.6}%
\[
J_{7}(-\exp (-\pi \sqrt{\frac{361}{3}}))=6^{6}(6617733763+2480036604\cdot
19^{1/3}+929409036\cdot 19^{2/3}) 
\]%
\[
\sum_{n=0}^{\infty }\binom{2n}{n}^{2}\binom{3n}{n}\frac{\{168587-36442\cdot
19^{1/3}-2128\cdot 19^{2/3}+(1512951+20634\cdot 19^{1/3}+51756\cdot
19^{2/3})n\}}{(6^{6}(6617733763+2480036604\cdot 19^{1/3}+9294090366\cdot
19^{2/3}))^{n}} 
\]%
\[
=\frac{40500\sqrt{19}}{\pi } 
\]

\textbf{7.7.}%
\[
J_{7}(\exp (\frac{\pi i}{3}-\pi \sqrt{\frac{35}{9}}))=288+160i\sqrt{7} 
\]%
\[
\sum_{n=0}^{\infty }\binom{2n}{n}^{2}\binom{3n}{n}\frac{9i+21\sqrt{7}%
+(35i+135\sqrt{7})n}{(288+160i\sqrt{7})^{n}}=\frac{192}{\pi } 
\]

\textbf{7.8}%
\[
J_{7}(\exp (\frac{2\pi i}{3}-\pi \sqrt{\frac{32}{9}}))=-146+322i 
\]%
\[
\sum_{n=0}^{\infty }\binom{2n}{n}^{2}\binom{3n}{n}\frac{12+(68-4i)n}{%
(-146+322i)^{n}}=\frac{33-6i}{\pi } 
\]%
\textbf{.}

\textbf{Example 8.}

Let%
\[
S=x+y+z+\frac{1}{xz^{2}}+\frac{1}{yz^{2}} 
\]%
Then%
\[
A_{n}=C.T.(S^{4n})=\binom{2n}{n}^{2}\binom{4n}{2n} 
\]%
so we are in case $\ s=\dfrac{1}{4}.$We have%
\[
J_{8}=\frac{1}{q}+104+4372q+96256q^{2}+1240002q^{3}+...\text{ \ \ \ \ 2A} 
\]%
which gives $13$ integer values, $41$ square roots and $13$ cube roots of
which the following is simple

\textbf{8.1}%
\[
J_{8}(-\exp (-\pi \sqrt{27}))=-16(299+234\cdot 2^{1/3}+178\cdot 2^{2/3})^{2} 
\]%
\[
\sum_{n=0}^{\infty }\binom{2n}{n}^{2}\binom{4n}{2n}\frac{\{-35253+30024\cdot
2^{1/3}+20610\cdot 2^{2/3}+(-191865+314400\cdot 2^{1/3}+234240\cdot
2^{2/3})n\}}{(-16(299+234\cdot 2^{1/3}+178\cdot 2^{2/3})^{2})^{n}} 
\]%
\[
=\frac{64009\sqrt{3}}{\pi } 
\]

\textbf{8.2.}%
\[
J_{8}(\exp (-\pi \sqrt{\frac{14}{3}}))=(84+32\sqrt{6})^{4} 
\]%
\[
\sum_{n=0}^{\infty }\binom{2n}{n}^{2}\binom{4n}{2n}\frac{\{-139+138\sqrt{6}%
+140(2+11\sqrt{6})n\}}{(84+32\sqrt{6})^{4n}}=\frac{361\sqrt{3}}{\pi } 
\]

\textbf{Example 9.}

Let%
\[
S=x+y+yz+\frac{1}{x^{3}y^{2}z} 
\]%
with%
\[
A_{n}=C.T.(S^{6n})=\binom{2n}{n}\binom{3n}{n}\binom{6n}{3n} 
\]%
so we are in case $s=\dfrac{1}{6}.$ we have%
\[
J_{9}=\frac{1}{q}+744+196884q+21493760q^{2}+...\text{ \ \ \ \ 1A} 
\]%
with $12$ integer values, $28$ square roots and $10$ cube roots of which the
following is simple

\textbf{9.1.}%
\[
J_{9}(\exp (-\pi \sqrt{108}))=6000(8389623817+6658848836\cdot
2^{1/3}+5285131824\cdot 2^{2/3}) 
\]%
We have not been able to find the corresponding formula for $\dfrac{1}{\pi }%
. $ However we found the "squared" formula for $\dfrac{1}{\pi ^{2}}$%
\[
\sum_{n=0}^{\infty }B_{n}\frac{(c_{0}+c_{1}n+c_{2}n^{2})}{%
(6000(8389623817+6658848836\cdot 2^{1/3}+5285131824\cdot 2^{2/3}))^{n}}=%
\frac{(5\cdot 11\cdot 17)^{3}}{36\pi ^{2}} 
\]%
where%
\[
B_{n}=\sum_{k=0}^{n}\binom{2k}{k}\binom{3k}{k}\binom{6k}{3k}\binom{2n-2k}{n-k%
}\binom{3n-3k}{n-k}\binom{6n-6k}{3n-3k} 
\]%
\[
c_{0}=13739282-39647980\cdot 2^{1/3}+24262676\cdot 2^{2/3} 
\]%
\[
c_{1}=171995949-202768920\cdot 2^{1/3}+99903102\cdot 2^{2/3} 
\]%
\[
c_{2}=761257259-157058640\cdot 2^{1/3}+160025472\cdot 2^{2/3} 
\]

\textbf{9.2.}%
\[
J_{9}(-\exp (-\pi \sqrt{243})) 
\]%
\[
=2^{12}\cdot 1000\cdot (151022371885959+104713064226304\cdot
3^{1/3}+72603983653110\cdot 3^{2/3}) 
\]%
\[
\sum_{n=0}^{\infty }A_{n}\frac{1041847953-338114880\cdot
3^{1/3}-167858240\cdot 3^{2/3}+10041470758n}{(2^{12}\cdot 1000\cdot
(151022371885959+104713064226304\cdot 3^{1/3}+72603983653110\cdot
3^{2/3}))^{n}} 
\]%
\[
=\frac{203239368141920\sqrt{2530}}{\pi \sqrt{88880634564507-28907873679132%
\cdot 3^{1/3}+98392731804162\cdot 3^{2/3}}} 
\]%
\[
\sum_{n=0}^{\infty }B_{n}\frac{c_{0}+c_{1}n+c_{2}n^{2}}{(2^{12}\cdot
1000\cdot (151022371885959+104713064226304\cdot 3^{1/3}+72603983653110\cdot
3^{2/3}))^{n}} 
\]%
\[
=\frac{(40\cdot 11\cdot 23)^{3}}{\pi ^{2}} 
\]%
where%
\[
c_{0}=-7900578369198-2982917384712\cdot 3^{1/3}+5916926155212\cdot 3^{2/3} 
\]%
\[
c_{1}=-26429839595271-15744422803284\cdot 3^{1/3}+26095087504854\cdot
3^{2/3} 
\]%
\[
c_{2}=29626878188169-9635957893044\cdot 3^{1/3}+32797577268054\cdot 3^{2/3} 
\]

\textbf{Example 10.}

Let%
\[
S=x+y+z+\frac{1}{z}+\frac{1}{xyz^{3}} 
\]%
with%
\[
A_{n}=C.T.(S^{2n})=\dbinom{2n}{n}\sum_{k=0}^{[n/3]}\dbinom{n}{3k}\frac{(3k)!%
}{k!^{3}} 
\]%
and%
\[
y_{0}=\sum_{n=0}^{\infty }A_{n}x^{n}=1+2x+6x^{2}+140x^{3}+1750x^{4}+... 
\]%
Let%
\[
y_{0}^{2}=\sum_{n=0}^{\infty }c_{n}x^{n} 
\]

We get%
\[
J_{10}=\frac{1}{q}+4+248q+4124q^{2}+34752q^{3}+213126q^{4}+... 
\]%
which is 3C after the substitution \ $q\rightarrow q^{1/3}$.

\textbf{10.1.}%
\[
J_{10}(\exp (-\pi \sqrt{\frac{8}{3}}))=24 
\]%
\[
\sum_{n=0}^{\infty }A_{n}\frac{4+14n}{24^{n}}=\frac{9\sqrt{6}}{\pi } 
\]

\textbf{10.2}%
\[
J_{10}(\exp (-\pi \sqrt{\frac{16}{9}}))=70 
\]%
\[
\sum_{n=0}^{\infty }A_{n}\frac{8+36n}{70^{n}}=\frac{5\sqrt{35}}{\pi } 
\]

\textbf{10.3.}%
\[
J_{10}(\exp (-\pi \sqrt{\frac{28}{9}}))=259 
\]%
\[
\sum_{n=0}^{\infty }A_{n}\frac{180+1026n}{259^{n}}=\frac{37\sqrt{259}}{\pi } 
\]

\textbf{10.4}%
\[
J_{10}(-\exp (-\pi \sqrt{\frac{11}{9}}))=-28 
\]%
\[
\sum_{n=0}^{\infty }A_{n}\frac{7+22n}{(-28)^{n}}=\frac{6\sqrt{7}}{\pi } 
\]

\textbf{10.5.}%
\[
J_{10}(-\exp (-\pi \sqrt{\frac{19}{9}}))=-92 
\]%
\[
\sum_{n=0}^{\infty }A_{n}\frac{79+342n}{(-92)^{n}}=\frac{46\sqrt{23}}{\pi } 
\]

\textbf{10.6.}%
\[
J_{10}(-\exp (-\pi \sqrt{\frac{43}{9}}))=-956 
\]%
\[
\sum_{n=0}^{\infty }A_{n}\frac{2391+16254n}{(-956)^{n}}=\frac{478\sqrt{239}}{%
\pi } 
\]

\textbf{10.7}%
\[
J_{10}(-\exp (-\pi \sqrt{\frac{67}{9}}))=-5276 
\]%
\[
\sum_{n=0}^{\infty }A_{n}\frac{30607+261702n}{(-5276)^{n}}=\frac{2638\sqrt{%
1319}}{\pi } 
\]

\textbf{10.8.}%
\[
J_{10}(-\exp (-\pi \sqrt{\frac{163}{9}}))=-640316 
\]%
\[
\sum_{n=0}^{\infty }A_{n}\frac{40775675+545140134n}{(-640316)^{n}}=\frac{%
320158\sqrt{160079}}{\pi } 
\]

\textbf{10.9}%
\[
J_{10}(\exp (-\pi \sqrt{\frac{4}{3}}))=4+30\cdot 2^{1/3} 
\]%
\[
\sum_{n=0}^{\infty }c_{n}\frac{\left\{ -816+792\cdot 2^{1/3}+16\cdot
2^{2/3}+(-1296+3348\cdot 2^{1/3})n+4356\cdot 2^{1/3}n^{2}\right\} }{%
(4+30\cdot 2^{1/3})^{n}} 
\]%
\[
=\frac{4050+10137\cdot 2^{1/3}+2702\cdot 2^{2/3}}{\pi ^{2}} 
\]

\textbf{10.10.}%
\[
J_{10}(-\exp (-\pi \sqrt{3}))=4-160\cdot 3^{1/3} 
\]%
\[
\sum_{n=0}^{\infty }c_{n}\frac{\left\{ 19494+22\cdot 3^{1/3}+2158\cdot
3^{2/3}+(105921+3918\cdot 3^{1/3})n+192027n^{2}\right\} }{(4-160\cdot
3^{1/3})^{n}} 
\]%
\[
=\frac{191999+120\cdot 3^{1/3}-4800\cdot 3^{2/3}}{\pi ^{2}} 
\]

\textbf{10.11.}%
\[
J_{10}(\exp (\frac{\pi i}{3}-\pi \sqrt{\frac{7}{9}}))=\frac{1}{2}(23+15i%
\sqrt{3}) 
\]%
\[
\sum_{n=0}^{\infty }A_{n}\frac{82i+46\sqrt{3}+(207i+135\sqrt{3})n}{(\frac{1}{%
2}(23+15i\sqrt{3}))^{n}}=\frac{86\sqrt{23+15i\sqrt{3}}}{\pi \sqrt{2}} 
\]

\textbf{10.12}%
\[
J_{10}(\exp (\frac{\pi i}{3}-\pi \sqrt{\frac{11}{9}}))=20+16i\sqrt{3} 
\]%
\[
\sum_{n=0}^{\infty }A_{n}\frac{239i+162\sqrt{3}+(770i+616\sqrt{3})n}{(20+16i%
\sqrt{3})^{n}}=\frac{438\sqrt{5+4i\sqrt{3}}}{\pi } 
\]

\textbf{10.13}%
\[
J_{10}(\exp (\frac{\pi i}{3}-\pi \sqrt{\frac{19}{9}}))=52+48i\sqrt{3} 
\]%
\[
\sum_{n=0}^{\infty }A_{n}\frac{1021i+850\sqrt{3}+(4446i+4104\sqrt{3})n}{%
(52+48i\sqrt{3})^{n}}=\frac{1202\sqrt{13+12i\sqrt{3}}}{\pi } 
\]

\textbf{10.14}%
\[
J_{10}(\exp (\frac{\pi i}{3}-\pi \sqrt{\frac{43}{9}}))=484+480i\sqrt{3} 
\]%
\[
\sum_{n=0}^{\infty }A_{n}\frac{41325i+40164\sqrt{3}+(280962i+278640\sqrt{3})n%
}{(484+480i\sqrt{3})^{n}}=\frac{16526\sqrt{121+120i\sqrt{3}}}{\pi } 
\]

\textbf{10.15}%
\[
J_{10}(\exp (\frac{\pi i}{3}-\pi \sqrt{\frac{67}{9}}))=2644+2640i\sqrt{3} 
\]%
\[
\sum_{n=0}^{\infty }A_{n}\frac{2890159i+2871466\sqrt{3}+(24712146i+24674760%
\sqrt{3})n}{(2644+2640i\sqrt{3})^{n}}=\frac{498206\sqrt{661+660i\sqrt{3}}}{%
\pi } 
\]

\textbf{10.16}%
\[
J_{10}(\exp (\frac{\pi i}{3}-\pi \sqrt{\frac{163}{9}}))=320164+320160i\sqrt{3%
} 
\]%
\[
\sum_{n=0}^{\infty }A_{n}\frac{466246542929i+466207604348\sqrt{3}%
+(6233365923642i+6433288046480\sqrt{3})n}{(320164+320160i\sqrt{3})^{n}} 
\]%
\[
=\frac{7321647566\sqrt{80041+80040i\sqrt{3}}}{\pi } 
\]

\textbf{10.17}%
\[
J_{10}(\exp (\frac{2\pi i}{3}-\pi \sqrt{\frac{8}{9}}))=-6+10i\sqrt{3} 
\]%
\[
\sum_{n=0}^{\infty }A_{n}\frac{3+7i\sqrt{3}+(12+20i\sqrt{3})n}{(-6+10i\sqrt{3%
})^{n}}=\frac{18\sqrt{-3+5i\sqrt{3}}}{\pi \sqrt{2}} 
\]

\textbf{10.18}%
\[
J_{10}(\exp (\frac{2\pi i}{3}-\pi \sqrt{\frac{16}{9}}))=-29+33i\sqrt{3} 
\]%
\[
\sum_{n=0}^{\infty }A_{n}\frac{800+1052i\sqrt{3}+(3654+4158i\sqrt{3})n}{%
(-29+33i\sqrt{3})^{n}}=\frac{1027\sqrt{-29+33i\sqrt{3}}}{\pi \sqrt{2}} 
\]

\textbf{10.19}%
\[
J_{10}(\exp (\frac{2\pi i}{3}-\pi \sqrt{\frac{28}{9}}))=\frac{1}{2}(-247+255i%
\sqrt{3}) 
\]%
\[
\sum_{n=0}^{\infty }A_{n}\frac{155394+169758i\sqrt{3}+(886977+915705i\sqrt{3}%
)n}{(\frac{1}{2}(-247+255i\sqrt{3}))^{n}}=\frac{64021\sqrt{-247+255i\sqrt{3}}%
}{\pi \sqrt{2}} 
\]

\textbf{K3-surfaces.}

Duco van Straten sent me a list consisting of 24 third order differential
equations coming from K3-surfaces with Picard number 19. They were obtained
from polytopes but the computations no longer existed. The list contained
the four hypergeometric equations as well $(\alpha ),(\beta ),...(\kappa )$
and some others which we now show (partly).

\textbf{Example 11.}

Consider%
\[
A_{n}=\sum_{k=0}^{n}\dbinom{n}{k}^{2}\dbinom{n+k}{n}\dbinom{2k}{n} 
\]%
Then%
\[
y_{0}=\sum_{n=0}^{\infty }A_{n}x^{n}=1+4x+48x^{2}+760x^{3}+13840x^{4}+... 
\]%
satisfies 
\[
\theta ^{3}-x(2\theta +1)(13\theta ^{2}+13\theta +4)-3x^{2}(\theta
+1)(3\theta +2)(3\theta +4) 
\]%
and%
\[
J_{11}=\frac{1}{q}+9+51q+204q^{2}+681q^{3}+1956q^{4}+... 
\]%
which is 7A in [13]. We obtain

\textbf{11.1.}%
\[
J_{11}(\exp (-\pi \sqrt{\frac{4}{7}}))=27 
\]%
\[
\sum_{n=0}^{\infty }A_{n}(6-n)\frac{1}{27^{n}}=\frac{4}{\pi } 
\]

\textbf{11.2.}%
\[
J_{11}(\exp (-\pi \sqrt{\frac{16}{7}}))=125 
\]%
\[
\sum_{n=0}^{\infty }A_{n}(32+168n)\frac{1}{125^{n}}=\frac{125}{\pi } 
\]

\textbf{11.3.}%
\[
J_{11}(-\exp (-\pi \sqrt{\frac{13}{7}}))=-64 
\]%
\[
\sum_{n=0}^{\infty }A_{n}(10+39n)\frac{1}{(-64)^{n}}=\frac{64}{\pi \sqrt{7}} 
\]

\textbf{11.4.}%
\[
J_{11}(-\exp (-\pi \sqrt{\frac{61}{7}}))=-10648 
\]%
\[
\sum_{n=0}^{\infty }A_{n}(1286+39n)\frac{1}{(-10648)^{n}}=\frac{10648}{\pi 
\sqrt{7}} 
\]

\textbf{11.5.}%
\[
J_{11}(-\exp (-\pi \sqrt{\frac{85}{7}}))=-(20+2\sqrt{85})^{3} 
\]%
\[
\sum_{n=0}^{\infty }A_{n}\frac{68\sqrt{119}-106\sqrt{35}+(198\sqrt{119}-153%
\sqrt{35})n}{(-(20+2\sqrt{85}))^{3n}}=\frac{360}{\pi } 
\]

\textbf{11.6.}%
\[
J_{11}(-\exp (-\pi \sqrt{\frac{109}{7}}))=-(31+3\sqrt{109})^{3} 
\]%
\[
\sum_{n=0}^{\infty }A_{n}\frac{308034-23048\sqrt{109}+835485n}{(-(31+3\sqrt{%
109}))^{3n}}=\frac{278495000}{\pi \sqrt{-22060520573+2113180384\sqrt{109}}} 
\]

\textbf{Example 12.}

Consider%
\[
\theta ^{3}-3x(2\theta +1)(9\theta ^{2}+9\theta +4)+81x^{2}(\theta
+1)(3\theta +2)(3\theta +4) 
\]%
with%
\[
J_{12}=\frac{1}{q}+15+54q-76q^{2}-243q^{3}+1188q^{4}+... 
\]%
which is 3B in [13]. After finding a couple of formulas for $\dfrac{1}{\pi }$
I realized that 
\[
A_{n}=\sum_{k=0}^{n}\dbinom{2k}{k}\dbinom{3k}{k}\dbinom{2n-2k}{n-k}\dbinom{%
3n-3k}{n-k} 
\]%
which is treated in Z.W.Sun [20]. Of the numerous square root values of \ $%
J_{12}$ we choose one

\textbf{12.1.}%
\[
J_{12}(-\exp (-\pi \sqrt{\frac{89}{3}}))=-27000(500+53\sqrt{89}) 
\]%
\[
\sum_{n=0}^{\infty }A_{n}\frac{1654\sqrt{89}+(133500+14151\sqrt{89})n}{%
(-27000(500+53\sqrt{89}))^{n}}=\frac{9000}{\pi }\sqrt{\frac{89}{3}} 
\]

\textbf{12.2.}

\[
J_{12}(-\exp (-\pi \sqrt{12}))=3^{3}(4+3\cdot 2^{1/3}+3\cdot 2^{2/3})^{3} 
\]%
\[
\sum_{n=0}^{\infty }A_{n}\frac{12-6\cdot 2^{2/3}+27n}{(3(4+3\cdot
2^{1/3}+3\cdot 2^{2/3}))^{3n}}=\frac{(1-2^{1/3}+3\cdot 2^{2/3})\sqrt{3}}{\pi 
} 
\]

\textbf{12.3.}%
\[
J_{12}(-\exp (-\pi \sqrt{27}))=-6^{3}(13+9\cdot 3^{1/3}+6\cdot 3^{2/3})^{3} 
\]%
\[
\sum_{n=0}^{\infty }A_{n}\frac{108-108\cdot 3^{1/3}+18\cdot 3^{2/3}-81\cdot
3^{2/3}n}{(-6(13+9\cdot 3^{1/3}+6\cdot 3^{2/3}))^{3n}}=\frac{(-56+72\cdot
3^{1/3}-32\cdot 3^{2/3})\sqrt{3}}{\pi } 
\]

\textbf{12.4.}%
\[
J_{12}(\exp (\frac{\pi i}{3}-\pi \sqrt{\frac{7}{4}}))=32+8i\sqrt{11} 
\]%
\[
\sum_{n=0}^{\infty }A_{n}\frac{6i+(11i+4\sqrt{11})n}{(32+8i\sqrt{11}))^{n}}=%
\frac{24}{\pi } 
\]

\textbf{Example 13.}

Consider%
\[
\theta ^{3}-24x(2\theta +1)(18\theta ^{2}+18\theta +5)+20736x^{2}(\theta
+1)(3\theta +1)(3\theta +5) 
\]%
with%
\[
J_{13}=\frac{1}{q}+312+10260q+901120q^{2}+91676610q^{3}+... 
\]%
which is not in [13]. We have%
\[
A_{n}=\sum_{k=0}^{n}\dbinom{3k}{k}\dbinom{6k}{3k}\dbinom{3n-3k}{n-k}\dbinom{%
6n-6k}{3n-3k} 
\]%
((also in Sun [20]). We show only

\textbf{13.1.}%
\[
J_{13}(-\exp (-\pi \sqrt{67}))=-69120(1064800+7161\sqrt{22110}) 
\]%
\[
\sum_{n=0}^{\infty }A_{n}\frac{10177\sqrt{22110}+(19456800+130851\sqrt{22110}%
)n}{(-69120(1064800+7161\sqrt{22110}))^{n}}=\frac{580800}{\pi }\sqrt{67} 
\]

\textbf{Example 14.}

Consider%
\[
\theta ^{3}-5x(2\theta +1)(26\theta ^{2}+26\theta +5)+x^{2}(\theta
+1)(774\theta ^{2}+1548\theta +823) 
\]%
\[
-386x^{3}(\theta +1)(\theta +2)(2\theta +3)+257x^{4}(\theta +1)(\theta
+2)(\theta +3) 
\]%
with%
\[
J_{14}=\frac{1}{q}+105+4372q+96256q^{2}+1240002q^{3}+10698752q^{4}+... 
\]%
which we identify as 2A but also as \ $J_{1/4}+1$. We get%
\[
A_{n}=\sum_{k=0}^{n}\dbinom{n}{k}\dbinom{2k}{k}^{2}\dbinom{4k}{2k} 
\]%
and arguments like%
\[
\frac{1}{2459125787}=\frac{1}{2^{8}99^{4}+1} 
\]%
Three further equations result in \ $J_{1/2}+1,J_{1/3}+1,J_{1/6}+1$ \
respectively.

\textbf{14.1.}%
\[
J_{14}(i\exp (-\pi \sqrt{\frac{15}{4}}))=\frac{1}{2}(209+495i\sqrt{3}) 
\]%
\[
\sum_{n=0}^{\infty }A_{n}\frac{659031048+1590088778i\sqrt{3}%
+(4081716765+712863585i\sqrt{3})n}{(\frac{1}{2}(209+495i\sqrt{3}))^{n}}=%
\frac{22^{2}1609^{2}}{\pi }\sqrt{3} 
\]

\textbf{Example 15.}

Consider%
\[
(625\theta ^{3}-250x(2\theta +1)(17\theta ^{2}+17\theta +6)-1400x^{2}(\theta
+1)(11\theta ^{2}+22\theta +12) 
\]%
\[
-630x^{3}(\theta +1)(\theta +2)(2\theta +3)-1504x^{4}(\theta +1)(\theta
+2)(\theta +3) 
\]%
with%
\[
J_{15}=\frac{1}{q}+\frac{22}{5}+17q+46q^{2}+116q^{3}+252q^{4}+... 
\]%
which is 11A.We have%
\[
y_{0}=1+\frac{12}{5}x+\frac{444}{25}x^{2}+\frac{20028}{125}x^{3}+\frac{%
1037436}{625}x^{4}+... 
\]

\textbf{15.1.}%
\[
J_{15}(\exp (-\pi \sqrt{\frac{8}{11}}))=\frac{102}{5} 
\]%
\[
\sum_{n=0}^{\infty }A_{n}(19+70n)\frac{1}{(\dfrac{102}{5})^{n}}=\frac{51^{2}%
}{5\pi \sqrt{11}} 
\]

\textbf{15.2.}%
\[
J_{15}(-\exp (-\pi \sqrt{\frac{17}{11}}))=-\frac{228}{5} 
\]%
\[
\sum_{n=0}^{\infty }A_{n}(307+1105n)\frac{1}{(-\dfrac{228}{5})^{n}}=\frac{%
12996}{5\pi \sqrt{11}} 
\]

\textbf{15.3}

\[
J_{15}(-\exp (-\pi \sqrt{\frac{137}{11}}))=-\frac{163278}{5}-2790\sqrt{137} 
\]%
\[
\sum_{n=0}^{\infty }A_{n}\frac{187348870-13008519\sqrt{137}+388888885n}{(-%
\dfrac{163278}{5}-2790\sqrt{137})^{n}}=\frac{2^{2}3^{4}17^{2}19^{2}7070707}{%
\pi \sqrt{-15807661556999+1366619918190\sqrt{137}}} 
\]

\textbf{Example 16.}

This example is obtained by scaling off \ $\dbinom{2n}{n}$ \ from \# 183 in
the Big Table [3] of Calabi-Yau differential equations%
\[
A_{n}=3\sum_{k=0}^{n}(-1)^{k}\frac{n-2k}{2n-3k}\dbinom{n}{k}\dbinom{2k}{k}%
\dbinom{2n-2k}{n-k}\dbinom{2n-3k}{n} 
\]%
\[
\theta ^{3}-2x(2\theta +1)(7\theta ^{2}+7\theta +3)+12x^{2}(\theta
+1)(4\theta +3)(4\theta +5) 
\]%
\[
J_{16}=\frac{1}{q}%
+8+16q-8q^{2}+128q^{4}+28q^{5}+576q^{7}-64q^{8}+2048q^{10}+... 
\]%
which seems to be missing from [13]. Observe that $J_{16}$ does not have any
terms of degree $3n.$

\textbf{16.1.}%
\[
J_{16}(\exp (-\pi \sqrt{\frac{10}{9}}))=36 
\]%
\[
\sum_{n=0}^{\infty }A_{n}(3+20n)\frac{1}{(36)^{n}}=\frac{18\sqrt{3}}{\pi } 
\]

\textbf{16.2.}%
\[
J_{16}(\exp (-\pi \sqrt{\frac{22}{9}}))=144 
\]%
\[
\sum_{n=0}^{\infty }A_{n}(15+88n)\frac{1}{(144)^{n}}=\frac{36\sqrt{3}}{\pi } 
\]

\textbf{16.3.}%
\[
J_{16}(\exp (-\pi \sqrt{\frac{58}{9}}))=2916 
\]%
\[
\sum_{n=0}^{\infty }A_{n}(789+6380n)\frac{1}{(2916)^{n}}=\frac{1458\sqrt{3}}{%
\pi } 
\]

\textbf{16.4.}%
\[
J_{16}(-\exp (-\pi \sqrt{\frac{13}{9}}))=-36 
\]%
\[
\sum_{n=0}^{\infty }A_{n}(9+26n)\frac{1}{(-36)^{n}}=\frac{9\sqrt{3}}{\pi } 
\]

\textbf{16.5.}%
\[
J_{16}(-\exp (-\pi \sqrt{\frac{25}{9}}))=-180 
\]%
\[
\sum_{n=0}^{\infty }A_{n}(6+28n)\frac{1}{(-180)^{n}}=\frac{9\sqrt{3}}{\pi } 
\]

\textbf{16.6.}%
\[
J_{16}(-\exp (-\pi \sqrt{\frac{37}{9}}))=-576 
\]%
\[
\sum_{n=0}^{\infty }A_{n}(171+1036n)\frac{1}{(-576)^{n}}=\frac{288\sqrt{3}}{%
\pi } 
\]

\textbf{16.7.}%
\[
J_{16}(\exp (\frac{\pi i}{3}-\pi \sqrt{\frac{7}{9}}))=\frac{1}{2}(33+15i%
\sqrt{3}) 
\]%
\[
\sum_{n=0}^{\infty }A_{n}\frac{66-2i\sqrt{3}+(143-65i\sqrt{3})n}{(\dfrac{1}{2%
}(33+15i\sqrt{3}))^{n}}=\frac{546\sqrt{42}}{\pi \sqrt{-111+25i\sqrt{3}}} 
\]

\textbf{16.8}%
\[
J_{16}(\exp (\frac{\pi i}{3}-\pi \sqrt{\frac{13}{9}}))=30+6i\sqrt{39} 
\]%
\[
\sum_{n=0}^{\infty }A_{n}\frac{117+7i\sqrt{39}+520n}{(30+6i\sqrt{39})^{n}}=%
\frac{74880}{\pi \sqrt{18486+1794i\sqrt{39}}} 
\]

\textbf{16.9}%
\[
J_{16}(\exp (\frac{\pi i}{3}-\pi \sqrt{\frac{25}{9}}))=102+42i\sqrt{15} 
\]%
\[
\sum_{n=0}^{\infty }A_{n}\frac{165+5i\sqrt{15}+920n}{(102+42i\sqrt{15})^{n}}=%
\frac{423936}{\pi \sqrt{542790+18102i\sqrt{15}}} 
\]

\textbf{16.10.}%
\[
J_{16}(\exp (\frac{\pi i}{3}-\pi \sqrt{\frac{37}{9}}))=300+48i\sqrt{111} 
\]%
\[
\sum_{n=0}^{\infty }A_{n}\frac{22977+101i\sqrt{111}+150220n}{(300+48i\sqrt{%
111})^{n}}=\frac{1623051990}{\pi \sqrt{468132066+1843377i\sqrt{111}}} 
\]

\textbf{16.11.}%
\[
J_{16}(\exp (\frac{2\pi i}{3}-\pi \sqrt{\frac{10}{9}}))=-6+6i\sqrt{15} 
\]%
\[
\sum_{n=0}^{\infty }A_{n}\frac{15+i\sqrt{15}+40n}{(-6+6i\sqrt{15})^{n}}=%
\frac{360}{\pi \sqrt{90+30i\sqrt{15}}} 
\]

\textbf{16.12.}%
\[
J_{16}(\exp (\frac{2\pi i}{3}-\pi \sqrt{\frac{22}{9}}))=-60+48i\sqrt{6} 
\]%
\[
\sum_{n=0}^{\infty }A_{n}\frac{339+17i\sqrt{6}+1540n}{(-60+48i\sqrt{6})^{n}}=%
\frac{838530}{\pi \sqrt{789921+58542i\sqrt{6}}} 
\]

\textbf{16.13}%
\[
J_{16}(\exp (\frac{2\pi i}{3}-\pi \sqrt{\frac{58}{9}}))=-1446+270i\sqrt{87} 
\]%
\[
\sum_{n=0}^{\infty }A_{n}\frac{29319+37i\sqrt{87}+232232n}{(-1446+270i\sqrt{%
87})^{n}}=\frac{10200326136}{\pi \sqrt{12492362538+11171670i\sqrt{87}}} 
\]%
\textbf{.}

\textbf{Mirror maps.}

In B.H.Lian-S.T.Yau [14] we found a table with some third order differential
equations.

\textbf{Example 17.}

Consider 
\[
\theta ^{3}-6x(2\theta +1)(2\theta ^{2}+2\theta +1)+4x^{2}(\theta
+1)(44\theta ^{2}+88\theta +51)-12x^{3}(\theta +1)(\theta +2)(2\theta +3) 
\]%
with%
\[
J_{17}=\frac{1}{q}+6+7q+15q^{2}+71q^{3}+106q^{4}+... 
\]%
which is 12C. We find

\textbf{17.1.}%
\[
J_{17}(\exp (-\pi \sqrt{\frac{3}{6}}))=16 
\]%
\[
\sum_{n=0}^{\infty }A_{n}n\frac{1}{(16)^{n}}=\frac{8}{\pi \sqrt{3}} 
\]

\textbf{17.2.}%
\[
J_{17}(\exp (-\pi \sqrt{\frac{5}{6}}))=24 
\]%
\[
\sum_{n=0}^{\infty }A_{n}(1+10n)\frac{1}{(24)^{n}}=\frac{12\sqrt{3}}{\pi } 
\]

\textbf{17.3.}%
\[
J_{17}(\exp (-\pi \sqrt{\frac{7}{6}}))=36 
\]%
\[
\sum_{n=0}^{\infty }A_{n}(2+14n)\frac{1}{(36)^{n}}=\frac{27}{\pi \sqrt{2}} 
\]

\textbf{17.4.}%
\[
J_{17}(\exp (-\pi \sqrt{\frac{13}{6}}))=108 
\]%
\[
\sum_{n=0}^{\infty }A_{n}(22+130n)\frac{1}{(108)^{n}}=\frac{81}{\pi }\sqrt{%
\frac{3}{2}} 
\]

\textbf{17.5.}%
\[
J_{17}(\exp (-\pi \sqrt{\frac{17}{6}}))=204 
\]%
\[
\sum_{n=0}^{\infty }A_{n}(23+40n)\frac{1}{(204)^{n}}=\frac{51\sqrt{3}}{\pi } 
\]

\textbf{17.6.}%
\[
J_{17}(-\exp (-\pi \sqrt{\frac{7}{6}}))=-24 
\]%
\[
\sum_{n=0}^{\infty }A_{n}(3+7n)\frac{1}{(-24)^{n}}=\frac{3\sqrt{2}}{\pi } 
\]

\textbf{17.7.}%
\[
J_{17}(-\exp (-\pi \sqrt{\frac{13}{6}}))=-96 
\]%
\[
\sum_{n=0}^{\infty }A_{n}(17+65n)\frac{1}{(-96)^{n}}=\frac{16\sqrt{6}}{\pi } 
\]

\textbf{17.8.}%
\[
J_{17}(-\exp (-\pi \sqrt{\frac{17}{6}}))=-192 
\]%
\[
\sum_{n=0}^{\infty }A_{n}(127+595n)\frac{1}{(-192)^{n}}=\frac{192\sqrt{3}}{%
\pi } 
\]

\bigskip \textbf{17.9.}%
\[
J_{17}(i\exp (-\pi \sqrt{\frac{11}{12}}))=6+6i\sqrt{11} 
\]%
\[
\sum_{n=0}^{\infty }A_{n}\frac{71i+11\sqrt{11}+(110i+50\sqrt{11})n}{(6+6i%
\sqrt{11})^{n}}=\frac{6^{3}}{\pi } 
\]

\textbf{17.10}%
\[
J_{17}(i\exp (-\pi \sqrt{\frac{19}{12}}))=6+30i\sqrt{3} 
\]%
\[
\sum_{n=0}^{\infty }A_{n}\frac{243i+227\sqrt{3}+(390i+962\sqrt{3})n}{(6+30i%
\sqrt{3})^{n}}=\frac{1368}{\pi } 
\]

\textbf{17.11.}%
\[
J_{17}(i\exp (-\pi \sqrt{\frac{31}{12}}))=6+90i\sqrt{3} 
\]%
\[
\sum_{n=0}^{\infty }A_{n}\frac{11691i+28673\sqrt{3}+(19530i+146258\sqrt{3})n%
}{(6+90i\sqrt{3})^{n}}=\frac{158184}{\pi } 
\]

\textbf{17.12}%
\[
J_{17}(i\exp (-\pi \sqrt{\frac{59}{12}}))=6+138i\sqrt{59} 
\]%
\[
\sum_{n=0}^{\infty }A_{n}\frac{411221i+1186903\sqrt{59}+(719210i+8270650%
\sqrt{59})n}{(6+138i\sqrt{59})^{n}}=\frac{306^{3}}{\pi } 
\]%
\textbf{.}

There are 22 real square root values of \ $J_{17}$

\textbf{Example 18.}

Consider%
\[
\theta ^{3}-x(2\theta +1)(11\theta ^{2}+11\theta +5)+x^{2}(\theta
+1)(121\theta ^{2}+242\theta +141)-98x^{3}(\theta +1)(\theta +2)(2\theta +3) 
\]%
with%
\[
J_{18}=\frac{1}{q}+6+11q+20q^{2}+57q^{3}+92q^{4}+... 
\]%
which is 14A. We find

\textbf{18.1.}%
\[
J_{18}(\exp (-\pi \sqrt{\frac{4}{7}}))=18 
\]%
\[
\sum_{n=0}^{\infty }A_{n}(1+7n)\frac{1}{(18)^{n}}=\frac{27}{\pi } 
\]

\textbf{18.2}%
\[
J_{18}(\exp (-\pi \sqrt{\frac{6}{7}}))=25 
\]%
\[
\sum_{n=0}^{\infty }A_{n}(4+24n)\frac{1}{(25)^{n}}=\frac{125}{\pi \sqrt{7}} 
\]

\textbf{18.3.}%
\[
J_{18}(\exp (-\pi \sqrt{\frac{10}{7}}))=49 
\]%
\[
\sum_{n=0}^{\infty }A_{n}(22+120n)\frac{1}{(49)^{n}}=\frac{49\sqrt{7}}{\pi } 
\]

\textbf{18.4.}%
\[
J_{18}(-\exp (-\pi \sqrt{\frac{19}{7}}))=-171 
\]%
\[
\sum_{n=0}^{\infty }A_{n}(73+340n)\frac{1}{(-171)^{n}}=\frac{513}{\pi \sqrt{7%
}} 
\]%
There are 20 square root values of \ $J_{18}.$

\textbf{Example 19.}

Consider%
\[
(\theta ^{3}-x(2\theta +1)(7\theta ^{2}+7\theta +4)+x^{2}(\theta
+1)(61\theta ^{2}+122\theta +88) 
\]%
\[
-8x^{3}(2\theta +3)(7\theta ^{2}+21\theta +18)+64x^{4}(\theta +2)^{3} 
\]%
with%
\[
J_{19}=\frac{1}{q}+3+3q+4q^{2}+9q^{3}+12q^{4}+... 
\]%
which is 14B. We show only one formula

\textbf{19.1.}%
\[
J_{19}(\exp (-\pi \sqrt{\frac{19}{7}}))=(3+\sqrt{7})^{3} 
\]%
\[
\sum_{n=0}^{\infty }A_{n}\frac{23416-8850\sqrt{7}+(12825-48245\sqrt{7})n}{(3+%
\sqrt{7})^{3n}}=\frac{4}{\pi } 
\]

\textbf{Example 20.}

Consider%
\[
\theta ^{3}-x(2\theta +1)(7\theta ^{2}+7\theta +3)+x^{2}(\theta +1)(29\theta
^{2}+58\theta +33)-30x^{3}(\theta +1)(\theta +2)(2\theta +3) 
\]%
with%
\[
J_{20}=\frac{1}{q}+4+8q+22q^{2}+42q^{3}+70q^{4}+... 
\]%
which is 15A. We obtain

\textbf{20.1.}%
\[
J_{20}(\exp (-\pi \sqrt{\frac{8}{15}}))=15 
\]%
\[
\sum_{n=0}^{\infty }A_{n}(2+8n)\frac{1}{(15)^{n}}=\frac{15\sqrt{3}}{\pi } 
\]

\textbf{20.2.}%
\[
J_{20}(\exp (-\pi \sqrt{\frac{16}{15}}))=30 
\]%
\[
\sum_{n=0}^{\infty }A_{n}(6+26n)\frac{1}{(30)^{n}}=\frac{15\sqrt{5}}{\pi } 
\]

\textbf{20.3.}%
\[
J_{20}(-\exp (-\pi \sqrt{\frac{13}{15}}))=-15 
\]%
\[
\sum_{n=0}^{\infty }A_{n}(11+26n)\frac{1}{(-15)^{n}}=\frac{5\sqrt{15}}{\pi } 
\]

\textbf{20.4.}%
\[
J_{20}(-\exp (-\pi \sqrt{\frac{29}{15}}))=-75 
\]%
\[
\sum_{n=0}^{\infty }A_{n}(251+986n)\frac{1}{(-75)^{n}}=\frac{375\sqrt{3}}{%
\pi } 
\]

\textbf{20.5.}%
\[
J_{20}(-\exp (-\pi \sqrt{\frac{37}{15}}))=-135 
\]%
\[
\sum_{n=0}^{\infty }A_{n}(113+518n)\frac{1}{(-135)^{n}}=\frac{81\sqrt{15}}{%
\pi } 
\]

\textbf{20.6.}

\[
J_{20}(-\exp (-\pi \sqrt{\frac{53}{15}}))=-363 
\]

\[
\sum_{n=0}^{\infty }A_{n}(2327+13250n)\frac{1}{(-363)^{n}}=\frac{3993\sqrt{3}%
}{\pi } 
\]

\textbf{20.7.}%
\[
J_{20}(\exp (\frac{\pi i}{3}-\pi \sqrt{\frac{31}{45}}))=11+2i\sqrt{31} 
\]%
\[
\sum_{n=0}^{\infty }A_{n}\frac{2697+144i\sqrt{31}+7750n}{(11+2i\sqrt{31})^{n}%
}=\frac{643125\sqrt{155}}{\pi \sqrt{280331+69838i\sqrt{31}}} 
\]%
There are 17 real square root values of $J_{20.}$

\textbf{Example 21.}

Consider%
\[
\theta ^{3}-x(2\theta +1)(7\theta ^{2}+7\theta +4)+x^{2}(\theta +1)(53\theta
^{2}+106\theta +72) 
\]%
\[
-14x^{3}(\theta +1)(\theta +2)(2\theta +3)-24x^{4}(\theta +2)(2\theta
+3)(2\theta +5) 
\]%
with%
\[
J_{21}=\frac{1}{q}+3+4q+2q^{2}+6q^{3}+10q^{4}+... 
\]%
which is 30B. We get

\textbf{21.1}%
\[
J_{21}(\exp (-\pi \sqrt{\frac{4}{15}}))=9 
\]%
\[
\sum_{n=0}^{\infty }A_{n}n\frac{1}{(9)^{n}}=\frac{81}{\pi \sqrt{5}} 
\]

\textbf{21.2.}%
\[
J_{21}(\exp (-\pi \sqrt{\frac{14}{15}}))=24 
\]%
\[
\sum_{n=0}^{\infty }A_{n}(6+35n)\frac{1}{(24)^{n}}=\frac{36\sqrt{2}}{\pi } 
\]

\textbf{21.3.}%
\[
J_{21}(\exp (-\pi \sqrt{\frac{22}{15}}))=48 
\]%
\[
\sum_{n=0}^{\infty }A_{n}(29+154n)\frac{1}{(48)^{n}}=\frac{192}{\pi }\sqrt{%
\frac{3}{5}} 
\]

\textbf{21.4.}%
\[
J_{21}(-\exp (-\pi \sqrt{\frac{11}{15}}))=-12 
\]%
\[
\sum_{n=0}^{\infty }A_{n}(6+11n)\frac{1}{(-12)^{n}}=\frac{18}{\pi \sqrt{5}} 
\]

\textbf{21.5.}%
\[
J_{21}(-\exp (-\pi \sqrt{\frac{23}{15}}))=-46 
\]%
\[
\sum_{n=0}^{\infty }A_{n}(101+315n)\frac{1}{(-46)^{n}}=\frac{46\sqrt{23}}{%
\pi } 
\]%
There are 17 square root values of $J_{21}.$

\textbf{Example 22.}

Consider 
\[
\theta ^{3}-4x(2\theta +1)(2\theta ^{2}+2\theta +1)+16x^{2}(\theta
+1)(4\theta ^{2}+8\theta +5)-8x^{3}(2\theta +3)^{3} 
\]%
with%
\[
J_{22}=\frac{1}{q}+4+6q+8q^{2}+17q^{3}+32q^{4}+... 
\]%
which is 20A. We get

\textbf{22.1.}%
\[
J_{22}(\exp (-\pi \sqrt{\frac{3}{5}}))=16 
\]%
\[
\sum_{n=0}^{\infty }A_{n}(1+6n)\frac{1}{(16)^{n}}=\frac{16}{\pi } 
\]

\textbf{22.2.}%
\[
J_{22}(-\exp (-\pi \sqrt{\frac{38}{5}}))=-2884-912\sqrt{10} 
\]%
\[
\sum_{n=0}^{\infty }A_{n}\frac{-209203+66156\sqrt{10}+(-151164+47804\sqrt{10}%
)n}{(-(2884+912\sqrt{10}))^{n}}=\frac{2}{\pi } 
\]%
There are 10 more square root values of $J_{22}.$

\textbf{Example 23.}

Consider%
\[
\theta ^{3}-8x(2\theta +1)(8\theta ^{2}+8\theta +3)+1024x^{2}(\theta
+1)(2\theta +1)(2\theta +3) 
\]%
with%
\[
J_{23}=\frac{1}{q}+40+276q-2048q^{2}+11202q^{3}-49152q^{4}+... 
\]%
which is 2B and is also recognized as the $\ J-$function of \ $\dbinom{2n}{n}%
\ast (e).$Hence%
\[
A_{n}=\dbinom{2n}{n}\sum_{k=0}^{n}4^{n-k}\dbinom{2k}{k}^{2}\dbinom{2n-2k}{n-k%
} 
\]%
We find

\textbf{23.1}%
\[
J_{23}(\exp (-\pi ))=72 
\]%
\[
\sum_{n=0}^{\infty }A_{n}(1+n)\frac{1}{(72)^{n}}=\frac{9}{\pi } 
\]

\textbf{23.2.}%
\[
J_{23}(\exp (-\pi \sqrt{2}))=128 
\]%
\[
\sum_{n=0}^{\infty }A_{n}n\frac{1}{(128)^{n}}=\frac{\sqrt{2}}{\pi } 
\]

\textbf{23.3.}%
\[
J_{23}(\exp (-2\pi ))=576 
\]%
\[
\sum_{n=0}^{\infty }A_{n}(1+8n)\frac{1}{(576)^{n}}=\frac{9}{2\pi } 
\]

\textbf{23.4.}%
\[
J_{23}(-\exp (-\pi \sqrt{3}))=-192 
\]%
\[
\sum_{n=0}^{\infty }A_{n}(1+4n)\frac{1}{(-192)^{n}}=\frac{\sqrt{3}}{\pi } 
\]

\textbf{23.5.}%
\[
J_{23}(-\exp (-\pi \sqrt{6}))=-4032 
\]%
\[
\sum_{n=0}^{\infty }A_{n}(1+8n)\frac{1}{(-4032)^{n}}=\frac{9\sqrt{7}}{8\pi } 
\]

\textbf{23.6}%
\[
J_{23}(i\exp (-\pi \sqrt{\frac{7}{4}}))=\frac{1}{2}(81+45i\sqrt{7}) 
\]%
\[
\sum_{n=0}^{\infty }A_{n}\frac{16i+(35i+9\sqrt{7})n}{(\frac{1}{2}(81+45i%
\sqrt{7}))^{n}}=\frac{36}{\pi } 
\]%
\textbf{.}

While searching for Calabi-Yau differential equations there was some
spin-off which we show below

\textbf{Example 24.}

Let%
\[
A_{n}=\sum_{k=0}^{n}\dbinom{n}{2k}\dbinom{2k}{k}^{2}\dbinom{2n-4k}{n-2k} 
\]%
with%
\[
y_{0}=\sum_{n=0}^{\infty
}A_{n}x^{n}=1+2x+10x^{2}+44x^{3}+250x^{3}+1412x^{4}+... 
\]%
satisfying%
\[
\theta ^{3}-2x(2\theta +1)(2\theta ^{2}+2\theta +1)-4x^{2}(\theta
+1)(4\theta ^{2}+8\theta +5)+64x^{3}(\theta +1)(\theta +2)(2\theta +3) 
\]%
Then%
\[
J_{24}=\frac{1}{q}+2+7q+15q^{3}+71q^{5}+106q^{7}+... 
\]%
which is 12C.

\textbf{24.1.}%
\[
J_{24}(\exp (-\pi \sqrt{\frac{1}{2}}))=12 
\]%
\[
\sum_{n=0}^{\infty }A_{n}(1+4n)\frac{1}{(12)^{n}}=\frac{6\sqrt{3}}{\pi } 
\]

\textbf{24.2.}%
\[
J_{24}(\exp (-\pi \sqrt{\frac{5}{6}}))=20 
\]%
\[
\sum_{n=0}^{\infty }A_{n}(1+4n)\frac{1}{(20)^{n}}=\frac{10\sqrt{3}}{3\pi } 
\]

\textbf{24.3.}%
\[
J_{24}(\exp (-\pi \sqrt{\frac{7}{6}}))=32 
\]%
\[
\sum_{n=0}^{\infty }A_{n}(5+21n)\frac{1}{(32)^{n}}=\frac{16\sqrt{2}}{\pi } 
\]

\textbf{24.4.}%
\[
J_{24}(\exp (-\pi \sqrt{\frac{13}{6}}))=104 
\]%
\[
\sum_{n=0}^{\infty }A_{n}(1+5n)\frac{1}{(104)^{n}}=\frac{13\sqrt{6}}{9\pi } 
\]

\textbf{24.5.}%
\[
J_{24}(\exp (-\pi \sqrt{\frac{17}{6}}))=200 
\]%
\[
\sum_{n=0}^{\infty }A_{n}(43+238n)\frac{1}{(200)^{n}}=\frac{250\sqrt{3}}{%
3\pi } 
\]

\textbf{24.6.}%
\[
J_{24}(-\exp (-\pi \sqrt{\frac{5}{6}}))=-16 
\]%
\[
\sum_{n=0}^{\infty }A_{n}(2+5n)\frac{1}{(-16)^{n}}=\frac{8\sqrt{3}}{3\pi } 
\]

\textbf{24.7.}%
\[
J_{24}(-\exp (-\pi \sqrt{\frac{7}{6}}))=-28 
\]%
\[
\sum_{n=0}^{\infty }A_{n}(1+3n)\frac{1}{(-28)^{n}}=\frac{7\sqrt{2}}{4\pi } 
\]

\textbf{24.8.}%
\[
J_{24}(-\exp (-\pi \sqrt{\frac{13}{6}}))=-100 
\]%
\[
\sum_{n=0}^{\infty }A_{n}(3+13n)\frac{1}{(-100)^{n}}=\frac{125\sqrt{6}}{%
36\pi } 
\]

\textbf{24.9.}%
\[
J_{24}(-\exp (-\pi \sqrt{\frac{17}{6}}))=-196 
\]%
\[
\sum_{n=0}^{\infty }A_{n}(67+340n)\frac{1}{(-196)^{n}}=\frac{343\sqrt{3}}{%
3\pi } 
\]

\textbf{24.10.}%
\[
J_{24}(\exp (-\pi \sqrt{\frac{11}{6}}))=38+6\sqrt{33} 
\]%
\[
\sum_{n=0}^{\infty }A_{n}\frac{123-11\sqrt{33}+(231+9\sqrt{33})n}{(38+9\sqrt{%
33})^{n}}=\frac{128\sqrt{3}}{\pi } 
\]

\textbf{24.11.}%
\[
J_{24}(\exp (-\pi \sqrt{\frac{83}{6}}))=(99\sqrt{3}+19\sqrt{83})^{2} 
\]%
\[
\sum_{n=0}^{\infty }A_{n}\frac{-23658+1856\sqrt{249}+(6723+3743\sqrt{249})n}{%
(99\sqrt{3}+19\sqrt{83})^{2n}}=\frac{7614\sqrt{3}+494\sqrt{83}}{\pi } 
\]

\textbf{24.12.}%
\[
J_{24}(-\exp (-\pi \sqrt{\frac{85}{6}}))=-59362-3762\sqrt{249} 
\]%
\[
\sum_{n=0}^{\infty }A_{n}\frac{-4929+655\sqrt{249}+(23904+2488\sqrt{249})n}{%
(-59362-3762\sqrt{249})^{n}}=\frac{5097\sqrt{3}+895\sqrt{83}}{\pi } 
\]

\textbf{24.13}%
\[
J_{24}(i\exp (-\pi \sqrt{\frac{7}{12}}))=2+6i\sqrt{3} 
\]%
\[
\sum_{n=0}^{\infty }A_{n}\frac{13i+11\sqrt{3}+(18i+26\sqrt{3})n}{(2+6i\sqrt{3%
})^{n}}=\frac{56}{\pi } 
\]

\textbf{24.14.}%
\[
J_{24}(i\exp (-\pi \sqrt{\frac{11}{12}}))=2+6i\sqrt{11} 
\]%
\[
\sum_{n=0}^{\infty }A_{n}\frac{133i+99\sqrt{11}+(198i+294\sqrt{11})n}{(2+6i%
\sqrt{11})^{n}}=\frac{1000}{\pi } 
\]

\textbf{24.15}%
\[
J_{24}(i\exp (-\pi \sqrt{\frac{19}{12}}))=2+30i\sqrt{3} 
\]%
\[
\sum_{n=0}^{\infty }A_{n}\frac{1073i+3251\sqrt{3}+(1710i+1280\sqrt{3})n}{%
(2+30i\sqrt{3})^{n}}=\frac{26^{3}}{\pi } 
\]

\textbf{24.16}%
\[
J_{24}(i\exp (-\pi \sqrt{\frac{31}{12}}))=2+90i\sqrt{3} 
\]%
\[
\sum_{n=0}^{\infty }A_{n}\frac{2103i+15647\sqrt{3}+(3510i+78962\sqrt{3})n}{%
(2+90i\sqrt{3})^{n}}=\frac{31\cdot 14^{3}}{\pi } 
\]

\textbf{24.17.}%
\[
J_{24}(i\exp (-\pi \sqrt{\frac{59}{12}}))=2+138i\sqrt{59} 
\]%
\[
\sum_{n=0}^{\infty }A_{n}\frac{712281i+6169689\sqrt{59}+(1245726i+42977394%
\sqrt{59})n}{(2+138i\sqrt{59})^{n}}=\frac{530^{3}}{\pi } 
\]

\textbf{Example 25. }

Let%
\[
A_{n}=\dbinom{2n}{n}\sum_{k=0}^{n}\dbinom{n}{2k}\dbinom{2k}{k}\dbinom{2n-4k}{%
n-2k} 
\]%
with

\[
y_{0}=\sum_{n=0}^{\infty
}A_{n}x^{n}=1+2x+10x^{2}+44x^{3}+250x^{4}+1412x^{5}+... 
\]%
satisfying%
\[
\theta ^{3}-4x(2\theta +1)(3\theta ^{2}+3\theta +1)-16x^{2}(\theta
+1)(2\theta +1)(2\theta +3)+192x^{3}(2\theta +1)(2\theta +3)(2\theta +5) 
\]%
Then%
\[
J_{24}=\frac{1}{q}+8+52q+834q^{3}+4760q^{5}+24703q^{7}+... 
\]%
which is 4B (a translation by $8$ of \ $\dbinom{2n}{n}\ast $(d) ).

\textbf{25.1.}%
\[
J_{25}(\exp (-\pi \sqrt{\frac{3}{2}}))=56 
\]%
\[
\sum_{n=0}^{\infty }A_{n}(5+24n)\frac{1}{(56)^{n}}=\frac{7\sqrt{14}}{\pi } 
\]

\textbf{25.2.}%
\[
J_{25}(\exp (-\pi \sqrt{\frac{5}{2}}))=152 
\]%
\[
\sum_{n=0}^{\infty }A_{n}(11+60n)\frac{1}{(152)^{n}}=\frac{19\sqrt{19}}{2\pi 
} 
\]

\textbf{25.3.}%
\[
J_{25}(\exp (-\pi \sqrt{\frac{9}{2}}))=792 
\]%
\[
\sum_{n=0}^{\infty }A_{n}(41+280n)\frac{1}{(792)^{n}}=\frac{33\sqrt{33}}{\pi 
\sqrt{2}} 
\]

\textbf{25.4.}%
\[
J_{25}(\exp (-\pi \sqrt{\frac{11}{2}}))=1592 
\]%
\[
\sum_{n=0}^{\infty }A_{n}(1237+9240n)\frac{1}{(1592)^{n}}=\frac{199\sqrt{398}%
}{\pi } 
\]

\textbf{25.5.}%
\[
J_{25}(\exp (-\pi \sqrt{\frac{29}{2}}))=156824 
\]%
\[
\sum_{n=0}^{\infty }A_{n}(873354+91045040n)\frac{1}{(156824)^{n}}=\frac{%
19603^{3/2}}{\pi } 
\]

\textbf{25.6.}%
\[
J_{25}(-\exp (-\pi \sqrt{\frac{3}{2}}))=-40 
\]%
\[
\sum_{n=0}^{\infty }A_{n}(7+24n)\frac{1}{(-40)^{n}}=\frac{5\sqrt{10}}{\pi } 
\]

\textbf{25.7.}%
\[
J_{25}(-\exp (-\pi \sqrt{\frac{5}{2}}))=-136 
\]%
\[
\sum_{n=0}^{\infty }A_{n}(13+60n)\frac{1}{(-136)^{n}}=\frac{17\sqrt{17}}{%
2\pi } 
\]

\textbf{25.8.}%
\[
J_{25}(-\exp (-\pi \sqrt{\frac{9}{2}}))=-776 
\]%
\[
\sum_{n=0}^{\infty }A_{n}(43+280n)\frac{1}{(-776)^{n}}=\frac{97\sqrt{97}}{%
3\pi \sqrt{6}} 
\]

\textbf{25.9.}%
\[
J_{25}(-\exp (-\pi \sqrt{\frac{11}{2}}))=-1576 
\]%
\[
\sum_{n=0}^{\infty }A_{n}(1271+9240n)\frac{1}{(-1576)^{n}}=\frac{197\sqrt{394%
}}{\pi } 
\]

2\textbf{5.10.}%
\[
J_{25}(-\exp (-\pi \sqrt{\frac{29}{2}}))=-156808 
\]%
\[
\sum_{n=0}^{\infty }A_{n}(873798+10450440n)\frac{1}{(-156808)^{n}}=\frac{%
(17\cdot 1153)^{3/2}}{\pi } 
\]

\textbf{25.11}%
\[
J_{25}(\exp (-\pi \sqrt{\frac{51}{2}}))=3879224+940800\sqrt{17} 
\]%
\[
\sum_{n=0}^{\infty }A_{n}\frac{10932560482+2724143108\sqrt{17}%
+(174961501680+42847402080\sqrt{17})n}{(3879224+940800\sqrt{17})^{n}} 
\]%
\[
=\frac{24999409}{\pi }\sqrt{3879224+940800\sqrt{17}} 
\]

\textbf{25.12.}%
\[
J_{25}(-\exp (-\pi \sqrt{\frac{51}{2}}))=-3879208-940800\sqrt{17} 
\]%
\[
\sum_{n=0}^{\infty }A_{n}\frac{10051483618+2520838532\sqrt{17}%
+(161321188560+39539139360\sqrt{17})n}{(-3879208+940800\sqrt{17})^{n}} 
\]%
\[
=\frac{23059801}{\pi }\sqrt{3879208+940800\sqrt{17}} 
\]

\textbf{25.13}%
\[
J_{25}(i\exp (-\pi \sqrt{\frac{5}{4}}))=8+32i 
\]%
\[
\sum_{n=0}^{\infty }A_{n}(6+i+20n)\frac{1}{(8+32i)^{n}}=\frac{34\sqrt{17}}{%
\pi \sqrt{52+47i}} 
\]

\textbf{25.14.}%
\[
J_{25}(i\exp (-\pi \sqrt{\frac{7}{4}}))=8+63i 
\]%
\[
\sum_{n=0}^{\infty }A_{n}(336+44i+1365n)\frac{1}{(8+63i)^{n}}=\frac{8066%
\sqrt{37\cdot 109}}{\pi \sqrt{237951+94744i}} 
\]

\textbf{25.15}%
\[
J_{25}(i\exp (-\pi \sqrt{\frac{13}{4}}))=8+288i 
\]%
\[
\sum_{n=0}^{\infty }A_{n}(138+7i+780n)\frac{1}{(8+288i)^{n}}=\frac{2594\sqrt{%
1297}}{\pi \sqrt{46548+3887i}} 
\]

\textbf{25.16.}%
\[
J_{25}(i\exp (-\pi \sqrt{\frac{37}{4}}))=8+14112i 
\]%
\[
\sum_{n=0}^{\infty }A_{n}(47166+101i+450660n)\frac{1}{(8+14112i)^{n}}=\frac{%
6223394\sqrt{17\cdot 183041}}{\pi \sqrt{5489026452+9335087i}} 
\]

\textbf{25.17.}%
\[
J_{25}(i\exp (-\pi \sqrt{\frac{9}{4}}))=8+64i\sqrt{3} 
\]%
\[
\sum_{n=0}^{\infty }A_{n}(-18-i\sqrt{3}-84n)\frac{1}{(8+64i\sqrt{3})^{n}}=%
\frac{(1+i)(1+8i\sqrt{3})^{3/2}}{3^{1/4}\pi } 
\]

\textbf{25.18.}%
\[
J_{25}(i\exp (-\pi \sqrt{\frac{25}{4}}))=8+1152i\sqrt{5} 
\]%
\[
\sum_{n=0}^{\infty }A_{n}(246+i\sqrt{5}+1932n)\frac{1}{(8+1152i\sqrt{5})^{n}}%
=\frac{\sqrt{2}(1+i)(1+144i\sqrt{5})^{3/2}}{10\cdot 5^{1/4}\pi } 
\]

\textbf{25.19}%
\[
J_{25}(\exp (\frac{\pi i}{4}-\pi \sqrt{\frac{15}{16}}))=\frac{1}{2}(49+15i%
\sqrt{3}) 
\]%
\[
\sum_{n=0}^{\infty }A_{n}(186i+398\sqrt{3}+1365n\sqrt{3})\frac{1}{(\dfrac{1}{%
2}(49+15i\sqrt{3}))^{n}}=\frac{39988\sqrt{3\cdot 769}}{\pi \sqrt{%
121547+143852i\sqrt{3}}} 
\]%
\textbf{.}

\textbf{Exempel 26.}

Consider%
\[
A_{n}=\sum_{k=0}^{n}\dbinom{2k}{k}\dbinom{4k}{2k}\dbinom{2n-2k}{n-k}\dbinom{%
4n-4k}{2n-2k}
\]%
with%
\[
y_{0}=\sum_{n=0}^{\infty
}A_{n}x^{n}=1+24x+984x^{2}+47040x^{3}+2421720x^{3}+...
\]%
satisfying%
\[
\theta ^{3}-8x(2\theta +1)(8\theta ^{2}+8\theta +3)+1024x^{2}(\theta
+1)(2\theta +1)(2\theta +3)
\]%
Then%
\[
J_{26}=\frac{1}{q}+40+276q-2048q^{2}+11202q^{3}-49152q^{4}+...\text{ \ \ \ \
2B}
\]

\textbf{26.1.}%
\[
J_{26}(-\exp (-\pi \sqrt{27}))=-4096192-3251200\cdot 2^{1/3}-2580480\cdot
2^{2/3} 
\]%
\[
\sum_{n=0}^{\infty }A_{n}\frac{\{-45-18\cdot 2^{1/3}+45\cdot
2^{2/3}+12(1+2\cdot 2^{1/3}+2^{2/3})n\}}{(-4096192-3251200\cdot
2^{1/3}-2580480\cdot 2^{2/3})^{n}}=\frac{\sqrt{3}}{\pi }(41+22\cdot
2^{1/3}-39\cdot 2^{2/3}) 
\]

\textbf{26.2.} \ $J_{26}(\exp (-\pi \sqrt{42}))$ \ satisfies the equation%
\[
x^{4}-695296512x^{3}+736297058304x^{2}-88550154436608x+2833604941971456=0 
\]%
with the solution%
\[
384(452667+184800\sqrt{6}+120980\sqrt{14}+98780\sqrt{21}) 
\]%
\[
w_{0}=\frac{1}{J_{26}(-\exp (-\pi \sqrt{42}))}=\frac{1}{128}-\frac{5(\sqrt{14%
}+11\sqrt{21})}{2^{5}\cdot 3\cdot 19^{2}}=\frac{1}{2^{7}3^{2}19^{2}}(6+14%
\sqrt{6}-12\sqrt{14}+\sqrt{21})^{2} 
\]%
\[
\sum_{n=0}^{\infty }A_{n}\left\{ -278+276\sqrt{6}+(280+1540\sqrt{6}+1083%
\sqrt{14})n\right\} w_{0}^{n}=\frac{722\sqrt{3}}{\pi } 
\]

This Example is a transformation of Example 6, which will be treated in [4].

\textbf{26.3.}%
\[
J_{26}(i\exp (-\pi \sqrt{\frac{7}{4}}))=\frac{1}{2}(81+45i\sqrt{7}) 
\]%
\[
\sum_{n=0}^{\infty }A_{n}\frac{16i+(35i+9\sqrt{7})n}{(\dfrac{1}{2}(81+45i%
\sqrt{7}))^{n}}=\frac{36}{\pi } 
\]

\textbf{Example 27.}

Consider%
\[
A_{n}=8^{n}\sum_{k=0}^{n}\dbinom{2k}{k}\dbinom{2n-2k}{n-k}\dbinom{-1/4}{k}%
\dbinom{-3/4}{n-k} 
\]%
with%
\[
y_{0}=\sum_{n=0}^{\infty
}A_{n}x^{n}=1-16x+360x^{2}-9088x^{3}+243160x^{4}-... 
\]%
satisfying%
\[
\theta ^{3}+16x(2\theta +1)(2\theta ^{2}+2\theta +1)+64x^{2}(\theta
+1)(4\theta +3)(4\theta +5) 
\]%
Then%
\[
J_{26}=\frac{1}{q}-16+76q-702q^{2}+5224q^{3}-23425q^{4}+... 
\]%
This example occurs in [20] as part of Conjecture 4 (i). Observe the power \ 
$8^{n}$, the more natural \ $16^{n}$ does not work.

\textbf{27.1.}%
\[
J_{27}(\exp (-\pi \sqrt{\frac{177}{4}}))=144(4135531+310845\sqrt{177}) 
\]%
\[
\sum_{n=0}^{\infty }A_{n}\frac{6654067259-271083501\sqrt{177}%
+(64663632548-73359420\sqrt{177})n}{(144(4135531+310845\sqrt{177}))^{n}}=%
\frac{2^{14}\cdot 3\cdot 11^{4}}{\pi }\sqrt{177} 
\]

\textbf{27.2.}%
\[
J_{27}(\exp (-\pi \sqrt{\frac{59}{12}}))=144(4135531-310845\sqrt{177}) 
\]%
\[
\sum_{n=0}^{\infty }A_{n}\frac{6654067259+271083501\sqrt{177}%
+(64663632548+73359420\sqrt{177})n}{(144(4135531-310845\sqrt{177}))^{n}}=%
\frac{2^{14}\cdot 3^{2}\cdot 11^{4}}{\pi }\sqrt{177} 
\]

\textbf{27.3.}%
\[
J_{27}(\exp (-\pi \sqrt{\frac{27}{4}}))=4(9+10\cdot 2^{1/3}+7\cdot
2^{2/3})^{2} 
\]%
\[
\sum_{n=0}^{\infty }A_{n}\frac{18(5-7\cdot 2^{1/3}+4\cdot
2^{2/3})+3(-111+51\cdot 2^{1/3}-186\cdot 2^{2/3})n}{(2(9+10\cdot
2^{1/3}+7\cdot 2^{2/3}))^{2n}}=\frac{(-186+14\cdot 2^{1/3}+156\cdot 2^{2/3})%
\sqrt{3}}{\pi } 
\]

\textbf{27.4.}%
\[
J_{27}(-\exp (-\pi \sqrt{\frac{27}{4}}))=-4(11+6\cdot 2^{1/3}+7\cdot
2^{2/3})^{2} 
\]%
\[
\sum_{n=0}^{\infty }A_{n}\frac{6(-7+3\cdot 2^{1/3}+7\cdot
2^{2/3})+3(3+29\cdot 2^{1/3}+59\cdot 2^{2/3})n}{(-4(11+6\cdot 2^{1/3}+7\cdot
2^{2/3})^{2})^{n}}=\frac{(18+14\cdot 2^{1/3}+34\cdot 2^{2/3})\sqrt{3}}{\pi } 
\]

\textbf{27.5}%
\[
J_{27}(i\exp (-\pi \sqrt{\frac{3}{2}}))=-16+32i\sqrt{2} 
\]%
\[
\sum_{n=0}^{\infty }A_{n}\frac{-5i+\sqrt{2}+(-8i+7\sqrt{2})n}{(-16+32i\sqrt{2%
}))^{n}}=\frac{6\sqrt{3}}{\pi } 
\]

\textbf{27.6.}%
\[
J_{27}(i\exp (-\pi \sqrt{\frac{5}{2}}))=-16+64i\sqrt{5} 
\]%
\[
\sum_{n=0}^{\infty }A_{n}\frac{-24i+15\sqrt{5}+(-40i+79\sqrt{5})n}{(-16+64i%
\sqrt{5}))^{n}}=\frac{3^{4}\sqrt{2}}{\pi } 
\]

\textbf{27.7}%
\[
J_{27}(i\exp (-\pi \sqrt{\frac{9}{2}}))=-16+320i\sqrt{6} 
\]%
\[
\sum_{n=0}^{\infty }A_{n}\frac{1077-69i\sqrt{6}+(7197-120i\sqrt{6})n}{%
(-16+320i\sqrt{6}))^{n}}=\frac{7^{4}\sqrt{2}}{\pi } 
\]

\textbf{27.8.}%
\[
J_{27}(i\exp (-\pi \sqrt{\frac{11}{2}}))=-16+1120i\sqrt{2} 
\]%
\[
\sum_{n=0}^{\infty }A_{n}\frac{-159i+1329\sqrt{2}+(-280i+9799\sqrt{2})n}{%
(-16+1120i\sqrt{2}))^{n}}=\frac{22\cdot 3^{4}\sqrt{11}}{\pi } 
\]

\textbf{27.9}%
\[
J_{27}(i\exp (-\pi \sqrt{\frac{29}{2}}))=-16+29120i\sqrt{29} 
\]%
\[
\sum_{n=0}^{\infty }A_{n}\frac{-57192i+8029839\sqrt{29}+(-105560i+96059599%
\sqrt{29})n}{(-16+29120i\sqrt{29}))^{n}}=\frac{99^{4}\sqrt{2}}{\pi } 
\]%
\textbf{.}

\textbf{Example 28.}

Finally we take an example where the $J-$function does not seem to work. Let%
\[
A_{n}=9^{n}\sum_{k=0}^{n}\dbinom{2k}{k}\dbinom{2n-2k}{n-k}\dbinom{-1/3}{k}%
\dbinom{-2/3}{n-k} 
\]%
with%
\[
y_{0}=\sum_{n=0}^{\infty
}A_{n}x^{n}=1-18x+450x^{2}-12636x^{3}+376650x^{4}+... 
\]%
satisfying%
\[
\theta ^{3}+18x(2\theta +1)(2\theta ^{2}+2\theta +1)+36x^{2}(\theta
+1)(6\theta +5)(6\theta +7) 
\]%
Then%
\[
J_{28}=\frac{1}{q}-54+891q-111537q^{3}+14033979q^{5}-... 
\]

\textbf{28.1.}%
\[
\sum_{n=0}^{\infty }A_{n}\frac{-15805+87665\sqrt{89}+(-28302+750006\sqrt{89}%
)n}{(-18-954\sqrt{89}))^{n}}=\frac{1500000\sqrt{3}}{\pi } 
\]

\textbf{28.2.}%
\[
\sum_{n=0}^{\infty }A_{n}\frac{15805+87665\sqrt{89}+(28302+750006\sqrt{89})n%
}{(-18+954\sqrt{89}))^{n}}=\frac{1500000\sqrt{3}}{\pi } 
\]%
In [20] there are the formulas 4.2-4.6 with $\dfrac{1}{x_{0}}=9\cdot
16,-9\cdot 20,9\cdot 108.-9\cdot 112$ but I could not find any reasonable
rational $t$ such that $J(\pm \exp (-\pi \sqrt{t}))$ assumed these values,
nor $\dfrac{1}{x_{0}}=-18\pm 954\sqrt{89}.$ Why the case $s=1/3$ (and $s=1/6$
) should be so different from $s=1/4$ is mysterious.

\textbf{Example 29.}

Let%
\[
W_{1}(x)=\frac{1}{\pi ^{3}}\int_{0}^{\pi }\int_{0}^{\pi }\int_{0}^{\pi }%
\frac{d\theta _{1}d\theta _{2}d\theta _{3}}{1-\dfrac{x}{3}(\cos (\theta
_{1})\cos (\theta _{2})+\cos (\theta _{1})\cos (\theta _{3})+\cos (\theta
_{2})\cos (\theta _{3}))} 
\]%
be the Lattice Green Function of the face-centered cubic lattice (after
talks by David Broadhurst and Tony Guttman).Then 
\[
y_{0}=W_{1}(x)=\sum_{n=0}^{\infty
}A_{n}x^{n}=1+12x^{2}+48x^{3}+540x^{4}+4320x^{5}+... 
\]%
where%
\[
A_{n}=\sum_{j=0}^{n}\sum_{k=0}^{n}(-4)^{n-k}\dbinom{n}{j}\dbinom{j}{k}^{2}%
\dbinom{2k}{k}\dbinom{2j-2k}{j-k} 
\]%
satisfies%
\[
\theta ^{3}-2x\theta (\theta +1)(2\theta +1)-16x^{2}(\theta +1)(5\theta
^{2}+10\theta +6)-96x^{3}(\theta +1)(\theta +2)(2\theta +3) 
\]%
Then we have%
\[
J_{29}=\frac{1}{q}+2+15q+32q^{2}+87q^{3}+192q^{4}+343q^{5}+... 
\]%
which is 6C. We note that \ $J_{29}=J_{\alpha }-4.$

\textbf{29.1.}%
\[
J_{29}(\exp (-\pi \sqrt{\frac{5}{3}}))=60 
\]%
\[
\sum_{n=0}^{\infty }A_{n}(1+4n)\frac{1}{60^{n}}=\frac{15\sqrt{3}}{8\pi } 
\]

\textbf{29.2.}%
\[
J_{29}(-\exp (-\pi \sqrt{\frac{4}{3}}))=-36 
\]%
\[
\sum_{n=0}^{\infty }A_{n}(1+4n)\frac{1}{(-36)^{n}}=\frac{27}{8\pi } 
\]

\textbf{29.3.}%
\[
J_{29}(\exp (-\pi \sqrt{3}))=76+64\cdot 2^{1/3}+48\cdot 2^{2/3} 
\]%
\[
\sum_{n=0}^{\infty }A_{n}\frac{\left\{ -144+252\cdot 2^{1/3}-66\cdot
2^{2/3}+(-224+392\cdot 2^{1/3}+64\cdot 2^{2/3})n\right\} }{(76+64\cdot
2^{1/3}+48\cdot 2^{2/3})^{n}}=\frac{125\sqrt{3}}{\pi } 
\]

\textbf{29.4.}%
\[
J_{29}(\exp (-\pi \sqrt{\frac{25}{3}}))=2892+1344\cdot 10^{1/3}+624\cdot
10^{2/3} 
\]%
\[
\sum_{n=0}^{\infty }A_{n}\frac{\left\{ 80-34\cdot 10^{1/3}+5\cdot
10^{2/3}+(80+20\cdot 10^{1/3}+32\cdot 10^{2/3})n\right\} }{(2892+1344\cdot
10^{1/3}+624\cdot 10^{2/3})^{n}}=\frac{243\sqrt{15}}{10\pi } 
\]

\textbf{29.5}%
\[
J_{29}(\exp (-\pi \sqrt{\frac{49}{3}}))=108876+23616\cdot 98^{1/3}+\frac{%
35856}{7}\cdot 98^{2/3} 
\]%
\[
\sum_{n=0}^{\infty }A_{n}\frac{\left\{ 125598-58156\cdot 98^{1/3}-26698\cdot
98^{2/3}+(1773408+423864\cdot 98^{1/3}+75712\cdot 98^{2/3})n\right\} }{%
(108876+23616\cdot 98^{1/3}+\dfrac{35856}{7}\cdot 98^{2/3})^{n}} 
\]%
\[
=\frac{3\cdot 55^{3}\sqrt{7}}{\pi } 
\]

\textbf{29.6.}%
\[
J_{29}(\exp (-\pi \sqrt{\frac{59}{3}}))=561804+324360\sqrt{3} 
\]%
\[
\sum_{n=0}^{\infty }A_{n}\frac{\left\{ 3085074-593849\sqrt{3}%
+(143226002+8270650\sqrt{3})n\right\} }{(561804+324360\sqrt{3})^{n}}=\frac{%
3^{2}\cdot 23^{3}\cdot 59}{\pi } 
\]

\textbf{29.7.}%
\[
J_{29}(i\exp (-\pi \sqrt{\frac{25}{12}}))=24(v-\frac{5}{v})\text{ \ \ where
\ }v=10^{1/3}\exp (\frac{2\pi i}{3}) 
\]%
\[
\sum_{n=0}^{\infty }A_{n}\frac{\left\{ 61+28v+\dfrac{130}{v}+3n\right\} }{%
(24(v-\dfrac{5}{v}))^{n}}=\frac{108}{\pi \sqrt{1815-390v-\dfrac{2100}{v}}} 
\]

\textbf{29.8}%
\[
J_{29}(i\exp (-\pi \sqrt{\frac{41}{12}}))=48(-v+\frac{16}{v})\text{ \ \
where \ }v=(2+10\sqrt{41})^{1/3}\exp (\frac{2\pi i}{3}) 
\]%
\[
\sum_{n=0}^{\infty }A_{n}\frac{\left\{ 53\sqrt{41}-(164v+\dfrac{2624}{v}-205%
\sqrt{41})n\right\} }{(48(-v+\dfrac{16}{v}))^{n}}=\frac{96\sqrt{123}}{\pi } 
\]

\textbf{.29.9.}%
\[
J_{29}(i\exp (-\pi \sqrt{\frac{49}{12}}))=72(-v+\frac{21}{v})\text{ \ \
where \ }v=98^{1/3}\exp (\frac{2\pi i}{3}) 
\]%
\[
\sum_{n=0}^{\infty }A_{n}\frac{\left\{ 29497+6396v+\dfrac{135954}{v}%
+55n\right\} }{(72(-v+\dfrac{21}{v}))^{n}}=\frac{5940}{\pi \sqrt{31759-3486v-%
\dfrac{72324}{v}}} 
\]

\textbf{29.10.}%
\[
J_{29}(i\exp (-\pi \sqrt{\frac{89}{12}}))=300(-v+\frac{100}{v})\text{ \ \
where \ }v=(2+106\sqrt{89})^{1/3}\exp (\frac{2\pi i}{3}) 
\]%
\[
\sum_{n=0}^{\infty }A_{n}\frac{\left\{ 827\sqrt{89}-(2225v+\dfrac{2222%
{\acute{}}%
500}{v}-4717\sqrt{89})n\right\} }{(300(-v+\dfrac{100}{v}))^{n}}=\frac{1500%
\sqrt{267}}{\pi } 
\]

\textbf{Remark.}

Also the other three Lattice Green Functions, $W_{2},W_{3},W_{4}$ \ give
formulas for \ $\dfrac{1}{\pi }$ But these functions are already known as, $%
\dbinom{2n}{n}\ast (c),$ $s=\dfrac{1}{2},$ $(\alpha )$ \ respectively.

\textbf{Harmonic numbers.}

In [3] there are many coefficients that are only expressed in Harmonic
Numbers%
\[
H_{n}=\sum_{k=1}^{n}\frac{1}{k} 
\]%
A similar search for third order equations is not so successful.We show
three examples.

\textbf{Example 30.}

Let 
\[
A_{n}=\dbinom{2n}{n}\sum_{k=0}^{n}\dbinom{n}{k}\dbinom{2k}{k}\dbinom{2n-2k}{%
n-k}\dbinom{3k}{k}\dbinom{3n-3k}{n-k}\left\{
1+k(-4H_{k}+4H_{n-k}+3H_{3k}-3H_{3n-3k})\right\} 
\]%
Then%
\[
y_{0}=\sum_{n=0}^{\infty }A_{n}x^{n}=1+42x+2870x^{2}+243420x^{3}+... 
\]%
satisfying%
\[
\theta ^{3}-6x(2\theta +1)(18\theta ^{2}+18\theta +7)+2916x^{2}(\theta
+1)(2\theta +1)(2\theta +3) 
\]%
with%
\[
J_{30}=\frac{1}{q}+66+783q-8672q^{2}+65367q^{3}+352q^{4}+870q^{5}+... 
\]%
which is 3A. We show only one formula of each kind, the $J$ -value being an
integer, cubic root, complex number but no quadratic roots.

\textbf{30.1}%
\[
J_{30}(\exp (-\pi \sqrt{\frac{67}{3}}))=27000108 
\]

\[
\sum_{n=0}^{\infty }A_{n}(87659+1500000n)\frac{1}{27000108^{n}}=\frac{16584}{%
\pi }\sqrt{267} 
\]

\textbf{30.2.}%
\[
J_{30}(\exp (-\pi \sqrt{\frac{361}{3}}))=12(292858+1097504\cdot
19^{1/3}+411296)\cdot 19^{2/3})^{2} 
\]%
\[
\sum_{n=0}^{\infty }A_{n}(a+bn)\frac{1}{(12(292858+1097504\cdot
19^{1/3}+411296)\cdot 19^{2/3})^{2})^{n}}=\frac{6\cdot (5\cdot 11\cdot
17\cdot 29)^{2}\sqrt{3}}{\pi } 
\]%
\newline
where%
\[
a=37089976361+878939544\cdot 19^{1/3}-5196816312\cdot 19^{2/3} 
\]%
\[
b=136652316192+5491939968\cdot 19^{1/3}-9478664064\cdot 19^{2/3} 
\]

\textbf{30.3}%
\[
J_{30}(\exp (\frac{\pi }{3}-\pi \sqrt{\frac{20}{9}}))=124+88i 
\]%
\[
\sum_{n=0}^{\infty }A_{n}(6+33i+(96+28in)\frac{1}{(124+88i)^{n}}=\frac{%
(45-24i)\sqrt{5}}{\pi } 
\]

\textbf{30.4.}%
\[
J_{30}(\exp (\frac{2\pi }{3}-\pi \sqrt{\frac{32}{9}}))=254+322i 
\]%
\[
\sum_{n=0}^{\infty }A_{n}(96+120i+(936+352i)n)\frac{1}{(254+322i)^{n}}=\frac{%
(432+51i)}{\pi }\sqrt{2} 
\]

\textbf{Example 31.}

Let 
\[
A_{n}=\dbinom{2n}{n}\sum_{k=0}^{n}\dbinom{n}{k}\dbinom{2k}{k}\dbinom{2n-2k}{%
n-k}\dbinom{4k}{2k}\dbinom{4n-4k}{2n-2k} 
\]%
\[
\cdot \left\{
1+k(-3H_{k}+3H_{n-k}-2H_{2k}+2H_{2n-2k}+4H_{4k}-4H_{4n-4k})\right\} 
\]%
Then%
\[
y_{0}=\sum_{n=0}^{\infty }A_{n}x^{n}=1+104x+17880x^{2}+3536960x^{3}+... 
\]%
satisfying%
\[
\theta ^{3}-8x(2\theta +1)(32\theta ^{2}+32\theta +13)+16484x^{2}(\theta
+1)(2\theta +1)(2\theta +3) 
\]%
with%
\[
J_{30}=\frac{1}{q}+152+4372q-96256q^{2}+65367q^{3}+352q^{4}+870q^{5}+... 
\]%
which is 3A.

\textbf{31.1.}%
\[
J_{31}(-\exp (-\pi \sqrt{58}))=-24591257600 
\]%
\[
\sum_{n=0}^{\infty }A_{n}\left\{ 8029841+192119202n\right\} \frac{1}{%
(-24591257600)^{n}}=\frac{1820^{2}\sqrt{58}}{\pi } 
\]

\textbf{31.2.}%
\[
J_{31}(\exp (-\pi \sqrt{27}))=400\cdot (57+46\cdot 2^{1/3}+38\cdot
2^{2/3})^{2} 
\]%
\[
\sum_{n=0}^{\infty }A_{n}\frac{a+bn}{(400\cdot (57+46\cdot 2^{1/3}+38\cdot
2^{2/3})^{2})^{n}}=\frac{(5\cdot 29\cdot 53)^{2}\sqrt{3}}{\pi } 
\]%
where%
\[
a=53732277-305352\cdot 2^{1/3}-13097682\cdot 2^{2/3} 
\]%
\[
b=530408529+18202176\cdot 2^{1/3}-13745664\cdot 2^{2/3} 
\]

\textbf{31.3.}%
\[
J_{31}(i\exp (-\pi \sqrt{\frac{15}{4}}))=\frac{5}{2}(61+99i\sqrt{3}) 
\]%
\[
\sum_{n=0}^{\infty }A_{n}\left\{ -84+108i+(216+423i)n\right\} \frac{1}{(%
\dfrac{5}{2}(61+99i\sqrt{3}))^{n}}=\frac{80\sqrt{15}+118i\sqrt{5}}{\pi } 
\]

\textbf{Example 32.}

Let 
\[
A_{n}=\dbinom{2n}{n}\sum_{k=0}^{n}\dbinom{n}{k}\dbinom{3k}{k}\dbinom{3n-3k}{%
n-k}\dbinom{6k}{3k}\dbinom{6n-6k}{3n-3k} 
\]%
\[
\left\{
1+k(-2H_{k}+2H_{n-k}-2H_{2k}+2H_{2n-2k}-3H_{3k}+3H_{3n-3k}+6H_{6k}-6H_{6n-6k})\right\} 
\]%
Then%
\[
y_{0}=\sum_{n=0}^{\infty }A_{n}x^{n}=1+744x+891864x^{2}+1218154560x^{3}+... 
\]%
satisfying%
\[
\theta ^{3}-24x(2\theta +1)(72\theta ^{2}+72\theta +31)+746496x^{2}(\theta
+1)(2\theta +1)(2\theta +3) 
\]%
with%
\[
J_{30}=\frac{1}{q}+984+196884q-21493760q^{2}+... 
\]%
which is 1A.

\textbf{32.1.}%
\[
J_{32}(\exp (-\pi \sqrt{163}))=262537412640769728 
\]%
\[
\sum_{n=0}^{\infty }A_{n}\left\{ 3787946075413+151931373056000n\right\} 
\frac{1}{262537412640769728^{n}}=\frac{3^{3}\cdot (7\cdot 11\cdot 19\cdot
127)^{2}\sqrt{163}}{\pi } 
\]

\textbf{32.2.}%
\[
J_{32}(-\exp (-\pi \sqrt{108}))=48(591239+469274\cdot 2^{1/3}+372458\cdot
2^{2/3})^{2} 
\]%
\[
\sum_{n=0}^{\infty }A_{n}(a+bn)\frac{1}{(48(591239+469274\cdot
2^{1/3}+372458\cdot 2^{2/3})^{2})^{n}}=\frac{(11\cdot 23\cdot 59\cdot
83\cdot 107)^{2}\sqrt{3}}{\pi } 
\]%
where%
\[
a=15352720352896986+30152156171347440\cdot 2^{1/3}-27499693962346020\cdot
2^{2/3} 
\]%
\[
b=288419215519028250+7821123473276000\cdot 2^{1/3}-44493090885216000\cdot
2^{2/3} 
\]

\textbf{Remark.}

The similar construction%
\[
A_{n}=\dbinom{2n}{n}\sum_{k=0}^{n}\dbinom{n}{k}\dbinom{2k}{k}^{2}\dbinom{%
2n-2k}{n-k}^{2}\left\{ 1+k(-5H_{k}+5H_{n-k}+4H_{2k}-4H_{2n-2k})\right\} 
\]%
is equal to the $A_{n}$ of Example 23, one of the numerous identities
involving harmonic numbers and binomial coefficients.

\textbf{Wronskians.}

If 
\[
y_{0}=\sum_{n=0}^{\infty }a_{n}x^{n} 
\]%
then%
\[
y_{1}=y_{0}\log (x)+\sum_{n=1}^{\infty }\frac{da_{n}}{dn}x^{n} 
\]%
If $a_{n}$ is a sum of a product of binomial coefficients, then $\dfrac{%
da_{n}}{dn}$ can be computed using the formula%
\[
\frac{d}{dn}\dbinom{an}{bn}=\dbinom{an}{bn}\left\{
aH_{an}-bH_{bn}-(a-b)H_{(a-b)n}\right\} 
\]

\textbf{Example 33.}

We consider the case $s=\dfrac{1}{2}$ where $y_{0}=\sum_{n=0}^{\infty }%
\dbinom{2n}{n}^{3}x^{n}.$Then the wronskian 
\[
w_{0}=x(y_{0}y_{1}^{\prime }-y_{0}^{\prime }y_{1})=\sum_{n=0}^{\infty
}A_{n}x^{n} 
\]%
with%
\[
A_{n}=\sum_{k=0}^{n}\dbinom{2k}{k}^{3}\dbinom{2n-2k}{n-k}^{3}\left\{
1+6k(-H_{k}+H_{n-k}+H_{2k}-H_{2n-2k})\right\} 
\]%
Then $w_{0}$ satisfies the differential equation%
\[
\theta ^{3}-8x(2\theta +1)(8\theta ^{2}+8\theta +5)+4096x^{2}(\theta +1)^{3} 
\]%
which happens to be the equation labelled $(\vartheta )$ in [1] with 
\[
A_{n}=\sum_{k=0}^{n}16^{n-k}\dbinom{2k}{k}^{3}\dbinom{2n-2k}{n-k} 
\]%
The $J-$function is invariant under forming wronskians, which follows from
that $y_{0}=u_{0}^{2},y_{1}=u_{0}u_{1}$ where $u_{0}$ and $u_{1}$ satisfy a
second order equation. So%
\[
J_{33}=\frac{1}{q}+24+276q+2048q^{2}+11202q^{3}+...=J_{1/2}=J_{6} 
\]

\textbf{33.1.}

\[
J_{33}(\exp (-\pi \sqrt{27}))=4^{4}(1+2^{1/3}+2^{2/3})^{8} 
\]%
\[
\sum_{n=0}^{\infty }A_{n}\frac{\{-792\cdot 2^{1/3}+630\cdot
2^{2/3}+(27-900\cdot 2^{1/3}+720\cdot 2^{2/3})n\}}{%
(2(1+2^{1/3}+2^{2/3}))^{8n}}=\frac{4\sqrt{3}}{\pi } 
\]

\textbf{Example 34.}

Consider the case $s=\dfrac{1}{3}.$Then the wronskian has 
\[
A_{n}=\sum_{k=0}^{n}\dbinom{2k}{k}^{2}\dbinom{2n-2k}{n-k}^{2}\dbinom{3k}{k}%
\dbinom{3n-3k}{n-k}\left\{
1+k(-5H_{k}+5H_{n-k}+2H_{2k}-2H_{2n-2k}+3H_{3k}-3H_{3n-3k})\right\} 
\]%
and%
\[
w_{0}=\sum_{n=0}^{\infty }A_{n}x^{n}=1+66x+5562x^{2}+508908x^{3}+... 
\]%
satisfies%
\[
\theta ^{3}-6x(2\theta +1)(18\theta ^{2}+18\theta +11)+324x^{2}(\theta
+1)(6\theta +5)(6\theta +7) 
\]%
with%
\[
J_{34}=\frac{1}{q}+42+783q^{2}+8672q^{2}+...=J_{1/3}=J_{7} 
\]

\textbf{34.1.}%
\[
J_{34}(\exp (-\pi \sqrt{\frac{89}{3}}))=-300^{3} 
\]%
\[
\sum_{n=0}^{\infty }A_{n}\frac{87665+1500006n}{(-300)^{3n}}=\frac{1500000%
\sqrt{3}}{\pi } 
\]

\textbf{34.2.}

\[
J_{34}(-\exp (-\pi \sqrt{\frac{361}{3}}))=-6^{6}(6617733763+2480036604\cdot
19^{1/3}+929409036\cdot 19^{2/3}) 
\]%
\[
\sum_{n=0}^{\infty }A_{n}\frac{a+bn}{(-6^{6}(6617733763+2480036604\cdot
19^{1/3}+9294090366\cdot 19^{2/3}))^{n}}=\frac{40500000\sqrt{3}}{\pi } 
\]%
where%
\[
a=-100930375+598500\cdot 19^{1/3}+17086500\cdot 19^{2/3} 
\]%
\[
b=478014198+23254632\cdot 19^{1/3}+32222088\cdot 19^{2/3} 
\]

\textbf{34.3.}

\[
J_{30}(\exp (\frac{2\pi }{3}-\pi \sqrt{\frac{32}{9}}))=-146+322i 
\]%
\[
\sum_{n=0}^{\infty }A_{n}(-12+39i+(17+144i)n)\frac{1}{(-146+322i)^{n}}=\frac{%
(132+351i)}{8\pi }\sqrt{2} 
\]

\textbf{Example 35.}

Consider the case $s=\dfrac{1}{4}.$Then the wronskian has 
\[
A_{n}=\sum_{k=0}^{n}\dbinom{2k}{k}^{2}\dbinom{2n-2k}{n-k}^{2}\dbinom{4k}{2k}%
\dbinom{4n-4k}{2n-2k}\left\{ 1+4k(-H_{k}+H_{n-k}+H_{4k}-H_{4n-4k})\right\} 
\]%
and%
\[
w_{0}=\sum_{n=0}^{\infty }A_{n}x^{n}=1+152x+30168x^{2}+6524864x^{3}+... 
\]%
satisfies%
\[
\theta ^{3}-8x(2\theta +1)(32\theta ^{2}+32\theta +19)+4096x^{2}(\theta
+1)(4\theta +3)(4\theta +5) 
\]%
with%
\[
J_{35}=\frac{1}{q}+104+4372q+96256q^{2}+...=J_{1/4}=J_{8} 
\]

\textbf{35.1.}%
\[
J_{35}(\exp (-\pi \sqrt{58}))=396^{4} 
\]%
\[
\sum_{n=0}^{\infty }A_{n}\frac{8029839+192119200n}{(396)^{4n}}=\frac{99^{4}}{%
\pi }\sqrt{\frac{2}{29}} 
\]

\textbf{35.2.}

\[
J_{35}(-\exp (-\pi \sqrt{27}))=-16(299+234\cdot 2^{1/3}+178\cdot
2^{2/3})^{2} 
\]%
\[
\sum_{n=0}^{\infty }A_{n}\frac{a+bn}{(-16(299+234\cdot 2^{1/3}+178\cdot
2^{2/3})^{2})^{n}}=\frac{(11\cdot 23)^{4}\sqrt{3}}{\pi } 
\]%
where%
\[
a=3727841133-1404744840\cdot 2^{1/3}+189863070\cdot 2^{2/3} 
\]%
\[
b=36821378025-2111995200\cdot 2^{1/3}+995827200\cdot 2^{2/3} 
\]

\textbf{35.3.}

\[
J_{35}(i\exp (-\pi \sqrt{\frac{15}{4}}))=\frac{1}{2}(207+495i\sqrt{3}) 
\]%
\[
\sum_{n=0}^{\infty }A_{n}(4i+20\sqrt{3}+(-59i+40\sqrt{3})n)\frac{1}{(\frac{1%
}{2}(207+495i\sqrt{3}))^{n}}=\frac{144\sqrt{15}+495i\sqrt{5}}{25\pi } 
\]

\textbf{.Example 36.}

Consider the case $s=\dfrac{1}{6}.$Then the wronskian has 
\[
A_{n}=\sum_{k=0}^{n}\dbinom{2k}{k}\dbinom{2n-2k}{n-k}\dbinom{3k}{k}\dbinom{%
3n-3k}{n-k}\dbinom{6k}{3k}\dbinom{6n-6k}{3n-3k} 
\]%
\[
\left\{ 1+k(-3H_{k}+3H_{n-k}-3H_{3k}+3H_{3n-3k}+6K_{6k}-6H_{6n-6k})\right\} 
\]%
and%
\[
w_{0}=\sum_{n=0}^{\infty }A_{n}x^{n}=1+984x+1306584x^{2}+1900332480x^{3}+... 
\]%
satisfies%
\[
\theta ^{3}-24x(2\theta +1)(72\theta ^{2}+72\theta +41)+331776x^{2}(\theta
+1)(3\theta +2)(3\theta +4) 
\]%
with%
\[
J_{36}=\frac{1}{q}+744+196884q+96256q^{2}+...=J_{1/4}=J_{8} 
\]

\textbf{36.1.}

\[
J_{36}(-\exp (-\pi \sqrt{163}))=-640320^{3} 
\]%
\[
\sum_{n=0}^{\infty }A_{n}\left\{ 3787946075414+15193137305600n\right\} \frac{%
1}{(-640320)^{3n}}=\frac{(16\cdot 5\cdot 23\cdot 29)^{2}}{\pi \sqrt{163}} 
\]

\textbf{36.2.}%
\[
J_{31}(\exp (-\pi \sqrt{108}))=6000\cdot (8389623817+6658848836\cdot
2^{1/3}+5285131824\cdot 2^{2/3}) 
\]%
\[
\sum_{n=0}^{\infty }A_{n}\frac{a+bn}{(6000\cdot (8389623817+6658848836\cdot
2^{1/3}+5285131824\cdot 2^{2/3}))^{n}}=\frac{(5\cdot 11\cdot 17)^{2}\sqrt{3}%
}{\pi } 
\]%
where%
\[
a=695116998-2158965360\cdot 2^{1/3}+1559571444\cdot 2^{2/3} 
\]%
\[
b=13702630662-2827055520\cdot 2^{1/3}+2880458496\cdot 2^{2/3} 
\]

We could not find any complex formula.

\textbf{Example 37.}

Taking the wronskian can be iterated. Let us start with the hypergeometric
case%
\[
a_{n}=C^{n}(1/2)_{n}(s)_{n}(1-s)_{n} 
\]%
where $s=\dfrac{1}{2},\dfrac{1}{3},\dfrac{1}{4},\dfrac{1}{6}$ and $%
C=64,108,256,1732$ respectively. Then if we iterate the wronskian $p$ times
we find%
\[
A_{n}=\sum_{k=0}^{n}(C/4)^{k}\frac{\dbinom{2k+2p}{2p}\dbinom{2k}{k}}{\dbinom{%
k+p}{p}}a_{n-k} 
\]%
Assume we have a formula%
\[
\sum_{n=0}^{\infty }a_{n}(a+bn)x_{0}^{n}=\frac{1}{\pi } 
\]%
Then after iterating the wronskian $p$ times we find the formula%
\[
\sum_{n=0}^{\infty }a_{n}(A+Bn)x_{0}^{n}=\frac{1}{\pi } 
\]%
with%
\[
A=\left\{ a-(p+\frac{1}{2})\frac{Cbx_{0}}{1-Cx_{0}}\right\}
(1-Cx_{0})^{p+1/2} 
\]%
\[
B=b(1-Cx_{0})^{p+1/2} 
\]

\textbf{37.1.}

Let us take $s=\dfrac{1}{3},C=108,p=11.$We find%
\[
A_{n}=\sum_{k=0}^{n}27^{k}\frac{\dbinom{2k+22}{22}\dbinom{2k}{k}}{\dbinom{%
k+11}{11}}\dbinom{2n-2k}{n-k}^{2}\dbinom{2n-3k}{n-k} 
\]%
and%
\[
Y_{0}=\sum_{n=0}^{\infty }A_{n}x^{n}=1+1254x+853794x^{2}+418202580x^{2}+... 
\]%
satisfying%
\[
\theta ^{3}-6x(54\theta ^{3}+648\theta ^{2}+634\theta
+209)+2^{2}3^{4}x^{2}(108\theta ^{3}+2592\theta ^{2}+16819\theta +14797) 
\]%
\[
-2^{4}3^{7}x^{3}(\theta +12)(6\theta +71)(6\theta +73) 
\]%
Then the formula 
\[
\sum_{n=0}^{\infty }\dbinom{2n}{n}^{2}\dbinom{3n}{n}\frac{827+14151n}{%
(-300)^{3n}}=\frac{1500\sqrt{3}}{\pi } 
\]%
is transformed to%
\[
\sum_{n=0}^{\infty }A_{n}\frac{7\cdot 83\cdot 151+6\cdot 53^{2}\cdot 89\cdot
n}{(-300)^{3n}}=\frac{3\cdot 2^{49}\cdot 5^{72}}{53^{22}\cdot 89^{11}\cdot
\pi }\sqrt{\frac{3}{89}} 
\]%
\ 

\textbf{Example 38.}

The formula%
\[
A_{n}=\sum_{k=0}^{[n/5]}(-1)^{k}\dbinom{n}{k}^{3}\left\{ \dbinom{4n-5k-1}{3n}%
+\dbinom{4n-5k}{3n}\right\} 
\]%
for the case $\ (\eta )$ \ was found during a search for formulas with
harmonic numbers. It is the only case where there are no rational \ $x_{0}.$
The closest to a rational number is \ $J_{31}(-\exp (-\pi \sqrt{\dfrac{4}{5}}%
)=-5\sqrt{5}$ which gives a convergent series but the convergence is so slow
that we could not find the formula even after summing 8000 terms.However
using the method explained in Example 61, we found formula \textbf{38.3 }%
below.

We have%
\[
y_{0}=\sum_{n=0}^{\infty }A_{n}x^{n}=1+5x+35x^{2}+275x^{3}+2275x^{4}+... 
\]%
satisfying%
\[
\theta ^{3}-x(2\theta +1)(11\theta ^{2}+11\theta +5)-125x^{2}(\theta +1)^{3} 
\]%
with%
\[
J_{38}=\frac{1}{q}+6+9q-10q^{2}-30q^{3}-6q^{4}-25q^{5}+... 
\]%
which is 5B.

\textbf{38.1.}%
\[
J_{38}(-\exp (-\pi \sqrt{\frac{8}{5}}))=-(\frac{5+\sqrt{5}}{2})^{3} 
\]%
\[
\sum_{n=0}^{\infty }A_{n}\frac{6-\sqrt{5}+12n}{(-\dfrac{5+\sqrt{5}}{2})^{3n}}%
=\frac{5\sqrt{5}}{2\pi \sqrt{5-2\sqrt{5}}} 
\]

\qquad\ 

\textbf{38.2}%
\[
J_{38}(\exp (-\pi \sqrt{\frac{47}{5}}))=(\frac{25+11\sqrt{5}}{2})^{3} 
\]

\[
\sum_{n=0}^{\infty }A_{n}\frac{141-50\sqrt{5}+282n}{(\dfrac{25+11\sqrt{5}}{2}%
)^{3n}}=\frac{5\sqrt{5}}{3\pi \sqrt{1525-682\sqrt{5}}} 
\]

\textbf{38.3.}%
\[
J_{38}(-\exp (-\pi \sqrt{\frac{4}{5}}))=-5\sqrt{5} 
\]%
\[
\sum_{n=0}^{\infty }A_{n}\frac{1+2n}{(-5\sqrt{5})^{n}}=\frac{25}{\pi \sqrt{%
250+110\sqrt{5}}} 
\]

\textbf{Example 39.}

Let%
\[
A_{n}=\dbinom{2n}{n}\sum_{k=0}^{n}4^{n-k}\dbinom{2k}{k}^{2}\dbinom{2n-2k}{n-k%
} 
\]%
Then%
\[
y_{0}=\sum_{n=0}^{\infty
}A_{n}x^{n}=1+24x+984x^{2}+47040x^{3}+2421720x^{4}+... 
\]%
satisfies%
\[
\theta ^{3}-8x(2\theta +1)(8\theta ^{2}+8\theta +3)+1024x^{2}(\theta
+1)(2\theta +1)(2\theta +3) 
\]%
with%
\[
J_{39}=\frac{1}{q}+40+276q-2048q^{2}+11202q^{3}-49152q^{4}+... 
\]%
which is 4A.

\textbf{39.1}%
\[
J_{39}(\exp (-\pi ))=72 
\]%
\[
\sum_{n=0}^{\infty }A_{n}\frac{-1+n}{72^{n}}=\frac{9}{\pi } 
\]

\textbf{39.2.}%
\[
J_{39}(\exp (-\pi \sqrt{2}))=128 
\]%
\[
\sum_{n=0}^{\infty }A_{n}\frac{n}{128^{n}}=\frac{\sqrt{2}}{\pi } 
\]

\textbf{39.3.}%
\[
J_{39}(\exp (-2\pi ))=576 
\]%
\[
\sum_{n=0}^{\infty }A_{n}\frac{2+16n}{576^{n}}=\frac{9}{\pi } 
\]

\textbf{39.4.}%
\[
J_{39}(-\exp (-\pi \sqrt{3}))=-192 
\]%
\[
\sum_{n=0}^{\infty }A_{n}\frac{1+4n}{(-192)^{n}}=\frac{\sqrt{3}}{\pi } 
\]

\textbf{39.5.}%
\[
J_{39}(-\exp (-\pi \sqrt{7}))=-4032 
\]%
\[
\sum_{n=0}^{\infty }A_{n}\frac{8+64n}{(-4032)^{n}}=\frac{9\sqrt{7}}{\pi } 
\]

\textbf{39.6.}%
\[
J_{39}(\exp (-\pi \sqrt{58}))=128(6930+1287\sqrt{29})^{2} 
\]%
\[
\sum_{n=0}^{\infty }A_{n}\frac{4412+(52780+9801\sqrt{29})n}{(128(6930+1287%
\sqrt{29})^{2})^{n}}=\frac{99^{2}}{\pi } 
\]

\textbf{39.7}%
\[
J_{39}(i\exp (-\pi ))=16i\sqrt{2} 
\]%
\textbf{.}

\[
\sum_{n=0}^{\infty }A_{n}\frac{64+3i\sqrt{2}+(176-22i\sqrt{2})n}{(16i\sqrt{2}%
)^{n}}=\frac{3872}{\pi \sqrt{464+190i\sqrt{2}}} 
\]

\textbf{32.9.}%
\[
J_{39}(-i\exp (-\pi \sqrt{\frac{7}{4}}))=\frac{1}{2}(81+45i\sqrt{7}) 
\]%
\[
\sum_{n=0}^{\infty }A_{n}\frac{16i+(35i+9\sqrt{7})n}{(\frac{1}{2}(81+45i%
\sqrt{7}))^{n}}=\frac{36}{\pi } 
\]

\textbf{Example 40.}

Let%
\[
A_{n}=\sum_{k=0}^{n}\dbinom{n-k}{k}\dbinom{2k}{k}^{2}\dbinom{2n-2k}{n-k} 
\]%
Then%
\[
y_{0}=\sum_{n=0}^{\infty }A_{n}x^{n}=1+2x+14x^{2}+68x^{3}+526x^{4}+... 
\]%
satisfies%
\[
\theta ^{3}-2x(8\theta ^{3}+6\theta ^{2}+4\theta +1)+2x^{2}(18\theta
^{3}-30\theta ^{2}-46\theta -18)+8x^{3}(58\theta ^{3}+273\theta
^{2}+353\theta +148) 
\]%
\[
-80x^{4}(36\theta ^{3}+144\theta ^{2}+183\theta +76)+32x^{5}(184\theta
^{3}+660\theta ^{2}+730\theta +255)-4096x^{6}(\theta +1)^{3} 
\]%
with%
\[
J_{40}=\frac{1}{q}+2+14q+40q^{3}+464q^{5}-1695q^{7}+... 
\]%
which is not in [13].

\textbf{40.1.}%
\[
J_{40}(\exp (-\pi \sqrt{\frac{7}{4}}))=2+10\sqrt{41} 
\]%
\[
\sum_{n=0}^{\infty }A_{n}\frac{-2665+2573\sqrt{41}+(-4305+10773\sqrt{41})n}{%
(2+10\sqrt{41})^{n}}=\frac{52480}{\pi } 
\]

\textbf{40.2}%
\[
J_{40}(-\exp (-\pi \sqrt{\frac{7}{4}}))=2+2\sqrt{65} 
\]%
\[
\sum_{n=0}^{\infty }A_{n}\frac{-130+34\sqrt{65}+(-195+99\sqrt{65})n}{(2+12%
\sqrt{65})^{n}}=\frac{1040}{\pi } 
\]

\textbf{Example 41.}

This example and the two following it, have essentially the same $J-$%
function, but the formulas for \ $\dfrac{1}{\pi }$ are rather different. Let%
\[
A_{n}=\sum_{k=0}^{n}\dbinom{2k}{k}^{3}\dbinom{2n-4k}{n-2k} 
\]%
which gives%
\[
y_{0}=\sum_{n=0}^{\infty
}A_{n}x^{n}=1+2x+14x^{2}+36x^{3}+334x^{4}+844x^{5}+... 
\]%
satisfying%
\[
\theta ^{3}-2x(6\theta ^{3}+3\theta ^{2}+3\theta +1)-4x^{2}(4\theta
^{3}+36\theta ^{2}+39\theta +15) 
\]%
\[
+8x^{3}(2\theta +3)(44\theta ^{2}+96\theta +69)-256x^{4}(12\theta
^{3}+48\theta ^{2}+69\theta +34)+512x^{5}(2\theta +3)^{3} 
\]%
We have 
\[
J_{41}=\frac{1}{q}+12q+66q^{3}+232q^{5}+639q^{7}+1596q^{9}+... 
\]%
which is 8B.

\textbf{41.1.}%
\[
J_{41}(\exp (-\pi \sqrt{\frac{3}{4}}))=16 
\]%
\[
\sum_{n=0}^{\infty }A_{n}\frac{1+6n}{16^{n}}=\frac{16}{\pi \sqrt{3}} 
\]

\textbf{41.2.}%
\[
J_{41}(\exp (-\pi \sqrt{\frac{7}{4}}))=64 
\]%
\[
\sum_{n=0}^{\infty }A_{n}\frac{43+210n}{64^{n}}=\frac{128}{\pi }\sqrt{\frac{5%
}{3}} 
\]

\textbf{41.3.}%
\[
J_{41}(-\exp (-\pi \sqrt{\frac{3}{4}}))=-16 
\]%
\[
\sum_{n=0}^{\infty }A_{n}\frac{13+30n}{(-16)^{n}}=\frac{16\sqrt{5}}{\pi } 
\]

\textbf{41.4.}%
\[
J_{41}(-\exp (-\pi \sqrt{\frac{7}{4}}))=-64 
\]%
\[
\sum_{n=0}^{\infty }A_{n}\frac{191+714n}{(-64)^{n}}=\frac{128\sqrt{17}}{\pi }
\]

\textbf{41.5.}%
\[
J_{41}(-\exp (-\pi \sqrt{\frac{27}{4}}))=-1168-928\cdot 2^{1/3}-736\cdot
2^{2/3} 
\]%
\[
\sum_{n=0}^{\infty }A_{n}\frac{-9077+27566\cdot 2^{1/3}-8078\cdot
2^{2/3}+(22770+45540\cdot 2^{1/3}+15180\cdot 2^{2/3})n}{(-1168-928\cdot
2^{1/3}-736\cdot 2^{2/3})^{n}} 
\]%
\[
=\frac{10120}{\pi \sqrt{-603-4362\cdot 2^{1/3}+3846\cdot 2^{1/3}}} 
\]

\textbf{41.6.}%
\[
J_{41}(i\exp (-\pi ))=16i\sqrt{2} 
\]%
\[
\sum_{n=0}^{\infty }A_{n}\frac{64+3i\sqrt{2}+(176-22i\sqrt{2})n}{(16i\sqrt{2}%
)^{n}}=\frac{3872}{\pi \sqrt{464+190i\sqrt{2}}} 
\]

\textbf{Example 42.}

Let%
\[
A_{n}=\sum_{k=0}^{n}2^{n-2k}\dbinom{n}{2k}\dbinom{2k}{k}^{3} 
\]%
Then%
\[
y_{0}=\sum_{n=0}^{\infty
}A_{n}x^{n}=1+2x+12x^{2}+56x^{3}+424x^{4}+2832x^{5}+... 
\]%
satisfies%
\[
\theta ^{3}-2x(2\theta +1)(2\theta ^{2}+2\theta +1)-4x^{2}(\theta
+1)(10\theta ^{2}+20\theta +9) 
\]%
\[
+112x^{3}(\theta +1)(\theta +2)(2\theta +3)-240x^{4}(\theta +1)(\theta
+2)(\theta +3) 
\]

We have 
\[
J_{42}=\frac{1}{q}+2+12q+66q^{3}+232q^{5}+639q^{7}+1596q^{9}+... 
\]

\textbf{42.1.}%
\[
J_{42}(\exp (-\pi \sqrt{\frac{3}{4}}))=18 
\]%
\[
\sum_{n=0}^{\infty }A_{n}\frac{1+4n}{18^{n}}=\frac{27}{4\pi } 
\]

\textbf{42.2.}%
\[
J_{42}(\exp (-\pi \sqrt{\frac{7}{4}}))=66 
\]%
\[
\sum_{n=0}^{\infty }A_{n}\frac{3+28n}{66^{n}}=\frac{363}{32\pi } 
\]

\textbf{42.3.}%
\[
J_{42}(-\exp (-\pi \sqrt{\frac{3}{4}}))=-14 
\]%
\[
\sum_{n=0}^{\infty }A_{n}\frac{5+12n}{(-14)^{n}}=\frac{7^{2}}{4\pi } 
\]

\textbf{42.4.}%
\[
J_{42}(-\exp (-\pi \sqrt{\frac{7}{4}}))=-62 
\]%
\[
\sum_{n=0}^{\infty }A_{n}\frac{352+1344n}{(-62)^{n}}=\frac{31^{2}}{\pi } 
\]

\textbf{42.5.}%
\[
J_{42}(\exp (-\pi \sqrt{\frac{27}{4}}))=1170+928\cdot 2^{1/3}+736\cdot
2^{2/3} 
\]%
\[
\sum_{k=0}^{n}\frac{762300-300288\cdot 2^{1/3}-223080\cdot
2^{2/3}+(721776-329760\cdot 2^{1/3}-38880\cdot 2^{2/3})n}{(1170+928\cdot
2^{1/3}+736\cdot 2^{2/3})^{n}}=\frac{307^{2}}{\pi } 
\]

\textbf{42.6.}%
\[
J_{42}(-\exp (-\pi \sqrt{\frac{27}{4}}))=-1170-928\cdot 2^{1/3}-736\cdot
2^{2/3} 
\]%
\[
\sum_{k=0}^{n}\frac{11305476+2838336\cdot 2^{1/3}-6752280\cdot
2^{2/3}+(16802544+22575840\cdot 2^{1/3}-7185120\cdot 2^{2/3})n}{%
(-1170-928\cdot 2^{1/3}-736\cdot 2^{2/3})^{n}}=\frac{3607^{2}}{\pi } 
\]

\textbf{Example 43.}

Let%
\[
A_{n}=\sum_{k=0}^{n}2^{n-2k}\dbinom{2k}{k}^{3}\dbinom{2n-4k}{n-2k} 
\]%
Then%
\[
y_{0}=\sum_{n=0}^{\infty
}A_{n}x^{n}=1+4x+32x^{2}+192x^{3}+1528x^{4}+10208x^{5}+... 
\]%
satisfies%
\[
\theta ^{3}-4x(2\theta +1)(\theta ^{2}+\theta +1)-16x^{2}(\theta +1)(4\theta
^{2}+8\theta +7)+64x^{3}(2\theta +3)^{3} 
\]%
We have%
\[
J_{43}=\frac{1}{q}+12q+66q^{3}+232q^{5}+639q^{7}+1596q^{9}+... 
\]

\textbf{43.1.}%
\[
J_{43}(\exp (-\pi \sqrt{\frac{3}{4}}))=16 
\]%
\[
\sum_{n=0}^{\infty }A_{n}\frac{1+6n}{16^{n}}=\frac{8\sqrt{2}}{\pi } 
\]

\textbf{43.2.}%
\[
J_{36}(\exp (-\pi \sqrt{\frac{7}{4}}))=64 
\]%
\[
\sum_{n=0}^{\infty }A_{n}\frac{1+6n}{64^{n}}=\frac{64}{7\pi }\sqrt{\frac{2}{7%
}} 
\]

\textbf{43.3.}%
\[
J_{43}(-\exp (-\pi \sqrt{\frac{3}{4}}))=-16 
\]%
\[
\sum_{n=0}^{\infty }A_{n}\frac{1+2n}{(-16)^{n}}=\frac{8}{3\pi }\sqrt{\frac{2%
}{3}} 
\]

\textbf{43.4.}%
\[
J_{43}(-\exp (-\pi \sqrt{\frac{7}{4}}))=-64 
\]%
\[
\sum_{n=0}^{\infty }A_{n}\frac{37+126n}{(-64)^{n}}=\frac{64\sqrt{2}}{\pi } 
\]

\textbf{43.5.}%
\[
J_{43}(\exp (-\pi \sqrt{\frac{27}{4}}))=1168+928\cdot 2^{1/3}+736\cdot
2^{2/3} 
\]%
\[
\sum_{n=0}^{\infty }A_{n}\frac{-459\sqrt{2}+78\cdot 2^{1/6}+318\cdot
2^{5/6}+(-366\sqrt{2}+60\cdot 2^{1/6}+276\cdot 2^{5/6})n}{(1168+928\cdot
2^{1/3}+736\cdot 2^{2/3})^{n}}=\frac{16}{\pi } 
\]

\textbf{43.6.}%
\[
J_{43}(-\exp (-\pi \sqrt{\frac{27}{4}}))=-1168-928\cdot 2^{1/3}-736\cdot
2^{2/3} 
\]%
\[
\sum_{n=0}^{\infty }A_{n}\frac{-315\sqrt{2}+1062\cdot 2^{1/6}-414\cdot
2^{5/6}+(-198\sqrt{2}+828\cdot 2^{1/6}-324\cdot 2^{5/6})n}{(-1168-928\cdot
2^{1/3}-736\cdot 2^{2/3})^{n}}=\frac{16\sqrt{3}}{\pi } 
\]

\textbf{The T-formulas by Z.W.Sun.}

In [19] Z.W.Sun considers%
\[
T_{n}(b,c)=\sum_{k=0}^{n}\dbinom{n}{k}\dbinom{n-k}{k}b^{n-2k}c^{k} 
\]%
We will only use \ $T_{n}(b)=T_{n}(b,1)$. The question is , how did Z.W.Sun
find his formulas? Let us take the most spectacular example, II.9%
\[
\sum_{n=0}^{\infty }\dbinom{2n}{n}\dbinom{3n}{n}T_{n}(287298)\frac{-7157+210n%
}{198^{3n}}=\frac{114345\sqrt{3}}{\pi } 
\]%
The sum converges like%
\[
(\frac{71825}{71874})^{n}=0.9993^{n} 
\]%
Already checking the formula numerically, seems impossible, since
computing,say 2000 terms does not give even the first digit correctly. Here
is Sun's answer (e-mail. March 2, 2011)

"I don't use any tricks for computation. I also have no access to fast
computer. My computation ability is preliminary. As I mentioned before, for
me, philosophy, intuition, inspiration, experience are more vital than
computation".

This, of course, does not give any information. Fortunately, H.H.Chan, J.Wan
and W.Zudilin found a way to prove II.9 by using Brafman's beautiful formula
to connect it to a Ramanujan-like formula ( $s=\dfrac{1}{3}$, $J_{7}(\exp
(-\pi \sqrt{\dfrac{136}{3}}))$

$=6^{3}(3555313+609700\sqrt{34})$ )%
\[
\sum_{n=0}^{\infty }\dbinom{2n}{n}^{2}\dbinom{3n}{n}\frac{-361+236\sqrt{34}%
+(1530+3420\sqrt{34})n}{6^{3n}(3555313+609700\sqrt{34})^{n}}=\frac{92\sqrt{6}%
+415\sqrt{51}}{\pi } 
\]%
(though this formula is not in [12]). Now we will use the machinery in [12]
to show how to find II.9, without going through the tedious proof involving
the modular equation for N=17. First we solve the equation%
\[
t(1-t)=\frac{27}{216(3555313+609700\sqrt{34})} 
\]%
We find%
\[
t=\frac{1}{2}(1-\rho _{0}-z_{0}) 
\]%
where%
\[
\rho _{0}=\frac{35\sqrt{2}}{71874}\text{ \ \ \ \ }z_{0}=\frac{8710\sqrt{17}}{%
35937} 
\]%
which gives%
\[
x_{0}=\frac{1-\rho _{0}^{2}-z_{0}^{2}}{2z_{0}}=\frac{143649\sqrt{17}}{592280}
\]%
and finally%
\[
b=\frac{2}{\sqrt{1-\dfrac{1}{x_{0}^{2}}}}=287298\text{ \ \ \ and \ }m=\frac{%
27b}{x_{0}z_{0}}=198^{3} 
\]

So now we know \ $b$ and $m$ \ and normally we should have used PSLQ to find
the coefficients $-7157$ \ and \ $210$, but that is hopeless with such a
slow convergence. Instead we sum to infinity by using integrals. Let \ $d=%
\sqrt{b^{2}-4}$. Then%
\[
T_{n}(b)=d^{n}P_{n}(\frac{b}{d}) 
\]%
where \ $P_{n}$ is the Legendre polynomial. Using the integral
representation of \ $P_{n}$ we get%
\[
T_{n}(b)=\frac{1}{\pi }\int_{0}^{\pi }(b-2\cos (\theta ))^{n}d\theta 
\]%
We have%
\[
\sum_{n=0}^{\infty }\dbinom{2n}{n}\dbinom{3n}{n}x^{n}=F(1/3,2/3;1;27x) 
\]%
where \ $F$ \ is the hypergeometric function. Furthermore%
\[
\sum_{n=0}^{\infty }n\dbinom{2n}{n}\dbinom{3n}{n}x^{n}=6xF(4/3,5/3;2;27x) 
\]%
so we have%
\[
a_{0}=\sum_{n=0}^{\infty }\dbinom{2n}{n}\dbinom{3n}{n}T_{n}(b)\frac{1}{m^{n}}%
=\frac{1}{\pi }\int_{0}^{\pi }F(1/3,2/3;1;\frac{27}{m}(b-2\cos (\theta
)))d\theta 
\]%
and%
\[
a_{1}=\sum_{n=0}^{\infty }n\dbinom{2n}{n}\dbinom{3n}{n}T_{n}(b)\frac{1}{m^{n}%
}=\frac{6}{\pi }\int_{0}^{\pi }(b.2\cos (\theta ))F(4/3,5/3;2;\frac{27}{m}%
(b-2\cos (\theta )))d\theta 
\]%
Using Maple, $a_{0}$ and $a_{1}$ can be computed with 60 digits in 8 sec on
a laptop. Now we can use PSLQ,preferable in the form%
\[
\text{PSLQ([a0\symbol{94}2,a0*a1,a1\symbol{94}2,1/Pi\symbol{94}2])} 
\]%
which also finds the square root on the right hand side.

Switching $\rho _{0}$ and $z_{0}$ we find the dual formula II.7%
\[
\sum_{n=0}^{\infty }\dbinom{2n}{n}\dbinom{3n}{n}T_{n}(198)\frac{15724+222105n%
}{198^{3n}}=\frac{114345\sqrt{3}}{4\pi } 
\]%
It is of the special type where \ $m=b^{3}.$ There are four other formulas
of this type, namely II.2, A1, A2, and II.6 with $b=10,18,30,102$
respectively. Allowing complex numbers, I found five more such cases with \ $%
b=6i\sqrt{3},6i\sqrt{11},30i\sqrt{3},90i\sqrt{3},138i\sqrt{59}.$They can all
be written in rational form by introducing the polynomial%
\[
U_{n}(x)=\sum_{k=0}^{n}\dbinom{n}{k}\dbinom{n-k}{k}x^{n+k} 
\]%
Then we have%
\[
\frac{T_{n}(b)}{b^{3n}}=\sum_{k=0}^{n}\dbinom{n}{k}\dbinom{n-k}{k}\frac{%
b^{n-2k}}{b^{3n}}=U_{n}(\frac{1}{b^{2}}) 
\]

\textbf{44.1.}%
\[
\sum_{n=0}^{\infty }\dbinom{2n}{n}\dbinom{3n}{n}U_{n}(-\frac{1}{108})(7+39n)=%
\frac{9\sqrt{3}}{\pi } 
\]

\textbf{44.2}%
\[
\sum_{n=0}^{\infty }\dbinom{2n}{n}\dbinom{3n}{n}U_{n}(-\frac{1}{396})(6+45n)=%
\frac{5\sqrt{11}}{\pi } 
\]

\textbf{44.3.}%
\[
\sum_{n=0}^{\infty }\dbinom{2n}{n}\dbinom{3n}{n}U_{n}(-\frac{1}{2700}%
)(1654+17157n)=\frac{2925\sqrt{3}}{\pi } 
\]

\textbf{44.4.}%
\[
\sum_{n=0}^{\infty }\dbinom{2n}{n}\dbinom{3n}{n}U_{n}(-\frac{1}{24300}%
)(7843+105339n)=\frac{14175\sqrt{3}}{\pi } 
\]

\textbf{44.5.}%
\[
\sum_{n=0}^{\infty }\dbinom{2n}{n}\dbinom{3n}{n}U_{n}(-\frac{1}{1123596}%
)(342786+6367095n)=\frac{140185\sqrt{59}}{\pi } 
\]

\textbf{Remark: }Observe that 
\[
\sum_{n=0}^{\infty }\dbinom{2n}{n}\dbinom{3n}{n}U_{n}(x)=\sum_{n=0}^{\infty
}A_{n}x^{n}=1+6x+90x^{2}+1860x^{3}+44730x^{4}+... 
\]%
where \ $A_{n}=\dbinom{2n}{n}\sum_{k=0}^{n}\dbinom{n}{k}^{2}\dbinom{2k}{k}$
which we recognize as case 1. But%
\[
\sum_{n=0}^{\infty }n\dbinom{2n}{n}\dbinom{3n}{n}%
U_{n}(x)=6x+180x^{2}+5400x^{3}+168840x^{4}+... 
\]%
is not equal to 
\[
\sum_{n=0}^{\infty }nA_{n}x^{n}=6x+180x^{2}+5580x^{3}+178920x^{4}+... 
\]

\textbf{Example 45.}

Doing the same trick to the AZ-numbers we define%
\[
y_{0}=\sum_{n=0}^{\infty }\sum_{k=0}^{n}(-1)^{k}3^{n-3k}\dbinom{n}{3k}%
\dbinom{n+k}{k}\frac{(3k)!}{k!^{3}}x^{n+k}=\sum_{n=0}^{\infty }A_{n}x^{n} 
\]%
\[
=1+3x+9x^{2}+27x^{3}+57x^{4}-117x^{5}-2511x^{6}-... 
\]%
Then we find%
\[
A_{n}=\sum_{k=0}^{n}(-1)^{k}3^{n-4k}\dbinom{n-k}{3k}\dbinom{n}{k}\frac{(3k)!%
}{k!^{3}} 
\]%
and the differential operator is%
\[
\theta ^{3}-3x(2\theta +1)(2\theta ^{2}+2\theta +1)+9x^{2}(\theta
+1)(6\theta ^{2}+12\theta +7) 
\]%
\[
-54x^{3}(\theta +1)(\theta +2)(2\theta +3)+337x^{4}(\theta +1)(\theta
+2)(\theta +3) 
\]%
with%
\[
J_{45}=\frac{1}{q}+3-26q^{3}+79q^{7}-326q^{11}+755q^{15}+... 
\]%
which is 8C. Most values of \ $J_{45}$ seem to satisfy equations of degree
four which all are solvable in simple radicals.We show a sample of formulas.

\textbf{45.1.}%
\[
J_{45}(\exp (-\pi \sqrt{\frac{5}{16}}))=3+4\sqrt{2} 
\]%
\[
\sum_{n=0}^{\infty }A_{n}\frac{636-411\sqrt{2}+(820-480\sqrt{2})n}{(3+4\sqrt{%
2})^{n}}=\frac{23^{2}}{\pi } 
\]

\textbf{45.2.}%
\[
J_{45}(-\exp (-\pi \sqrt{\frac{13}{16}}))=3+12\sqrt{2} 
\]%
\[
\sum_{n=0}^{\infty }A_{n}\frac{3372-1429\sqrt{2}+(8580-2080\sqrt{2})n}{(3+12%
\sqrt{2})^{n}}=\frac{93^{2}}{\pi } 
\]

\textbf{45.3.}%
\[
J_{45}(-\exp (-\pi \sqrt{\frac{25}{16}}))=3+24\cdot 20^{1/4} 
\]%
\[
\sum_{n=0}^{\infty }A_{n}\frac{a+bn}{(3+24\cdot 20^{1/4})^{n}}=\frac{9\cdot
81919^{2}}{\pi } 
\]%
where%
\[
a=594744820+8599866880\cdot \sqrt{5}-2647725120\cdot 20^{1/4}-5365775\cdot
20^{3/4} 
\]%
\[
b=791350420+33765383680\cdot \sqrt{5}-4220569920\cdot 20^{1/4}-6594560\cdot
20^{3/4} 
\]

\textbf{45.4.}%
\[
J_{45}(\exp (-\pi \sqrt{\frac{253}{16}}))=3+94164\sqrt{2}+28392\sqrt{22} 
\]%
\[
\sum_{n=0}^{\infty }A_{n}\frac{a+bn}{(3+94164\sqrt{2}+28392\sqrt{22})^{n}}=%
\frac{(3\cdot 3335791423)^{2}}{\pi } 
\]%
where%
\[
a=-3332721359191736134868-2299135623260036907019\cdot \sqrt{2} 
\]%
\[
+1014464862309261870848\cdot \sqrt{11}+693214440231976416094\cdot \sqrt{22} 
\]%
\[
b=-3076855543973703712220-2231460528289037325280\cdot \sqrt{2} 
\]%
\[
+10477795534221377888320\cdot \sqrt{11}+672808751556638831680\cdot \sqrt{22} 
\]

\textbf{45.5.}%
\[
J_{45}(-i\exp (-\pi \sqrt{\frac{5}{16}}))=3+4i\sqrt{2} 
\]%
\[
\sum_{n=0}^{\infty }A_{n}\frac{2516+1387i\sqrt{2}+(8060+2080i\sqrt{2})n}{%
(3+4i\sqrt{2})^{n}}=\frac{99^{2}}{\pi } 
\]

\textbf{45.6.}%
\[
J_{45}(-i\exp (-\pi \sqrt{\frac{25}{16}}))=3+24i\cdot 20^{1/4} 
\]%
\[
\sum_{n=0}^{\infty }A_{n}\frac{a+bn}{(3+24i\cdot 20^{1/4})^{n}}=\frac{%
(3\cdot 81919)^{2}}{\pi } 
\]%
where%
\[
a=-594744820+8599866880\cdot \sqrt{5}-2647725120i\cdot
20^{1/4}-5365775i\cdot 20^{3/4} 
\]%
\[
b=-791350420+33765383680\cdot \sqrt{5}+4220569920i\cdot
20^{1/4}-6594560i\cdot 20^{3/4} 
\]

\textbf{45.7.}%
\[
J_{45}(-i\exp (-\pi \sqrt{\frac{253}{16}}))=3+i(94164\sqrt{2}+28392\sqrt{22}%
) 
\]%
\[
\sum_{n=0}^{\infty }A_{n}\frac{a+bn}{(3+i(94164\sqrt{2}+28392\sqrt{22}))^{n}}%
=\frac{(3\cdot 2297\cdot 5410441)^{2}}{\pi } 
\]%
where%
\[
a=1371351262931220126164+1979346205791052668533i\sqrt{2} 
\]%
\[
-280069428474211619584\sqrt{11}-596780999185697244578i\sqrt{22} 
\]%
\[
b=3643075877893903390940+2231501382910522442720i\sqrt{2} 
\]%
\[
+568180426253290036160\sqrt{11}-672796434107683296320i\sqrt{22} 
\]

\textbf{45.8.}%
\[
J_{45}(\exp (-\frac{\pi i}{4}-\pi \sqrt{\frac{11}{8}}))=3+6(1+i)\sqrt{22} 
\]%
\[
\sum_{n=0}^{\infty }A_{n}\frac{102(-1+i)+150\sqrt{22}+560n\sqrt{22}}{%
(3+6(1+i)\sqrt{22})^{n}}=\frac{1}{\pi }\left\{ (1584-9i)\sqrt{2}+72(1-i)%
\sqrt{11}\right\} 
\]

\textbf{Example 46.}

Now we present all the new T-formulas we found.

$\mathbf{s=}\dfrac{\mathbf{1}}{\mathbf{3}}.$

\textbf{46.1.}%
\[
J_{7}(\exp (-\pi \sqrt{\frac{80}{3}}))=555378+3206385\sqrt{3} 
\]%
\[
\sum_{n=0}^{\infty }\dbinom{2n}{n}\dbinom{3n}{n}T_{n}(488)\frac{976+11310n}{%
33^{3n}}=\frac{4719\sqrt{3}}{\pi } 
\]%
\[
\sum_{n=0}^{\infty }\dbinom{2n}{n}\dbinom{3n}{n}T_{n}(843)\frac{48+2520n}{%
33^{3n}}=\frac{1573\sqrt{5}}{\pi } 
\]

\textbf{46.2.}%
\[
J_{7}(-\exp (-\pi \sqrt{\frac{65}{3}}))=-1121472-311040\sqrt{13} 
\]%
\[
\sum_{n=0}^{\infty }\dbinom{2n}{n}\dbinom{3n}{n}T_{n}(258)\frac{13+45n}{%
(-12)^{3n}}=\frac{2\sqrt{15}}{\pi } 
\]

\textbf{46.3.}%
\[
J_{7}(-\exp (-\pi \sqrt{\frac{145}{3}}))=-1529014752-126977760\sqrt{145} 
\]%
\[
\sum_{n=0}^{\infty }\dbinom{2n}{n}\dbinom{3n}{n}T_{n}(5778)\frac{7075+39585n%
}{(-72)^{3n}}=\frac{7344\sqrt{3}}{\pi } 
\]

\textbf{46.4.}%
\[
J_{7}(-\exp (-\pi \sqrt{\frac{185}{3}}))=-25891780032-4256582400\sqrt{37} 
\]%
\[
\sum_{n=0}^{\infty }\dbinom{2n}{n}\dbinom{3n}{n}T_{n}(18498)\frac{1043+7245n%
}{(-132)^{3n}}=\frac{605\sqrt{15}}{\pi } 
\]

\textbf{46.5.}%
\[
J_{7}(-\exp (-\pi \sqrt{\frac{209}{3}}))=-122165930592-16181272800\sqrt{57} 
\]%
\[
\sum_{n=0}^{\infty }\dbinom{2n}{n}\dbinom{3n}{n}T_{n}(93102)\frac{7467+23265n%
}{(-110592)^{n}}=\frac{448\sqrt{11}}{\pi } 
\]

\textbf{46.6.}%
\[
J_{7}(-\exp (-\pi \sqrt{\frac{265}{3}}))=-3327915564000-204432228000\sqrt{265%
} 
\]%
\[
\sum_{n=0}^{\infty }\dbinom{2n}{n}\dbinom{3n}{n}T_{n}(132498)\frac{%
7049+62403n}{(-360)^{3n}}=\frac{10800\sqrt{3}}{\pi } 
\]

\textbf{46.7.}%
\[
J_{7}(-\exp (-\pi \sqrt{\frac{11}{9}}))=-16736-2912\sqrt{33} 
\]%
\[
\sum_{n=0}^{\infty }\dbinom{2n}{n}\dbinom{3n}{n}T_{n}(\frac{730}{27})\frac{%
25+99n}{(-8)^{3n}}=\frac{16\sqrt{3}}{\pi } 
\]

$\mathbf{s=}\dfrac{\mathbf{1}}{\mathbf{4}}$

This case is treated in the same way as for $s=\dfrac{1}{3}$ , only $27$ is
replaced by $64.$ The cases with large \ $b$ \ have been computed with the
integrals%
\[
a_{0}=\sum_{n=0}^{\infty }\dbinom{2n}{n}\dbinom{4n}{2n}T_{n}(b)\frac{1}{m^{n}%
}=\frac{1}{\pi }\int_{0}^{\pi }F(1/4,3/4;1;\frac{64}{m}(b-2\cos (\theta
)))d\theta 
\]%
\[
a_{1}=\sum_{n=0}^{\infty }n\dbinom{2n}{n}\dbinom{4n}{2n}T_{n}(b)\frac{1}{%
m^{n}}=\frac{12}{\pi }\int_{0}^{\pi }(b.2\cos (\theta ))F(5/4,7/4;2;\frac{64%
}{m}(b-2\cos (\theta )))d\theta 
\]

\textbf{46.8.}%
\[
J_{7}(\exp (-\pi \sqrt{30}))=14858496+10506240\sqrt{2} 
\]%
\[
\sum_{n=0}^{\infty }\dbinom{2n}{n}\dbinom{4n}{2n}T_{n}(322)\frac{71+760n}{%
336^{2n}}=\frac{126\sqrt{7}}{\pi } 
\]%
\[
\sum_{n=0}^{\infty }\dbinom{2n}{n}\dbinom{4n}{2n}T_{n}(1442)\frac{-1+10n}{%
336^{2n}}=\frac{7}{2\pi }\sqrt{\frac{105}{2}} 
\]

\textbf{46.9.}%
\[
J_{7}(\exp (-\pi \sqrt{42}))=896645376+634023936\sqrt{2} 
\]%
\[
\sum_{n=0}^{\infty }\dbinom{2n}{n}\dbinom{4n}{2n}T_{n}(898)\frac{138+1540n}{%
912^{2n}}=\frac{95}{\pi }\sqrt{\frac{57}{2}} 
\]%
\[
\sum_{n=0}^{\infty }\dbinom{2n}{n}\dbinom{4n}{2n}T_{n}(12098)\frac{-139+280n%
}{912^{2n}}=\frac{95}{\pi }\sqrt{399} 
\]

\textbf{46.10.}%
\[
J_{7}(\exp (-\pi \sqrt{70}))=130065336576+41130270720\sqrt{10} 
\]%
\[
\sum_{n=0}^{\infty }\dbinom{2n}{n}\dbinom{4n}{2n}T_{n}(39202)\frac{68+920n}{%
3024^{2n}}=\frac{135}{\pi }\sqrt{\frac{21}{2}} 
\]%
\[
\sum_{n=0}^{\infty }\dbinom{2n}{n}\dbinom{4n}{2n}T_{n}(103682)\frac{-25+440n%
}{3024^{2n}}=\frac{378}{\pi }\sqrt{3} 
\]

\textbf{46.11.}%
\[
J_{7}(\exp (-\pi \sqrt{78}))=560809933056+396552499200\sqrt{2} 
\]%
\[
\sum_{n=0}^{\infty }\dbinom{2n}{n}\dbinom{4n}{2n}T_{n}(10402)\frac{%
12022+168740n}{10416^{2n}}=\frac{3689}{\pi }\sqrt{\frac{217}{2}} 
\]%
\[
\sum_{n=0}^{\infty }\dbinom{2n}{n}\dbinom{4n}{2n}T_{n}(1684802)\frac{%
-50087+8840n}{10416^{2n}}=\frac{7378}{\pi }\sqrt{8463} 
\]

\textbf{46.12.}%
\[
J_{7}(\exp (-\pi \sqrt{102}))=30095095654656+7299132825600\sqrt{17} 
\]%
\[
\sum_{n=0}^{\infty }\dbinom{2n}{n}\dbinom{4n}{2n}T_{n}(39202)\frac{%
732103+11657240n}{39216^{2n}}=\frac{80883}{\pi }\sqrt{817} 
\]%
\[
\sum_{n=0}^{\infty }\dbinom{2n}{n}\dbinom{4n}{2n}T_{n}(23990402)\frac{%
-58871+3080n}{39216^{2n}}=\frac{17974\sqrt{2451}}{\pi } 
\]

\textbf{46.13.}%
\[
J_{7}(\exp (-\pi \sqrt{130}))=1799922876293376+223252956241920\sqrt{65} 
\]%
\[
\sum_{n=0}^{\infty }\dbinom{2n}{n}\dbinom{4n}{2n}T_{n}(1684802)\frac{%
53378+691220n}{47376^{2n}}=\frac{14805}{\pi }\sqrt{\frac{329}{2}} 
\]%
\[
\sum_{n=0}^{\infty }\dbinom{2n}{n}\dbinom{4n}{2n}T_{n}(33385282)\frac{%
-160150+168740n}{47376^{2n}}=\frac{115479}{\pi }\sqrt{\frac{329}{2}} 
\]

\textbf{46.14.}%
\[
J_{7}(\exp (-\pi \sqrt{190}))=3203318942465423616+2265088546520616960\sqrt{2}
\]%
\[
\sum_{n=0}^{\infty }\dbinom{2n}{n}\dbinom{4n}{2n}T_{n}(33385282)\frac{%
1877581+27724840n}{440496^{2n}}=\frac{49266\sqrt{15295}}{\pi } 
\]%
\[
\sum_{n=0}^{\infty }\dbinom{2n}{n}\dbinom{4n}{2n}T_{n}(2998438562)\frac{%
-144982+36380n}{440496^{2n}}=\frac{156009}{\pi }\sqrt{\frac{161}{2}} 
\]

\textbf{46.15.}%
\[
J_{7}(\exp (-\pi \sqrt{7}))=8292456+31326750\sqrt{7} 
\]%
\[
\sum_{n=0}^{\infty }\dbinom{2n}{n}\dbinom{4n}{2n}T_{n}(4048)\frac{-784+455n}{%
513^{2n}}=\frac{2052\sqrt{57}}{\pi } 
\]

\textbf{46.16.}%
\[
J_{7}(-\exp (-\pi \sqrt{15}))=-\frac{1}{2}(192303+85995\sqrt{5}) 
\]%
\[
\sum_{n=0}^{\infty }\dbinom{2n}{n}\dbinom{4n}{2n}T_{n}(47)\frac{100+455n}{%
(-33^{2})^{n}}=\frac{22\sqrt{11}}{\pi } 
\]

\textbf{46.17.}%
\[
J_{7}(-\exp (-\pi \sqrt{21}))=-893952-516096\sqrt{3} 
\]%
\[
\sum_{n=0}^{\infty }\dbinom{2n}{n}\dbinom{4n}{2n}T_{n}(110)\frac{5+28n}{%
(-96^{2})^{n}}=\frac{3\sqrt{6}}{\pi } 
\]

\textbf{46.18.}%
\[
J_{7}(-\exp (-\pi \sqrt{85}))=-1896262493184-459911208960\sqrt{17} 
\]%
\[
\sum_{n=0}^{\infty }\dbinom{2n}{n}\dbinom{4n}{2n}T_{n}(103682)\frac{%
12197+97580n}{(-5472^{2})^{n}}=\frac{2736\sqrt{95}}{\pi } 
\]

\textbf{46.19.}%
\[
J_{7}(-\exp (-\pi \sqrt{93}))=-7186858066944-4149334425600\sqrt{3} 
\]%
\[
\sum_{n=0}^{\infty }\dbinom{2n}{n}\dbinom{4n}{2n}T_{n}(24302)\frac{%
71161+1071980n}{(-24288^{2})^{n}}=\frac{11385}{\pi }\sqrt{\frac{759}{2}} 
\]

\textbf{46.20.}%
\[
J_{7}(-\exp (-\pi \sqrt{133}))=-2714724158137344-622800434995200\sqrt{19} 
\]%
\[
\sum_{n=0}^{\infty }\dbinom{2n}{n}\dbinom{4n}{2n}T_{n}(5177198)\frac{%
88597+602140n}{(-26784^{2})^{n}}=\frac{4185\sqrt{1302}}{\pi } 
\]

\textbf{46.21.}%
\[
J_{7}(-\exp (-\pi \sqrt{253}))=-2515964063008169862144-758591707546926489600%
\sqrt{11} 
\]%
\[
\sum_{n=0}^{\infty }\dbinom{2n}{n}\dbinom{4n}{2n}T_{n}(6445234798)\frac{%
88597+602140n}{(-606816^{2})^{n}}=\frac{360297\sqrt{86}}{\pi } 
\]

$\mathbf{s}$\textbf{=}$\dfrac{\mathbf{1}}{\mathbf{2}}.$

Here we replace \ $27$ \ by \ $16.$ The following series for \ $\dfrac{1}{%
\pi }$ are all divergent but the corresponding integrals are convergent%
\[
a_{0}=\sum_{n=0}^{\infty }\dbinom{2n}{n}^{2}T_{n}(b)\frac{1}{m^{n}}=\frac{1}{%
\pi }\int_{0}^{\pi }F(1/2,1/2;1;\frac{16}{m}(b-2\cos (\theta )))d\theta 
\]%
\[
a_{1}=\sum_{n=0}^{\infty }n\dbinom{2n}{n}^{2}T_{n}(b)\frac{1}{m^{n}}=\frac{4%
}{\pi }\int_{0}^{\pi }(b-2\cos (\theta ))F(3/2,3/2;2;\frac{16}{m}(b-2\cos
(\theta )))d\theta 
\]

\textbf{46.22.}%
\[
J_{6}(-\exp (-\pi \sqrt{7}))=-1088-768\sqrt{2} 
\]%
\[
\sum_{n=0}^{\infty }\dbinom{2n}{n}^{2}T_{n}(10)\frac{2+6n}{(-64)^{n}}=\frac{%
\sqrt{3}}{\pi } 
\]

\textbf{46.23.}%
\[
J_{6}(-\exp (-\pi \sqrt{10}))=-10304-4608\sqrt{5} 
\]%
\[
\sum_{n=0}^{\infty }\dbinom{2n}{n}^{2}T_{n}(34)\frac{23+60n}{(-64)^{n}}=%
\frac{10}{\pi } 
\]

\textbf{46.24.}%
\[
J_{6}(-\exp (-\pi \sqrt{18}))=-307264-125440\sqrt{6} 
\]%
\[
\sum_{n=0}^{\infty }\dbinom{2n}{n}^{2}T_{n}(194)\frac{59+140n}{(-64)^{n}}=%
\frac{10}{\pi } 
\]

\textbf{46.25.}%
\[
J_{6}(-\exp (-\pi \sqrt{22}))=-1254464-887040\sqrt{2} 
\]%
\[
\sum_{n=0}^{\infty }\dbinom{2n}{n}^{2}T_{n}(394)\frac{142+330n}{(-64)^{n}}=%
\frac{5\sqrt{11}}{\pi } 
\]

\textbf{46.26.}%
\[
J_{6}(-\exp (-\pi \sqrt{28}))=-8290304-3133440\sqrt{7} 
\]%
\[
\sum_{n=0}^{\infty }\dbinom{2n}{n}^{2}T_{n}(898)\frac{35+102n}{(-64)^{n}}=%
\frac{24}{\pi } 
\]

\textbf{46.27.}%
\[
J_{6}(-\exp (-\pi \sqrt{58}))=-12295628864-2283240960\sqrt{29} 
\]%
\[
\sum_{n=0}^{\infty }\dbinom{2n}{n}^{2}T_{n}(39202)\frac{68403+149292n}{%
(-64)^{n}}=\frac{754}{\pi } 
\]

\textbf{Example 47.}

Let%
\[
A_{n}=\dbinom{2n}{n}\sum_{k=0}^{n}\dbinom{n}{k}\dbinom{2k}{k}\dbinom{2n-4k}{%
n-2k} 
\]%
Then%
\[
y_{0}=\sum_{n=0}^{\infty }A_{n}x^{n}=1+4x+60x^{2}+640x^{3}+10780x^{4}+... 
\]%
satisfies%
\[
\theta ^{3}-2x(2\theta +1)(5\theta ^{2}+5\theta +2)-64x^{2}(\theta
+1)(2\theta +1)(2\theta +3)+640x^{3}(2\theta +1)(2\theta +3)(2\theta +5) 
\]%
with%
\[
J_{47}=\frac{1}{q}+6+79q-352q^{2}+1431q^{3}-4160q^{4}+... 
\]%
which is 6A.

\textbf{47.1.}%
\[
J_{47}(-\exp (-\pi \sqrt{\frac{26}{3}}))=-10384 
\]%
\[
\sum_{n=0}^{\infty }A_{n}\frac{2157+19890n}{(-10384)^{n}}=\frac{649^{3/2}}{%
\pi \sqrt{6}} 
\]

\textbf{47.2.}%
\[
J_{47}(-\exp (-\pi \sqrt{\frac{34}{3}}))=-39184 
\]%
\[
\sum_{n=0}^{\infty }A_{n}\frac{22299+235620n}{(-39184)^{n}}=\frac{2449^{3/2}%
}{\pi \sqrt{3}} 
\]

\textbf{47.3.}%
\[
J_{47}(\exp (-\pi \sqrt{\frac{59}{3}}))=1123616 
\]%
\[
\sum_{n=0}^{\infty }A_{n}\frac{23693757+330117390n}{(1123616)^{n}}=\frac{%
4\cdot 70226^{3/2}}{\pi } 
\]

\textbf{47.4.}%
\[
J_{47}(\exp (\frac{2\pi i}{3}-\pi \sqrt{\frac{59}{3}}))=-32+64i 
\]%
\[
\sum_{n=0}^{\infty }A_{n}\frac{18-7i+34(2-i)n}{(-32+64i)^{n}}=\frac{24(2-i)%
\sqrt{10}}{\pi \sqrt{11+2i}} 
\]

\textbf{Example 48.}

Let%
\[
A_{n}=\sum_{k=0}^{n}(-1)^{k}\dbinom{n}{k}\dbinom{n-k}{k}\dbinom{2k}{k}%
\dbinom{3k}{k} 
\]%
Then%
\[
y_{0}=\sum_{n=0}^{\infty }A_{n}x^{n}=1+x-11x^{2}-35x^{3}+469x^{4}+... 
\]%
satisfies%
\[
\theta ^{3}-x(2\theta +1)(2\theta ^{2}+2\theta +1)+x^{2}(\theta
+1)(114\theta ^{2}+228\theta +103) 
\]%
\[
-110x^{3}(\theta +1)(\theta +2)(2\theta +3)+109x^{4}(\theta +1)(\theta
+2)(\theta +3) 
\]%
with%
\[
J_{48}=\frac{1}{q}+1-21q+171q^{2}-745q^{3}+2418q^{4}-... 
\]%
which is not in [8].

\textbf{48.1.}%
\[
J_{48}(\exp (-\pi \sqrt{\frac{49}{12}}))=1+216\sqrt{7} 
\]%
\[
\sum_{n=0}^{\infty }A_{n}\frac{-436654477+12839444353\sqrt{7}%
+(-754427520+81478421640\sqrt{7})n}{(1+216\sqrt{7})^{n}}=\frac{(19\cdot
17189)^{2}}{\pi } 
\]

\textbf{48.2.}%
\[
J_{48}(\exp (-\pi \sqrt{\frac{89}{12}}))=1+3000\sqrt{3} 
\]%
\[
\sum_{n=0}^{\infty }A_{n}\frac{426735012497-133974079944000\sqrt{3}%
+(764154000000-1146231042453000\sqrt{3})n}{(1+3000\sqrt{3})^{n}} 
\]%
\[
=\frac{(13\cdot 23\cdot 73\cdot 1237)^{2}}{\pi } 
\]

\textbf{48.3.}%
\[
J_{48}(i\exp (-\pi \sqrt{\frac{80}{12}}))=1+15i\sqrt{15} 
\]%
\[
\sum_{n=0}^{\infty }A_{n}\frac{5i+24\sqrt{15}+99\sqrt{15}n}{(1+15i\sqrt{15}%
)^{n}}=\frac{-150i\sqrt{3}+3374\sqrt{5}}{25\pi } 
\]

\textbf{Example 49.}

Let%
\[
A_{n}=\sum_{k=0}^{n}(-1)^{k}\dbinom{n}{k}\dbinom{n-k}{k}\dbinom{3k}{k}%
\dbinom{6k}{3k} 
\]%
Then%
\[
y_{0}=\sum_{n=0}^{\infty }A_{n}x^{n}=1+x-119x^{2}-359x^{3}+82441x^{4}+... 
\]%
satisfies%
\[
\theta ^{3}-x(2\theta +1)(2\theta ^{2}+2\theta +1)+x^{2}(\theta
+1)(1734\theta ^{2}+3468\theta +967) 
\]%
\[
-1730x^{3}(\theta +1)(\theta +2)(2\theta +3)+1729x^{4}(\theta +1)(\theta
+2)(\theta +3) 
\]%
with%
\[
J_{49}=\frac{1}{q}+1-372q+29250q^{3}+134120q^{5}+... 
\]%
We have the relation%
\[
(J_{42}(q^{1/2})-1)^{2}=\frac{1}{q}-42+783q-8672q^{2}+65367q^{3}+... 
\]%
which is 3A.

\textbf{49.1.}%
\[
J_{49}(i\exp (-\pi \sqrt{2}))=1+40i\sqrt{5} 
\]%
\[
\sum_{n=0}^{\infty }A_{n}\frac{155416i+1096000\sqrt{5}+(256000i+5119360\sqrt{%
5})n}{(1+40i\sqrt{5})^{n}}=\frac{3^{4}\cdot 7\cdot 127^{2}}{\pi } 
\]

\textbf{49.2.}%
\[
J_{49}(i\exp (-\pi \sqrt{3}))=1+60i\sqrt{15} 
\]%
\[
\sum_{n=0}^{\infty }A_{n}\frac{25271676i+233236800\sqrt{15}%
+(42768000i+1283016240\sqrt{15})n}{(1+60i\sqrt{15})^{n}}=\frac{54001^{2}}{%
\pi } 
\]

\textbf{49.3.}%
\[
J_{49}(i\exp (-2\pi ))=1+66i\sqrt{66} 
\]%
\[
\sum_{n=0}^{\infty }A_{n}\frac{1272i+15840\sqrt{66}+99792\sqrt{66}n}{(1+66i%
\sqrt{66})^{n}}=\frac{287495\sqrt{2}-264i\sqrt{33}}{\pi } 
\]

\textbf{49.4.}%
\[
J_{49}(i\exp (-\pi \sqrt{7}))=1+255i\sqrt{255} 
\]%
\[
\sum_{n=0}^{\infty }A_{n}\frac{752333103i+20600780625\sqrt{255}%
+(1343091375i+171244139985\sqrt{255})n}{(1+255i\sqrt{255})^{n}}=\frac{%
2^{15}\cdot 7\cdot 19\cdot 487^{2}}{\pi } 
\]

\textbf{Example 50.}

Let%
\[
A_{n}=\sum_{k=0}^{\infty }(-1)^{k}\dbinom{n}{3k}\dbinom{2k}{k}^{3} 
\]%
Then%
\[
y_{0}=\sum_{n=0}^{\infty }A_{n}x^{n}=1+x+x^{2}-7x^{3}-31x^{4}-79x^{5}+... 
\]%
satisfies%
\[
\theta ^{3}-x(6\theta ^{3}+6\theta ^{2}+4\theta +1)+3x^{2}(5\theta
^{3}+10\theta ^{2}+9\theta +3) 
\]%
\[
+x^{3}(44\theta ^{3}+228\theta ^{2}+364\theta +189)-3x^{4}(\theta
+1)(59\theta ^{2}+241\theta +249) 
\]%
\[
+6x^{5}(\theta +1)(\theta +2)(31\theta +78)-63x^{6}(\theta +1)(\theta
+2)(\theta +3) 
\]%
with%
\[
J_{50}=\frac{1}{q}+1-8q^{2}+28q^{5}-64q^{8}+134q^{11}-... 
\]%
which is 6F.

\textbf{50.1.}%
\[
J_{50}(\exp (-\pi \sqrt{\frac{2}{9}}))=5 
\]%
\[
\sum_{n=0}^{\infty }A_{n}\frac{11+16n}{5^{n}}=\frac{75}{2\pi } 
\]

\textbf{50.2.}%
\[
J_{50}(\exp (-\pi \sqrt{\frac{4}{9}}))=9 
\]%
\[
\sum_{n=0}^{\infty }A_{n}\frac{7+16n}{9^{n}}=\frac{81\sqrt{2}}{4\pi } 
\]

\textbf{50.3.}%
\[
J_{50}(-\exp (-\pi \sqrt{\frac{7}{9}}))=-15 
\]%
\[
\sum_{n=0}^{\infty }A_{n}\frac{89+224n}{(-15)^{n}}=\frac{225}{\pi } 
\]

\textbf{50.4.}%
\[
J_{50}(\exp (\frac{\pi i}{3}-\pi \sqrt{\frac{4}{9}}))=-3+4i\sqrt{3} 
\]%
\[
\sum_{n=0}^{\infty }A_{n}\frac{14i-18\sqrt{3}+32(i-\sqrt{3})n}{(-3+4i\sqrt{3}%
)^{n}}=\frac{39i\sqrt{2}-24\sqrt{6}}{\pi } 
\]

\textbf{50.5.}%
\[
J_{50}(\exp (\frac{\pi i}{3}-\pi \sqrt{\frac{7}{9}}))=9+8i\sqrt{3} 
\]%
\[
\sum_{n=0}^{\infty }A_{n}\frac{1931+235i\sqrt{3}+(5792+352i\sqrt{3})n}{(9+8i%
\sqrt{3})^{n}}=\frac{7098}{\pi } 
\]

\textbf{Exampe 51.}

Let%
\[
A_{n}=\sum_{k=0}^{n}(-1)^{k}\dbinom{n}{3k}\dbinom{2k}{k}^{2}\dbinom{3k}{k} 
\]%
Then

\[
y_{0}=\sum_{n=0}^{\infty }A_{n}x^{n}=1+x+x^{2}-11x^{3}-47x^{4}-119x^{5}+... 
\]%
satisfies%
\[
\theta ^{3}-x(\theta +1)(2\theta ^{2}+2\theta +1)+x^{2}(\theta +1)(6\theta
^{2}+12\theta +7) 
\]%
\[
+52x^{3}(\theta +1)(\theta +2)(2\theta +3)-107x^{4}(\theta +1)(\theta
+2)(\theta +3) 
\]%
with%
\[
J_{51}=\frac{1}{q}+1-14q^{2}+65q^{5}-156q^{8}+456q^{11}-1066q^{14}+... 
\]%
We have the transformation%
\[
(J_{44}(q^{1/3})-1)^{3}=\frac{1}{q}-42+783q-8672q^{2}+65367q^{3}-... 
\]%
which is 3A.

\textbf{51.1.}%
\[
J_{51}(\exp (-\pi \sqrt{\frac{41}{27}}))=49 
\]%
\[
\sum_{n=0}^{\infty }A_{n}\frac{399+1230n}{7^{2n}}=\frac{7^{4}}{4\pi } 
\]

\textbf{51.2.}%
\[
J_{51}(\exp (-\pi \sqrt{\frac{89}{27}}))=301 
\]%
\[
\sum_{n=0}^{\infty }A_{n}\frac{48842+283020n}{301^{n}}=\frac{301^{2}}{\pi } 
\]

\textbf{51.3.}%
\[
J_{51}(-\exp (-\pi \sqrt{\frac{16}{27}}))=1-9\cdot 2^{1/3} 
\]%
\[
\sum_{n=0}^{\infty }A_{n}\frac{632664+23332\cdot 2^{1/3}+122568\cdot
2^{2/3}+(1576800+29160\cdot 2^{1/3}+175020\cdot 2^{2/3})n}{(1-9\cdot
2^{1/3})^{n}}=\frac{(31\cdot 47)^{2}}{\pi } 
\]

\textbf{51.4.}%
\[
J_{51}(\exp (-\pi \sqrt{\frac{49}{27}}))=1+36\cdot 7^{1/3} 
\]%
\[
\sum_{n=0}^{\infty }A_{n}\frac{98-8148\cdot 7^{1/3}+78624\cdot
7^{2/3}+(-4620\cdot 7^{1/3}+332640\cdot 7^{2/3})n}{(1+36\cdot 7^{1/3})^{n}} 
\]%
\[
=\frac{653183\cdot 7^{1/6}+3888\cdot 7^{5/6}}{\pi } 
\]

\textbf{Example 52.}

Let%
\[
A_{n}=\sum_{j,k}\dbinom{n}{2j}\dbinom{j}{k}^{2}\dbinom{2k}{k}\dbinom{2j-2k}{%
j-k} 
\]%
Then

\[
y_{0}=\sum_{n=0}^{\infty }A_{n}x^{n}=1+x+5x^{2}+13x^{3}+53x^{4}+181x^{5}+... 
\]%
satisfies a differential equation of order three and degree sev\NEG{e}n
which we delete. We have%
\[
J_{52}=\frac{1}{q}+1+3q+3q^{3}+7q^{5}+18q^{7}+21q^{9}+... 
\]%
which is 24A.

\textbf{52.1.}%
\[
J_{52}(\exp (-\pi \sqrt{\frac{5}{12}}))=9 
\]%
\[
\sum_{n=0}^{\infty }A_{n}\frac{13+40n}{9^{n}}=\frac{54\sqrt{3}}{\pi } 
\]

\textbf{Example 53.}

Let%
\[
A_{n}=\sum_{j,k}\dbinom{n}{3j}\dbinom{j}{k}^{2}\dbinom{2k}{k}\dbinom{2j-2k}{%
j-k} 
\]%
Then%
\[
y_{0}=\sum_{n=0}^{\infty }A_{n}x^{n}=1+x+x^{2}+5x^{3}+17x^{4}+41x^{5}+... 
\]%
satisfies a differential equation of degree nine with%
\[
J_{53}=\frac{1}{q}+1+2q^{2}+q^{5}+4q^{8}+8q^{11}+6q^{14}+... 
\]%
We have the identity%
\[
(J_{46}(q^{1/3})-1)^{3}=\frac{1}{q}+6+15q+32q^{2}+87q^{3}+192q^{4}+... 
\]%
which is 6C.

\textbf{53.1.}%
\[
J_{53}(\exp (-\pi \sqrt{\frac{5}{27}}))=5 
\]%
\[
\sum_{n=0}^{\infty }A_{n}\frac{1+2n}{5^{n}}=\frac{5\sqrt{3}}{\pi } 
\]

\textbf{Example 54.}

Let%
\[
A_{n}=\sum_{j,k}(-1)^{k}3^{j-3k}\dbinom{n}{2j}\dbinom{j}{3k}\dbinom{j+k}{j}%
\dbinom{2k}{k}\dbinom{3k}{k} 
\]%
Then%
\[
J_{54}=\frac{1}{q}+1+2q-3q^{3}-8q^{5}-2q^{7}+6q^{9}+... 
\]%
which is 12G.

\textbf{54.1.}%
\[
J_{54}(\exp (-\pi \sqrt{\frac{31}{12}}))=79+45\sqrt{3} 
\]%
\[
\sum_{n=0}^{\infty }A_{n}\frac{2425371-1396077\sqrt{3}+(173166-78120\sqrt{3}%
)n}{(79+45\sqrt{3})^{n}}=\frac{13778\sqrt{3}}{\pi } 
\]

\textbf{Example 55.}

Let

\[
A_{n}=\sum_{j,k}(-1)^{j+k}3^{j-3k}\dbinom{n}{4j}\dbinom{j}{3k}\dbinom{j+k}{j}%
\dbinom{2k}{k}\dbinom{3k}{k} 
\]%
Then%
\[
J_{55}=\frac{1}{q}+1-q^{3}-2q^{7}+2q^{11}-q^{15}+5q^{19}+... 
\]%
which is 24G. The table in [13] gives only terms up to $q^{10}$ but higher
terms can be obtained from the product%
\[
J_{55}=1+\frac{\eta (q^{4})\eta (q^{8})}{\eta (q^{12})\eta (q^{24})} 
\]%
where 
\[
\eta (q)=q^{1/24}\dprod\limits_{n=1}^{\infty }(1-q^{n}) 
\]

\textbf{55.1.}%
\[
J_{55}(\exp (-\pi \sqrt{\frac{1}{8}}))=4 
\]%
\[
\sum_{n=0}^{\infty }A_{n}\frac{1+n}{4^{n}}=\frac{8}{\pi \sqrt{3}} 
\]

\textbf{Example 56.}

Let

\[
A_{n}=\sum_{j,k}(-1)^{j}\dbinom{n}{2j}\dbinom{j}{k}^{2}\dbinom{2j}{j}\dbinom{%
j+k}{j} 
\]%
Then%
\[
J_{56}=\frac{1}{q}+1-8q+35q^{3}-100q^{5}+260q^{7}-548q^{9}+... 
\]%
which is not in [13] but%
\[
(J_{56}(q^{1/2})-1)^{2}=\frac{1}{q}-16+134q-760q^{2}+3345q^{3}+... 
\]%
which is 5A.

\textbf{56.1.}%
\[
J_{56}(i\exp (-\pi \sqrt{\frac{3}{5}}))=1+7i\sqrt{3} 
\]%
\[
\sum_{n=0}^{\infty }A_{n}\frac{785i+1425\sqrt{3}+(1155i+4015\sqrt{3})n}{(1+7i%
\sqrt{3})^{n}}=\frac{8\cdot 37^{2}}{\pi } 
\]

\textbf{Example 57.}

Let

\[
A_{n}=\sum_{j,k}(-1)^{j}\dbinom{n}{2j}\dbinom{j}{k}^{2}\dbinom{2j}{j}\dbinom{%
2k}{k} 
\]%
Then%
\[
J_{57}=\frac{1}{q}+1-7q+15q^{3}-71q^{5}+106q^{7}-273q^{9}+... 
\]%
which is not in [13] but%
\[
(J_{57}(q^{1/2})-1)^{2}=\frac{1}{q}-14+79q-352q^{2}+1431q^{3}+... 
\]%
which is 6A.

\textbf{57.1.}%
\[
J_{57}(i\exp (-\pi \sqrt{\frac{17}{6}}))=1+198i 
\]%
\[
\sum_{n=0}^{\infty }A_{n}\frac{4128+8724i+(22848+45696i)n}{(1+198i)^{n}}=%
\frac{7999+15602i}{\pi } 
\]

\textbf{Example 58.}

Let

\[
A_{n}=\sum_{j,k}(-1)^{j}\dbinom{n}{2j}\dbinom{j}{k}\dbinom{2j}{j}\dbinom{2k}{%
k}\dbinom{2j-2k}{j-k} 
\]%
Then%
\[
J_{58}=\frac{1}{q}+1-8q-6q^{3}-48q^{5}+15q^{7}-168q^{9}+... 
\]%
which is not in [13] but%
\[
(J_{58}(q^{1/2})-1)^{2}=\frac{1}{q}-16+52q+834q^{3}+4760q^{5}+... 
\]%
which is 4B.

\textbf{58.1.}%
\[
J_{51}(i\exp (-\pi \sqrt{\frac{9}{8}}))=1+40i 
\]%
\[
\sum_{n=0}^{\infty }A_{n}\frac{488+34i+1848n}{(1+40i)^{n}}=\frac{1599-80i}{%
\pi } 
\]

\textbf{Example 59.}

Let

\[
A_{n}=\sum_{j,k}(-1)^{j}\dbinom{n}{4j}\dbinom{j}{k}\dbinom{2j}{j}\dbinom{2k}{%
k}\dbinom{2j-2k}{j-k} 
\]%
Then%
\[
J_{59}=\frac{1}{q}+1+4q^{3}-11q^{7}+68q^{11}-325q^{15}+2132q^{19}+... 
\]%
which is not in [13] but%
\[
(J_{59}(q^{1/4})-1)^{4}=\frac{1}{q}+16+52q+834q^{3}+4760q^{5}+... 
\]%
which is 4B.

\textbf{59.1.}%
\[
J_{52}(i\exp (-\pi \sqrt{\frac{11}{32}}))=1+\sqrt{40} 
\]%
\[
\sum_{n=0}^{\infty }A_{n}\frac{34800-4743\sqrt{10}+(63140-6160\sqrt{10})n}{%
(1+\sqrt{40})^{n}}=\frac{24\cdot 65^{2}}{\pi } 
\]

\textbf{Example 60.}

Let

\[
A_{n}=\sum_{k=0}^{n}\dbinom{n}{6k}\dbinom{2k}{k}^{2}\dbinom{3k}{k} 
\]%
Then%
\[
J_{60}=\frac{1}{q}+1+7q^{5}+8q^{11}+22q^{17}+42q^{23}+... 
\]%
which is not in [13] but%
\[
(J_{60}(q^{1/6})-1)^{6}=\frac{1}{q}+42+783q+8672q^{2}+65367q^{3}+... 
\]%
which is 3A.

\textbf{60.1.}%
\[
J_{60}(\exp (-\pi \sqrt{\frac{5}{27}}))=1+\sqrt{15} 
\]%
\[
\sum_{n=0}^{\infty }A_{n}\frac{335-48\sqrt{15}+(440-55\sqrt{15})n}{(1+\sqrt{%
15})^{n}}=\frac{490\sqrt{3}}{\pi } 
\]

\textbf{Example 61.}

Let

\[
A_{n}=\sum_{k=0}^{n}(-1)^{k}\dbinom{n}{6k}\dbinom{2k}{k}^{2}\dbinom{3k}{k} 
\]%
Then%
\[
J_{61}=\frac{1}{q}+1-7q^{5}+8q^{11}-22q^{17}+42q^{23}-... 
\]

\textbf{54.1.}

\[
J_{61}(-i\exp (-\pi \sqrt{\frac{5}{27}}))=1+i\sqrt{15} 
\]%
\[
\sum_{n=0}^{\infty }A_{n}\frac{265+47i\sqrt{15}+(385+55i\sqrt{15})n}{(1+i%
\sqrt{15})^{n}}=\frac{640\sqrt{3}}{\pi } 
\]

\textbf{Example 62.}

Let

\[
A_{n}=\sum_{k=0}^{n}\dbinom{n}{8k}\dbinom{2k}{k}^{2}\dbinom{4k}{2k} 
\]%
Then%
\[
J_{62}=\frac{1}{q}+1+13q^{7}-45q^{15}+748q^{23}-10359q^{31}+... 
\]%
which is not in [13] but%
\[
(J_{62}(q^{1/8})-1)^{8}=\frac{1}{q}+104+4372q+96256q^{2}+1240002q^{3}+... 
\]%
which is 2A.

\textbf{62.1.}%
\[
J_{62}(-i\exp (-\pi \sqrt{\frac{29}{32}}))=1+6i\sqrt{11} 
\]%
\[
\sum_{n=0}^{\infty }A_{n}\frac{227722968+13927178i\sqrt{11}%
+(687987300+20900880i\sqrt{11})n}{(1+6i\sqrt{11})^{n}}=\frac{27\cdot
(11\cdot 397)^{2}\sqrt{2}}{\pi } 
\]

\textbf{Example 63.}

Let%
\[
A_{n}=\sum_{k=0}^{n}\dbinom{n}{2k}a_{k} 
\]%
where%
\[
y_{0}=\sum_{k=0}^{\infty }a_{k}x^{k} 
\]%
is the solution of Example 20 (no formula for $a_{k}$ is known). Then 
\[
J_{63}=\frac{1}{q}+1-2q+2q^{3}-7q^{5}+5q^{7}-11q^{9}+21q^{11}-... 
\]%
with%
\[
(J_{63}(q^{1/2})-1)^{2}=\frac{1}{q}%
-4+8q-22q^{2}+42q^{3}-70q^{4}+155q^{5}-... 
\]%
which is 15A.

\textbf{63.1.}%
\[
J_{63}(\exp (-\pi \sqrt{\frac{37}{60}}))=1+3\sqrt{15} 
\]%
\[
\sum_{n=0}^{\infty }A_{n}\frac{-25183+23490\sqrt{15}+(-34965+52836\sqrt{15})n%
}{(1+3\sqrt{15})^{n}}=\frac{54\cdot 67^{2}}{\pi } 
\]

\textbf{63.2.}%
\[
J_{63}(\exp (-\pi \sqrt{\frac{53}{60}}))=1+11\sqrt{3} 
\]%
\[
\sum_{n=0}^{\infty }A_{n}\frac{14117796-1626937\sqrt{3}+(39789750-2404875%
\sqrt{3})n}{(1+11\sqrt{3})^{n}}=\frac{6\cdot (11\cdot 181)^{2}\sqrt{3}}{\pi }
\]

\textbf{Example 64.}

Let

\[
A_{n}=\sum_{j,k}\dbinom{n}{2j}\dbinom{2k}{k}^{3}\dbinom{2j-4k}{j-2k} 
\]%
Then%
\[
J_{64}=\frac{1}{q}+1+6q^{3}+15q^{7}+26q^{11}+51q^{15}+102q^{19}... 
\]%
which is not in [13] but%
\[
(J_{64}(q^{1/4})-1)^{4}=\frac{1}{q}+24+276q+2048q^{2}+11202q^{3}+... 
\]%
which is 2B.

\textbf{64.1.}%
\[
J_{64}(\exp (-\pi \sqrt{\frac{7}{16}}))=9 
\]%
\[
\sum_{n=0}^{\infty }A_{n}\frac{47+140n}{9^{n}}=\frac{72\sqrt{15}}{\pi } 
\]

\textbf{64.2.}%
\[
J_{64}(i\exp (-\pi \sqrt{\frac{7}{16}}))=1+8i 
\]%
\[
\sum_{n=0}^{\infty }A_{n}\frac{374432+139624i+(719712+182784i)n}{(1+8i)^{n}}=%
\frac{520^{2}\sqrt{17}}{\pi } 
\]

\textbf{Example 65.}

Let%
\[
A_{n}=\sum_{j,k}\dbinom{n}{2j}\dbinom{2j}{j}\dbinom{j}{k}^{2}\dbinom{j+k}{j} 
\]%
Then

\[
J_{65}=\frac{1}{q}+1+8q+35q^{3}+100q^{5}+260q^{7}+548q^{9}... 
\]%
which is not in [13] but%
\[
(J_{65}(q^{1/2})-1)^{2}=\frac{1}{q}+16+134q+760q^{2}+3345q^{3}+... 
\]%
which is 5A

\textbf{65.1}%
\[
J_{65}(i\exp (-\pi \sqrt{\frac{47}{20}}))=1+18i\sqrt{47} 
\]%
\[
\sum_{n=0}^{\infty }A_{n}\frac{3485370i+10805770\sqrt{47}+(5769720i+51924070%
\sqrt{47})n}{(1+18i\sqrt{47})^{n}}=\frac{(97\cdot 157)^{2}}{\pi } 
\]

\textbf{Example 66.}

Let%
\[
A_{n}=\sum_{j,k}(-1)^{j}\dbinom{n}{3j}\dbinom{2k}{j}\dbinom{j}{k}^{2}\dbinom{%
j+k}{j} 
\]%
Then%
\[
J_{66}=\frac{1}{q}+1-3q^{2}+8q^{5}-11q^{8}+25q^{11}-35q^{14}+... 
\]%
which is 21C.

\textbf{66.1.}%
\[
J_{66}(\exp (-\pi \sqrt{\frac{61}{63}}))=23 
\]%
\[
\sum_{n=0}^{\infty }A_{n}\frac{25613+87230n}{23^{n}}=\frac{(22\cdot 23)^{2}}{%
\pi \sqrt{7}} 
\]

\textbf{66.2.}%
\[
J_{66}(\exp (\frac{2\pi i}{3}-\pi \sqrt{\frac{61}{63}}))=-10+11i\sqrt{3} 
\]%
\[
\sum_{n=0}^{\infty }A_{n}\frac{190239553+17573591i\sqrt{3}%
+(563593030+26256230i\sqrt{3})n}{(-10+11i\sqrt{3})^{n}}=\frac{8\cdot
(11\cdot 463)^{2}\sqrt{7}}{\pi } 
\]

\textbf{Example 67.}

Let%
\[
A_{n}=\sum_{j,k}(-1)^{j+k}\frac{j-2k}{2j-3k}\dbinom{n}{2j}\dbinom{2k}{k}%
\dbinom{j}{k}\dbinom{2j-2k}{j-k}\dbinom{2j-3k}{j} 
\]%
Then

\[
J_{67}=\frac{1}{q}+1-4q+4q^{5}+16q^{7}+6q^{11}-40q^{13}+... 
\]%
which is not in [13] but%
\[
(J_{67}(q^{1/2})-1)^{2}=\frac{1}{q}-8+16q+8q^{2}-128q^{3}+...=J_{16} 
\]

\textbf{67.1.}%
\[
J_{67}(\exp (-\pi \sqrt{\frac{37}{36}}))=25 
\]%
\[
\sum_{n=0}^{\infty }A_{n}\frac{3757+12432n}{25^{n}}=\frac{7500\sqrt{3}}{\pi }
\]

\textbf{67.2.}%
\[
J_{67}(\exp (-\pi \sqrt{\frac{22}{36}}))=1+12i 
\]%
\[
\sum_{n=0}^{\infty }A_{n}\frac{180+29i+528n}{(1+12i)^{n}}=\frac{(429-72i)%
\sqrt{3}}{\pi } 
\]

\textbf{Example 68.}

Let%
\[
A_{n}=\sum_{i,j,k}\dbinom{n}{3i}\dbinom{i}{j}^{2}\dbinom{i}{k}\dbinom{j}{k}%
\dbinom{j+k}{i} 
\]%
Then%
\[
y_{0}=\sum_{n=0}^{\infty
}A_{n}x^{n}=1+x+x^{2}+4x^{3}+13x^{4}+31x^{5}+88x^{6}+295x^{7}+... 
\]%
satisfies a differential equation of order three and degree 9.We have%
\[
J_{68}=\frac{1}{q}+1+2q^{2}+5q^{5}+6q^{8}+12q^{11}+16q^{14}+27q^{17}+... 
\]%
with%
\[
(J_{68}(q^{1/3})-1)^{3}=\frac{1}{q}+6+27q+86q^{2}+243q^{3}+594q^{4}+... 
\]%
which is 9A.

\textbf{68.1.}%
\[
J_{68}(\exp (-\pi \sqrt{\frac{8}{81}}))=4 
\]%
\[
\sum_{n=0}^{\infty }A_{n}\frac{2+3n}{4^{n}}=\frac{18\sqrt{3}}{\pi } 
\]

This formula is not so easy to find. The reason is that the convergence is
very slow. If we expand to degre $100$ we find $\dfrac{A_{100}}{4^{n}}%
\approx 10^{-7}$ which is not good enough to use PSLQ. Instead we use the
formulas in [6] p.4. To get $b$ in%
\[
\sum_{n=0}^{\infty }A_{n}\frac{a+bn}{4^{n}}=\frac{1}{\pi } 
\]%
we need $P$, which in this case is not a polynomial. We need $e_{3}$ in%
\[
\theta ^{3}y_{0}=(e_{3}(x)\theta ^{2}+e_{1}(x)\theta +e_{0}(x))y_{0} 
\]%
If the differential operator is 
\[
L=\theta ^{3}+\sum_{j=1}^{9}x^{j}p_{j}(\theta ) 
\]%
then we have (using Maple notation)%
\[
e_{3}=-\frac{\sum_{j=1}^{9}coeff(p_{j},\theta ,2)\cdot x^{j}}{%
1+\sum_{j=1}^{9}coeff(p_{j},\theta ,3)\cdot x^{j}}=\frac{%
3x(2-12x+57x^{2}-157x^{3}+219x^{4}-66x^{5}-16x^{6})}{%
(1-x)(1-6x+15x^{2}-38x^{3}+69x^{4}-60x^{5}-8x^{6})} 
\]%
and%
\[
P=\exp (\int (-\frac{2e_{3}}{3x})dx)=\frac{%
1-6x+15x^{2}-38x^{3}+69x^{4}-60x^{5}-8x^{6}}{(1-x)^{2}} 
\]%
We obtain%
\[
b=\sqrt{\frac{8}{81}}\sqrt{P(\frac{1}{4})}=\frac{1}{6\sqrt{3}} 
\]%
To get $a$ we need to expand everything in $q$ and then put \ $q=q_{0}=\exp
(-\pi \sqrt{\dfrac{8}{81}})=0.372578...$We have%
\[
x=x(q)=\frac{1}{J_{68}}%
=q-q^{2}+q^{3}-3q^{4}+5q^{5}-7q^{6}+8q^{7}-13q^{8}+... 
\]%
and%
\[
y_{00}=y_{0}(x(q))=1+q+3q^{3}+q^{4}+9q^{6}+...+144q^{99}-2q^{100}+O(q^{101}) 
\]%
has very small coefficients compared with $y_{0}(x).$ We get%
\[
a=\frac{1}{\pi y_{00}}\left\{ 1-\frac{\pi q}{y_{00}}\sqrt{\frac{8}{81}}%
\dfrac{dy_{00}}{dq}\right\} _{\mid q=q_{0}}=0.0641500299099... 
\]%
and%
\[
\frac{1}{a^{2}}=242.999999999999999999999999999999999999... 
\]%
so%
\[
a=\frac{1}{9\sqrt{3}} 
\]

\textbf{Example 69.}

Consider the polytope (equivalent to the hexagon B7 in [18]) with Laurent
polynomial%
\[
S=x+y+\frac{1}{x}+\frac{1}{y}+\frac{x}{y}+\frac{y}{x} 
\]%
Then%
\[
u_{0}=\sum_{n=0}^{\infty
}c.t.(S^{n})x^{n}=1+6x^{2}+12x^{3}+90x^{4}+360x^{5}+... 
\]%
satisfies%
\[
\theta ^{2}-x\theta (\theta +1)-24x^{2}(\theta +1)^{2}-36x^{3}(\theta
+1)(\theta +2) 
\]%
with 
\[
J_{69}=\frac{1}{q}+1+6q+4q^{2}-3q^{3}-12q^{4}-8q^{5}+12q^{6}+... 
\]%
which is 6E. Then 
\[
y_{0}=u_{0}^{2}=\sum_{n=0}^{\infty
}A_{n}x^{n}=1+12x^{2}+24x^{3}+216x^{3}+864x^{4}+5304x^{5}+... 
\]%
satisfies%
\[
\theta ^{3}-x\theta (\theta +1)(2\theta +1)-x^{2}(\theta +1)(47\theta
^{2}+94\theta +96)-12x^{3}(2\theta +3)(\theta ^{2}+3\theta +8) 
\]%
\[
+648x^{4}(\theta +2)^{3}+864x^{5}(\theta +2)(\theta +3)(2\theta
+5)+1296x^{6}(\theta +2)(\theta +3)(\theta +4) 
\]

\textbf{69.1.}%
\[
J_{69}(\exp (-\pi \sqrt{\frac{4}{3}}))=18+12\sqrt{3} 
\]%
\[
\sum_{n=0}^{\infty }A_{n}\frac{1+(2+\sqrt{3})n}{(18+12\sqrt{3})^{n}}=\frac{27%
}{8\pi } 
\]

\textbf{69.2.}%
\[
J_{69}(-\exp (-\pi \sqrt{\frac{5}{3}}))=-30-12\sqrt{5} 
\]%
\[
\sum_{n=0}^{\infty }A_{n}\frac{1+(2+\sqrt{5})n}{(-30-12\sqrt{5})^{n}}=\frac{%
15\sqrt{3}}{8\pi } 
\]

\textbf{69.3.}%
\[
J_{69}(-\exp (-\pi \sqrt{3}))=-78-60\cdot 2^{1/3}-48\cdot 2^{2/3} 
\]%
\[
\sum_{n=0}^{\infty }A_{n}\frac{-144+252\cdot 2^{1/3}-66\cdot
2^{2/3}+(18+156\cdot 2^{1/3}+102\cdot 2^{2/3})n}{(-78-60\cdot
2^{1/3}-48\cdot 2^{2/3})^{n}}=\frac{125\sqrt{3}}{\pi } 
\]%
\[
\]

\textbf{Final remark.}

We have shown how Ramanujan-like formulas for $\dfrac{1}{\pi }$ are
connected to polytopes, K3-surfaces and "Moonshine". In many cases we have
no explicit formula for the coefficients, but the fastest way to compute
them is by solving the differential equation anyway.

In [16], Kreuzer and Skarke have listed all the 4319 reflexive polytopes in
dimension 3. Unfortunately only~the two with Picard number 19 are useful for
us.

In [15], Lian and Wieczer have computed $\ Q(x)$ \ such that one can find
the J-function to high degree in all 157 cases in [13] by solving%
\[
u^{\prime \prime }+Q(x)u=0 
\]%
E.g. in case 17A we have

\[
Q(x)=\frac{%
1+8x+58x^{2}+864x^{3}+7265x^{4}+30984x^{5}+79840x^{6}+84016x^{7}+41200x^{8}}{%
4x^{2}(1+4x-35x^{2}-254x^{3}-564x^{4}-424x^{5})^{2}} 
\]%
and we compute%
\[
J=\frac{1}{q}+7q+14q^{2}+29q^{3}+50q^{4}+92q^{5}+... 
\]%
We find%
\[
J(-\exp (-\pi \sqrt{\frac{11}{17}})=-13\text{ \ and \ }J(-\exp (-\pi \sqrt{%
\frac{83}{17}})=-515-126\sqrt{17} 
\]%
but we do not know any \ $A_{n}$ to use (you cannot use \ $y_{0}=u_{0}^{2}$
). So an interesting problem is to find the right third order differential
equations giving the various $J-$functions.

\textbf{Acknowledgements.}

I want to thank Jesus Guillera, who insisted that his formulas are correct.
\ 

\textbf{References.}

\textbf{1. }G.Almkvist, W.Zudilin, Differential equations, mirror maps and
zeta values,in \textit{Mirror symmetry V (}ed. N.Yui, S.-T.Yau and
J.D.Lewis), AMS/IP Studies in Advanced Mathematics, Volume 38, pp. 481-515
(AMS and Intern. Press, Providence, RI, 2006).

\textbf{2. }G.Almkvist, D.van Straten, W.Zudilin, Generalizations of
Clausen's formula and algebraic transformations of Calabi-Yau differential
equations, Proc. of the Edinburgh Math. Soc. 54 (2011), 273-295.

\textbf{3. }G.Almkvist, C.van Enckevort, D.van Straten, W.Zudilin, Tables of
Calabi-Yau equations, arXiv, AG/0507430.

\textbf{4. }G.Almkvist, A.Ackroy, appendix by A.Meurman, Proof of some
conjectured formulas for $\dfrac{1}{\pi }$ by Z.W.Sun,

arXiv NT/1112.3259.

\textbf{5. }G.Almkvist, J.Guillera, Ramanujan-like series for $\dfrac{1}{\pi
^{2}}$ and String Theory, to appear in Exp. Math., arXiv NT-1009.5202.

\textbf{6. }G.Almkvist, J.Guillera, Ramanujan-Sato-like series, arXiv,
NT/1201.5233

\textbf{7. }G.Almkvist, Str\"{a}ngar i m\aa nsken II, Normat 51 (2003),
63-79.

\textbf{8. }G.Almkvist, Ramanujan-like series for $\dfrac{1}{\pi ^{2}}$ \'{a}
la Guillera and Zudilin and Calabi-Yau differential equations, Computer
Science J. of Moldova, 17 (2009), 1-21.

\textbf{9. }N.D.Baruah, B.C.Berndt, Ramanujan's Eisenstein series and new
hypergeometric-like series for $\dfrac{1}{\pi ^{2}}$ , J. of Approximation
Theory, 160 (2008), 135-153.

\textbf{10. }M.Bogner, Differentielle Galoisgruppen und
Transformationstheorie fuer Calabi-Yau-Operatoren vierter Ordnung,
Diplomarbeit, Gutenberg-Universit\"{a}t Mainz, 2008.

\textbf{11. }H.H.Chan, Y.Tanigawa, Y.Yang, W.Zudilin, New analogues of
Clausen's identities .arising from the theory of \ \ modular forms, Adv.
Math., 228 (2011), 1294-1314.

\textbf{12. }H.H.Chan, J.Wan, W.Zudilin,Legendre polynomials and
Ramanujan-like series for $\dfrac{1}{\pi }$, Israel J.Math. (to appear).

\textbf{13. }J.H.Conway, S.P.Norton, "Monstrous Moonshine", Bull. of London
Math. Soc. 11 (1979), 308-339.

\textbf{14.} B.H.Lian, S.-T.Yau, Mirror maps, modular relations and
hypergeometric series I, Nuclear Phys. 176 (1991), 248-262, arXiv
hep-th/9507151.

\textbf{15.} B.H.Lian, J.L.Wiczes, Genus zero modular functions, arXiv
NT/0611291.

\textbf{16. }M.Kreuzer, H.Skarke, List of 4319 reflexive polytopes in
dimension 3, http//tph16.tuwien.ac.at/ \symbol{126}kreuzer/CY.html

\textbf{17. }K.Samol, D.van Straten, Frobenius polynomials for Calabi-Yau
equations, arXiv AG/0802-3994.

\textbf{18. }J.Stienstra, Hypergeometric systems in two variables, quivers,
dimers and dessins d'enfants, in Modular Forms and String Duality, Banff
2006, N.Yui, H.Verrill, C.F.Doran eds, Fields Inst. Comm. No. 54, 2008.

\textbf{19. }Z.W.Sun, List of conjectural series for powers of $\pi $ and
other constants, arXiv NT/11025649

\textbf{20. }Z.W.Sun, Some new series for $\dfrac{1}{\pi }$ and related
congruences, arXiv NT/11043856.

\textbf{21. }J.Wan, W.Zudilin, Generating functions of Legendre polynomials,
a tribute to Fred Brafman, preprint, Newcastle

\textbf{22. }J.D.Yu, Notes on Calabi-Yau ordinary differential equations,
arXiv AG/0810.4040.

\textbf{23. }W.Zudilin, More Ramanujan-type formulae for $\dfrac{1}{\pi ^{2}}
$, Russian Math. Surveys, 62 (2007), 634-636.

\bigskip

Institute of Algebraic Meditation

Fogdar\"{o}d 208

24333 H\"{o}\"{o}r, Sweden

gert.almkvist@yahoo.se

\end{document}